\documentclass{amsart}
\usepackage[utf8]{inputenc}

\usepackage{amsmath,amsthm,amssymb,amsrefs}
\usepackage{url,enumerate,epsfig,color,graphicx,epstopdf,mathrsfs}
\usepackage{hyperref}
\usepackage{enumitem}
\usepackage[all]{xy}
\usepackage{subcaption}

\usepackage{listings}
\usepackage{etoolbox}
\usepackage[a4paper, top=3cm, bottom=3cm, left=3cm, right=3cm]{geometry}
\usepackage{tikz}
\usetikzlibrary{matrix,arrows,intersections,tqft,decorations.pathreplacing,cd}

\usepackage[colorinlistoftodos,textsize=tiny,loadshadowlibrary,shadow]{todonotes}
\newcommand{\wt}[1]{\widetilde{#1}}
\newcommand{\wh}[1]{\widehat{#1}}
\newcommand{\ol}[1]{\overline{#1}}
\newcommand{\wH}{\wh{H}}

\usepackage{import}

\newtheorem{theorem}{Theorem}[section]
\newtheorem{proposition}[theorem]{Proposition}
\newtheorem{lemma}[theorem]{Lemma}

\newtheorem{corollary}[theorem]{Corollary}

\theoremstyle{definition}
\newtheorem{definition}[theorem]{Definition}
\newtheorem{example}[theorem]{Example}

\newtheorem{notation}[theorem]{Notation}

\theoremstyle{remark}
\newtheorem*{ack}{Acknowledgments}
\newtheorem*{ai}{AI statement}
\newtheorem{remark}[theorem]{Remark}
\newtheorem{question}[theorem]{Question}

\makeatletter
\DeclareFontEncoding{LS1}{}{}
\DeclareFontSubstitution{LS1}{stix}{m}{n}
\DeclareMathAlphabet{\mathscr}{LS1}{stixscr}{m}{n}
\makeatother
\numberwithin{equation}{section}

\usepackage{etoolbox}

\docsvlist{A,B,C,D,E,F,G,H,I,J,K,L,M,N,O,P,Q,R,S,T,U,V,W,X,Y,Z}

\docsvlist{N,Z,Q,C,R,F,A,W,M,K}

\docsvlist{a,b,e,d,j,m,f,o,q,D,E,I,J,K,P}
\def\mfp{\mathfrak{p}}

\docsvlist{Ann,Ass,ad,adj,Cl,codim,coker,colim,const,Cont,disc,Diff,deg,dist,div,End,FS,Frac,GL,Hom,id,Id,ind,inn,lcm,lk,mind,Nm,Ob,out,Pic,Rad,rb,Res,reg,ri,rk,rot,sign,tb,Tight,tw,trace,Tr}

\docsvlist{M,S}
\newcommand{\adic}[2]{\widehat{#1}_{#2}}
\newcommand{\ZZ}[1]{\adic{\Z}{#1}}
\newcommand{\QQ}[1]{\adic{\Q}{#1}}

\DeclareMathOperator{\im}{Im}

\DeclareMathOperator{\PP}{\mathfrak{P}}
\DeclareMathOperator{\pp}{\mathfrak{p}}

\DeclareMathOperator{\Gal}{Gal}

\DeclareMathOperator{\Bl}{B \ell}

\DeclareMathOperator{\sgn}{sign}

\DeclareMathOperator{\Fitt}{Fitt}

\newcounter{nparcount}
\newcounter{sparcount}

\title{Norm groups of orders and the algebraic unknotting number}
\author{Maciej Borodzik}
\address{Institute of Mathematics\\ University of Warsaw
\\Warsaw, Poland}
\email{mcboro@mimuw.edu.pl}
\author{Wojciech Politarczyk}
\address{Institute of Mathematics\\ University of Warsaw
\\Warsaw, Poland}
\email{w.politarczyk@mimuw.edu.pl}

\date{\today}

\begin{document}

\begin{abstract}
  We use number-theoretic tools to find an explicit and effective obstruction for a knot to have algebraic unknotting number equal to $1$.
  As an example, we show that the connected sum of $T(2,5)$ (the cinquefoil) and $-T(2,3)$ (the negative trefoil) has algebraic
  unknotting number $1$, even though $T(2,5)$ has algebraic unknotting number $2$.
\end{abstract}
\maketitle

\section{Introduction}
\label{sec:introduction}

\subsection{Overview}
Let $K$ be a knot. The unknotting number $u(K)$ is the minimal number of crossing changes needed to turn $K$ into the unknot. Determining
$u(K)$ is one of the hardest problems in knot theory, with only partial results available; see \cite[Problem 1.3]{Kirby}. It was recently shown \cite{BritHerm} that the unknotting number is not additive. It is still conjectured that $u(K_1\# K_2)\ge \max(u(K_1),u(K_2))$, although very little is known in that direction.

The algebraic unknotting number, $u_a(K)$, 
was defined and studied in the 1990s by Fogel, Murakami and Saeki \cite{Fogel_one,Burakami,Saeki}. It is the minimal
number of crossing changes needed to turn the knot $K$ into a knot with Alexander polynomial $1$. While the minimum is taken over all diagrams
of $K$, there is an algebraic interpretation of $u_a(K)$ using Blanchfield pairings; in particular, $u_a(K)$ is determined by
a single algebraic object associated with $K$. We recall the details in Subsection~\ref{sub:blanchfield}. 
The goal of this paper is to show the following statement.
\begin{theorem}\label{thm:alg}
  There is an explicit and effective algorithm determining whether a knot has algebraic unknotting number $1$.
\end{theorem}
We review the algorithm in Subsection~\ref{sub:review} below. The exact algorithm is provided in Subsection~\ref{sub:algorithm}. While the statement might look rather specific to knot theorists, we provide many examples of classes of knots for which the obstruction is computed. Also,
for each particular knot, the obstruction can be effectively computed.

As an exemplary application of Theorem~\ref{thm:alg}, we show how to compute the algebraic unknotting number of a connected sum of knots;
see Section~\ref{sec:connected}.
We obtain the following result.
\begin{theorem}\label{thm:discon}
  There exist knots $K_1$ and $K_2$, with $\max(u_a(K_1),u_a(K_2))=2$ and $u_a(K_1\# K_2)=1$.
\end{theorem}
The result extends the result of Fogel \cite{Fogel_one}, who showed
that the algebraic unknotting number is not additive under connected sums. 

The exact computations are given in Example~\ref{ex:connected-drop}. We can take $K_1=T(2,5)$ and $K_2=-T(2,3)$.

\subsection{Review of the algorithm}\label{sub:review}
The obstruction provides a hierarchy of invariants. Informally, it can be stated as `real algebraic unknotting number', `rational algebraic
unknotting number', `integral algebraic unknotting number'; see Subsections~\ref{prop:qunkn} and~\ref{sub:algorithm}. We sketch the details here.

Suppose $K$ has Nakanishi index $1$, that is, the Alexander module is cyclic \cite{Nakanishi}. We explain that detecting
this is algorithmic in Subsection~\ref{sub:nakanishi}, building on results of Fox--Smythe and Kearton--Wilson. 
Having Nakanishi index equal to $1$ is equivalent to saying that the Alexander module $M(K)$ is isomorphic to $R:=\Z[t^{\pm1}]/(\Delta)$,
where $\Delta$ is the Alexander polynomial. If that is the case, the Blanchfield pairing is determined by a single element $g\in R$
(obtained via $\Bl(1,1)=g/\Delta$, see Subsection~\ref{sub:rank-one-blanchfield}). It is a result going back to Fogel
that $u_a(K)=1$ if and only if $\epsilon g$ is a norm in $R$, that is, if it can be written as $f\ol{f}$. Here, $\epsilon=\pm 1$
can be interpreted as the sign of a crossing change transforming $K$ into a knot with Alexander polynomial $1$.

It can be easily argued that $\epsilon g$ is a norm if and only if its class in $\wH^0(\Z_2;R^\times)$ is zero, where $\wH^*$
is the Tate cohomology; see \cite{EBP_uni1}. We define $\cU(R)=\wH^0(\Z_2;R^\times)$ and attempt to compute it explicitly.

The first check is whether $\epsilon g$ is a norm in $\R\otimes_\Z R$. The algorithm was given in \cite{BF2} and involves jumps of the signature
function. We recover the classical statement that the Levine--Tristram signature bounds the unknotting number. Refer to \cite{BCP1} for
the relation between the Levine--Tristram signature function and the signature jumps. In the present paper, we give another explanation
of the signature jumps. A signature jump gives a local obstruction to being a norm at an infinite place; see Subsection~\ref{sub:infinite}.
Therefore it is a part of a bigger picture, namely of the `rational algebraic unknotting number', that is, the second level of the hierarchy
of invariants.

The second check is whether $\epsilon g$ is a norm in $\Q\otimes_\Z R$. To simplify the discussion here, assume that $\Delta$ is 
irreducible over $\Q$. We handle the general case in Subsection~\ref{sub:reducible-qnorm}. In that case, $E:=\Q\otimes_\Z R$ is a field
with involution. Being a norm in $E$ is fully determined by the Hasse norm principle, that is, $\epsilon g$ is a norm in $E$ if and only if it is a local norm
at each place of $F=E^\tau$. Here, a place is either a finite place (i.e.\ a prime ideal), or an infinite place: the latter correspond to (orbits of) roots of $\Delta$. While checking the infinite places recovers the real part, the finite places provide a computable obstruction.

Finally, suppose $\epsilon g$ is a norm in $\Q\otimes_\Z R$. Then, $\epsilon g$ need not be a norm in $R$. There is an obstruction
for that in a group $\Ob(R)$, which we define (at least in the irreducible case) as $\ker \wH^0(\Z_2;R^\times)\to\wH^0(\Z_2;(\Q\otimes_\Z R)^\times)$.
The key point is that the group $\Ob(R)$ is computable; see Section~\ref{sec:zorm}.

It is helpful to compare this list of invariants with the Lickorish--Jabuka obstruction \cite{Lickorish,Jabuka}. We refer to \cite[Section 4]{BF1}
for a comparison of the Jabuka obstruction with the Lickorish obstruction. We will henceforth refer to that obstruction as the Lickorish criterion. It is expressed
in terms of the linking form; see Theorem~\ref{thm:lickorish}. By carefully analyzing the local norms via the Hasse norm principle, we obtain the following statement. The formal and precise version is given in Theorem~\ref{thm:lick_rat}.
\begin{theorem}\label{thm:lick_intro}
  Suppose $K$ is a knot with Nakanishi index $1$, and real algebraic unknotting number $1$ with sign $-\epsilon$. If it passes the Lickorish criterion with the same sign, then $\epsilon g$
  is a norm in $\Q\otimes_\Z R$.
\end{theorem}
Practically, Theorem~\ref{thm:lick_intro} shows that if a knot has Nakanishi index $1$ and both the signature condition and the Lickorish criterion
are satisfied, then the last remaining obstruction is in $\Ob(R)$. Examples in Subsection~\ref{sec:exampl-where-class} show that this
obstruction can be nontrivial, that is, the group $\Ob(R)$ sees beyond the Lickorish criterion.

\subsection{Examples}\label{sub:examples_intro}
The present article answers a knot-theoretical question using tools from number theory, which might be unfamiliar to knot theorists. This is
why we provide more examples than is usual in such papers.
As a first remark, we note that if $\Delta\in\Z[t^{\pm1}]$ is symmetric, and $\Delta(1)=\pm 1$, then $\Delta$ is the Alexander polynomial of
a knot. Set $R=\Z[t^{\pm1}]/(\Delta)$ and let $(R^\tau)^\times$ be the group of units of the fixed subring $R^\tau$, where $\tau$ is the involution. The realization theorem (see \cite[Theorem 8.3]{Kearton_real}) 
states that for any $g\in (R^\tau)^\times$, there is a knot $K$ such that its Alexander module is $R$,
and the Blanchfield pairing is given by $g$; see Section~\ref{sub:realize}.

That is, many examples can be obtained by studying $\cU(R)$ and $\Ob(R)$ alone.
If $\Ob(R)$ is nontrivial, then there is a knot for which the obstruction is seen. On the other hand, if $\Ob(R)$ is trivial, then for all knots with Alexander module $R$, the Lickorish criterion and the signatures provide the only obstructions. This is the case of the connected sum of $T(2,5)$ and $-T(2,3)$ from Theorem~\ref{thm:discon}. That is, Example~\ref{ex:connected-drop} fits a bigger picture.

This said, we provide two types of examples. First, we provide families or classes of examples in Section~\ref{sec:families}. We first
discuss the case where $R$ is integrally closed; see Subsection~\ref{sub:int_clos}. Many calculations simplify in that case; moreover,
integral closedness is satisfied for many further families. This case is also best understood in the literature: Hermitian
forms over integrally closed rings were studied e.g. in \cite{Conner}, compare also \cite{EBP_uni1}. Our approach draws from \cite{EBP_uni1, EBP_uni2}.

Given this, we pass to our main application, that is, CM fields. This corresponds to polynomials $\Delta$ all of whose roots lie on the unit
circle. In that case, we can often compute $\cU(R)$ explicitly, using the unit signature deficiency $\delta$, counting roughly how many types of signature
jumps can occur on a knot with given Alexander polynomial, compare e.g.\ \cite{EBF_sig} for that question. As an application, we can
compute $\cU(R)$ for cyclotomic polynomials, and then, for Alexander polynomials of torus knots. 

We are also able to show explicit computations for $\Delta$ being a quadratic polynomial. Also, we provide some computations for specific
Alexander polynomials of degree $4$. Section~\ref{sec:families} does not appeal to knot theory; it is purely number-theoretic.

Examples in Section~\ref{sec:examples} are written in the language of knot theory. First, we provide examples of explicit calculations
of some small knots, like the trefoil or the figure-eight knot; see Subsection~\ref{sub:small}. Next, we discuss the knots for which
the Alexander polynomial is the Alexander polynomial of a torus knot. These calculations provide the example of Theorem~\ref{thm:discon}.
Also, we expect that more examples of that type can be found.

Finally, we provide examples explaining that the Lickorish criterion combined with the signature condition is stronger than the rational algebraic unknotting number, but in general not enough to obstruct the (integral) algebraic unknotting number.
Here, we do not just compute $\Ob(R)$, but single out a specific element in $\Ob(R)$ and show that it is nontrivial, while a knot with that element satisfies the Lickorish criterion.

\subsection{Connections to number theory}\label{sub:intro-number-theory}
The problem of detecting if the algebraic unknotting number is $1$
is related to classical problems in number theory. We briefly
sketch the connection.

In the paper, we consider $\Delta\in\Z[t^{\pm1}]$, a symmetric polynomial with $\Delta(1)=\pm 1$. We set
$R=\Z[t^{\pm1}]/(\Delta)$ and $\cA=\Q\otimes_\Z R$. The ring $R$ carries an involution $\tau$ generated by $t\mapsto t^{-1}$.
We write $R^\tau$ for the ring of fixed elements.
For either of these rings, write $\Nm(x)=x\tau(x)$ for the norm of a unit~$x$.
We consider the groups $\cU(R)=(R^\tau)^\times/\Nm(R^\times)$ and $\cU(\cA)=(\cA^\tau)^\times/\Nm(\cA^\times)$.
There is a map
$\iota\colon \cU(R)\to \cU(\cA)$, and we define
\[\Ob(R)=\ker\iota.\]
Then, $\Ob(R)$ consists of classes of units in $R^{\tau}$ that are norms of elements from $\cA$, modulo the norms from~\(R^{\times}\).

Suppose $\Delta$ is irreducible and monic, and $R$ is integrally closed. Then, $E=\cA$ is the field of fractions of $R$, and $R=\cO_E$ is the ring of integers.
Then \cite[Proposition~1]{Lemmermeyer2013Ambiguous} expresses $\Ob(R)$ as the group of ambiguous ideal classes of $E/E^\tau$
modulo the strongly ambiguous classes. We extend this result to the case where $\Delta$ is irreducible but not monic,
see Proposition~\ref{prop:Ob-integrally-closed-case}.

Assume further that $\Delta$ is the cyclotomic polynomial $\Phi_n$. 
Then, $\Ob(R)$ is controlled by signatures of units. Corollary~\ref{cor:Ob-groups-cyclotomic-polynomials} computes the group $\Ob(R)$ in terms of the Hasse unit index and the unit signature deficiency $\delta$ of the real subfield. As $\delta$
is unbounded \cite{DummitKisilevsky2021Unbounded}, the rank of $\Ob(R)$ can be arbitrarily large for cyclotomic $\Delta$. The smallest
$n$, with at least two prime divisors, for which $\Ob$ is non-trivial is $39$.

We emphasize that $\Ob(R)$ measures global obstructions for being a norm, as opposed to local obstructions, governed by the Hasse principle.
In particular, even if $\Delta$ is such that every unit in $R^\tau$ is a norm in $\cA$ (that is, local obstructions vanish, equivalently
$\iota=0$), we may still have $\Ob(R)\neq 0$.

\subsection{Structure of the paper}
Much of the structure of the paper was described in Subsections~\ref{sub:review} and~\ref{sub:examples_intro}. Section~\ref{sec:hermitian}
gives the necessary background on Hermitian forms. We also recall the definition of a Blanchfield pairing and provide precise definitions
of the real, rational and integral algebraic unknotting numbers.

Section~\ref{sec:number_theoretic} recalls some results from algebraic number theory.
We give a definition of a \emph{place} in Subsection~\ref{sub:places}, which is a central notion for the rational obstructions.
We also recall Hilbert reciprocity and the Hasse norm theorem, which are key tools in our further study.

The key mathematical content is in Sections~\ref{sec:qnorm} and~\ref{sec:zorm}, where we discuss obstructions for an element $g$ to
be a norm in $\Q\otimes_\Z R$ and in $R$. These obstructions are applied to various classes of Alexander polynomials in 
Section~\ref{sec:families}.

Section~\ref{sec:knot} returns to knot theory. We provide an algorithm for computing the Nakanishi index and then an algorithm for
computing the algebraic unknotting number. We also study non-additivity of the algebraic unknotting number under connected sums in Subsection~\ref{sec:connected}.

Finally, Section~\ref{sec:examples} contains knot-theoretical examples.

\begin{ai}
    The authors used an AI assistant during the preparation of this paper.
    The main use was proofreading and identifying gaps in the exposition, especially in Lemma~\ref{lem:reduction}.
    The AI assistant helped to develop several examples, including assistance with the PARI/GP library.
    All mathematical content was reviewed, verified and, in many cases, completely rewritten by the authors who take full responsibility for the final text.
\end{ai}
\begin{ack}
  The authors are extremely grateful to Eva Bayer--Fl\"uckiger for many ideas she shared with them and careful revision of the final document.
  MB was supported by the NCN grant OPUS 2024/53/B/ST1/03470.
\end{ack}

\section{Hermitian forms}\label{sec:hermitian}
Throughout this section, rings are commutative and have an identity. An involution on a ring~$R$ is a ring automorphism
$\tau\colon R\to R$ satisfying $\tau^2=\id$. We also write $\ol{x}=\tau(x)$.
Our main examples are Laurent polynomial rings $R=S[t^{\pm1}]$, where $S\subset\C$ is a subring containing~$1$
and preserved by complex conjugation, with involution $\sum a_jt^j\mapsto\sum\ol{a_j}t^{-j}$, and their involution-invariant quotients.

A \emph{Hermitian form} over~$R$ is a pair $(M,h)$, where $M$ is a free $R$-module of finite rank and
$h\colon M\times M\to R$ is linear in the first variable, antilinear in the second, and satisfies
$h(x,y)=\ol{h(y,x)}$. We call $h$ \emph{nonsingular} if the induced map
\[
    h^*\colon M\longrightarrow\ol{\Hom_R(M,R)},\qquad x\longmapsto h(-,x),
\]
is an isomorphism. Here the bar on a module means that we
consider the inverted module structure, i.e. $(r,m)\mapsto \ol{r}m$.

Unless stated otherwise, our Hermitian forms are nonsingular. See~\cite{Scharlau_book} for background.

Our goal is to study Hermitian forms of rank~$1$. We recall Tate cohomology in Subsection~\ref{sub:tate},
classify rank-one forms in Subsection~\ref{sub:dim_one}, and relate them to Blanchfield pairings in Subsection~\ref{sub:blanchfield}.

\subsection{Review of Tate cohomology}\label{sub:tate}
Let $M$ be an abelian group with an automorphism $\tau\colon M\to M$ satisfying $\tau^2=\id$.
Consider the two-periodic cochain complex
\[
    \cdots\longrightarrow M\xrightarrow{1+\tau}\underline{M}\xrightarrow{1-\tau}M\xrightarrow{1+\tau}M\longrightarrow\cdots,
\]
where the underlined copy of~$M$ is in cohomological degree~$0$.
\begin{definition}
    The \emph{Tate cohomology group} $\wH^i(\Z_2;M)$ is the cohomology of this complex in degree~$i$.
\end{definition}
\begin{example}\label{ex:normal}
    If $R$ is a ring with involution, its multiplicative group of units is an abelian group with involution. In multiplicative notation,
    \begin{equation}\label{eq:wh0}
        \wH^0(\Z_2;R^{\times})=\frac{\{x\in R^{\times}\colon\ol{x}=x\}}{\{y\ol{y}\colon y\in R^{\times}\}},\qquad
        \wH^1(\Z_2;R^{\times})=\frac{\{x\in R^{\times}\colon x\ol{x}=1\}}{\{y/\ol{y}\colon y\in R^{\times}\}}.
    \end{equation}
\end{example}
\begin{example}\label{ex:trivial}
    If $\tau\colon M\to M$ is the identity, then $\wH^0(\Z_2;M)=M/2M$ and
    $\wH^1(\Z_2;M)=\ker(M\xrightarrow{\cdot2}M)$.
\end{example}
\begin{lemma}\label{lem:snake}
    A short exact sequence $0\to M'\to M\to M''\to0$ of abelian groups with involutions and equivariant maps induces a long exact sequence
    \[
        \cdots\to\wH^i(\Z_2;M')\to\wH^i(\Z_2;M)\to\wH^i(\Z_2;M'')\to\wH^{i+1}(\Z_2;M')\to\cdots.
    \]
\end{lemma}
The following special case of Shapiro's lemma will be useful.
\begin{lemma}\label{lem:shapiro}
    Suppose $M=M_0\oplus M_0$ and $\tau$ interchanges the summands. Then $\wH^i(\Z_2;M)=0$ for every~$i$.
\end{lemma}
\begin{proof}
    The kernel of $1-\tau$ consists of $(a,a)=(1+\tau)(a,0)$, and the kernel of $1+\tau$
    consists of $(a,-a)=(1-\tau)(a,0)$.
\end{proof}

\begin{lemma}\label{lem:one_half}
    Let $(R,\tau)$ be a ring with involution. The following conditions are equivalent.
    \begin{enumerate}
        \item $\wH^0(\Z_2;R)=0$, where $R$ is regarded as an additive group;
        \item there exists $x\in R$ such that $x+\tau(x)=1$;
        \item every $R$-module~$M$ with an additive involution~$\tau_M$ satisfying
        $\tau_M(rm)=\tau(r)\tau_M(m)$ has $\wH^i(\Z_2;M)=0$ for every~$i$.
    \end{enumerate}
\end{lemma}
\begin{proof}
    Since $\tau(1)=1$, condition~(1) implies that $1$ lies in the image of $1+\tau$, proving~(2).
    Given~(2), if $\tau_M(m)=m$, then
    \[
        (1+\tau_M)(xm)=(x+\tau(x))m=m.
    \]
    If $\tau_M(m)=-m$, then $(1-\tau_M)(xm)=(x+\tau(x))m=m$.
    These identities prove~(3). The implication (3)$\to$(1) is obvious. 
\end{proof}

\begin{lemma}\label{lem:Delta_one_half}
    Suppose $\Delta\in\Lambda=\Z[t^{\pm1}]$ satisfies $\Delta(t)=\Delta(t^{-1})$ and $\Delta(1)$ is odd.
    Give $R=\Lambda/(\Delta)$ the involution induced by $t\mapsto t^{-1}$.
    Then $\wH^i(\Z_2;R)=0$ for every~$i$, with $R$ regarded as an additive group.
\end{lemma}
\begin{proof}
    Write $\Delta=a_0+\sum_{j=1}^d a_j(t^j+t^{-j})$. Since $\Delta(1)$ is odd, $a_0$ is odd.
    The element
    \[
        x=\frac{1-a_0}{2}-\sum_{j=1}^d a_jt^j\in\Lambda
    \]
    satisfies $x+\ol{x}=1-\Delta$, so its image in~$R$ has trace~$1$. Apply Lemma~\ref{lem:one_half}.
\end{proof}

\subsection{Hermitian forms in dimension~1}\label{sub:dim_one}
Suppose $M$ is a free $R$-module of finite rank.
After choosing a basis of~$M$, write a Hermitian form as $h(a,b)=a^TA\ol{b}$,
where $A=\ol{A}^T$. Nonsingularity means that $A$ is invertible. Under a change of basis by~$P\in\GL_n(R)$,
the matrix becomes $P^TA\ol{P}$, so its determinant is multiplied by the norm $\det(P)\ol{\det(P)}$.
Also, $\det A\in(R^\tau)^{\times}$, where $R^\tau$ is the fixed subring. Define
\begin{equation}\label{eq:ur_def}
    \cU(R):=(R^\tau)^{\times}/\Nm(R^{\times}),\qquad\Nm(f)=f\ol{f}.
\end{equation}
We write this quotient group additively and denote the class of a fixed unit~$x$ by~$[x]$.
Thus $[x]=0$ means that $x$ is the norm of a unit.
\begin{definition}
    The \emph{discriminant} of a rank-$n$ Hermitian form is
    \[
        \disc h=\left[(-1)^{n(n-1)/2}\det A\right]\in\cU(R).
    \]
\end{definition}
With this normalization, a rank-one form with coefficient~$s$ has discriminant~$[s]$.
We now specialize to rank-one forms.
\begin{lemma}\label{lem:rank_one}
    Two nonsingular rank-one Hermitian forms $(M_1,h_1)$ and $(M_2,h_2)$ over~$R$ are isometric
    if and only if $\disc h_1=\disc h_2$ in~$\cU(R)$.
\end{lemma}
\begin{proof}
    Choose generators $v_i$ and put $s_i=h_i(v_i,v_i)\in(R^\tau)^{\times}$.
    An isometry $M_2\to M_1$ sends $v_2$ to $fv_1$ for some $f\in R^{\times}$ and satisfies
    $s_2=f\ol{f}s_1$. Conversely, this equality makes that map an isometry.
    It is equivalent to $[s_1]=[s_2]$.
\end{proof}

\subsection{Blanchfield pairings}\label{sub:blanchfield}
Let $K\subset S^3$ be an oriented knot, let $X$ be its exterior, and let $\wt{X}$ be the infinite cyclic cover
associated to the abelianization $\pi_1(X)\to H_1(X;\Z)\cong\Z$, with a positive meridian mapping to~$1$.
The deck transformation corresponding to~$1$ gives $H_1(\wt{X};\Z)$ a module structure over $\Lambda=\Z[t^{\pm1}]$.
\begin{definition}
    The \emph{Alexander module} of~$K$ is
    \[
        M(K)=H_1(\wt{X};\Z)\cong H_1(X;\Lambda).
    \]
    Its minimal number of $\Lambda$-module generators is the \emph{Nakanishi index}, denoted by~$m(K)$.
\end{definition}
The \emph{algebraic unknotting number} $u_a(K)$ is the minimal number of crossing changes needed to turn~$K$
into a knot with Alexander polynomial~$1$. It satisfies
\[
    m(K)\le u_a(K)\le u(K),
\]
where $u(K)$ is the unknotting number; see~\cite{Nakanishi,BF1,BF3}.

The Alexander module carries a Hermitian linking pairing called the \emph{Blanchfield pairing}.
For its construction, see~\cite{Blanchfield,Kearton_Bla,FP_Bla}; see also~\cite{CPT,Conway_exp}.
\begin{theorem}\label{thm:existence_Blanchfield}
    For every knot~$K$, there is a nonsingular Hermitian linking pairing
    \[
        \Bl\colon M(K)\times M(K)\longrightarrow\Q(t)/\Lambda.
    \]
    Nonsingularity means that its adjoint $M(K)\to\ol{\Hom_\Lambda(M(K),\Q(t)/\Lambda)}$ is an isomorphism.
\end{theorem}
We use the convention that $\Bl$ is linear in its first variable. Thus our pairing is the conjugate of the pairing
in~\cite{BF1,BF3}, which is linear in its second variable. A formula in terms of a Seifert matrix is given in~\cite{FP_Bla}.
The following matrix description is proved in~\cite{BF1}.
\begin{proposition}\label{prop:represents}
    For every knot~$K$, there is a Hermitian matrix $A=\ol{A}^T$ over~$\Lambda$, with nonzero determinant, and an isomorphism
    $M(K)\cong\Lambda^n/A\Lambda^n$ under which
    \[
        \Bl(a,b)=a^T\ol{A}^{-1}\ol{b}\pmod\Lambda.
    \]
\end{proposition}
\begin{definition}
    A matrix~$A$ \emph{represents} the Blanchfield pairing if there is an isometry between $\Bl$ and the pairing on
    $\Lambda^n/A\Lambda^n$ in Proposition~\ref{prop:represents}.
\end{definition}
The conjugation on $A^{-1}$ is required by our linearity convention: since $A^T=\ol{A}$, changing~$a$ by~$Ac$
or~$b$ by~$Ad$ changes the displayed expression by an element of~$\Lambda$.

For a subring $S\subset\C$ containing~$1$ closed under conjugation, we 
put $\Lambda_S=S[t,t^{-1}]$ and $\F_S=\Frac(S)$.
Extension of scalars gives $M_S(K)=M(K)\otimes_\Lambda\Lambda_S$ and a pairing~$\Bl_S$ with values in $\F_S(t)/\Lambda_S$.
\begin{definition}
    The invariant $n_S(K)$ is the minimal size of a Hermitian matrix~$A$ over~$\Lambda_S$ representing~$\Bl_S$
    such that $A(1)$ is diagonalizable with diagonal terms~\(\pm 1\) by congruence over~$S$.
\end{definition}
Here diagonalizability means that $P^TA(1)\ol{P}$ is diagonal with diagonal entries~\(\pm1\), for some $P\in\GL_n(S)$.
For~$\R$ this condition is automatic.
Extension of scalars gives
\[
    n_\R(K)\le n_\Q(K)\le n_\Z(K)=u_a(K),
\]
where the last equality is~\cite[Theorem~1.1]{BF3}.
The congruence condition on $A(1)$ in the definition of~$n_\Z(K)$ is part of the characterization of~$u_a(K)$ used in that theorem.

\subsection{Blanchfield pairings of rank~1}\label{sub:rank-one-blanchfield}
Suppose throughout this subsection that $m(K)=1$. Let $\Delta$ be the symmetrized Alexander polynomial, normalized by
$\Delta(t)=\Delta(t^{-1})$ and $\Delta(1)=1$. Then $M(K)\cong R=\Lambda/(\Delta)$.
Fix such an isomorphism, or equivalently a generator of~$M(K)$. We now show a classical way to pass from a pairing
with values in $\Q(t)/\Lambda$ to a pairing with values in $R$.

Note that $g_0=\Bl(1,1)\in\Q(t)/\Lambda$ is annihilated by $\Delta$, hence $g_0=q/\Delta$ for some $q$
defined up to a multiple of $\Delta$.
Hence there is a unique $g\in R$ such that
\[
    \Bl(a,b)=\frac{ag\ol{b}}{\Delta}\pmod\Lambda.
\]
Hermitian symmetry implies $\ol{g}=g$, and nonsingularity implies $g\in R^{\times}$: under the displayed identification,
the adjoint is multiplication by~$g$, up to the conjugation identifying the dual with~$R$.
Consequently,
\[
    \wt{\Bl}\colon R\times R\longrightarrow R,\qquad\wt{\Bl}(a,b)=ag\ol{b},
\]
is a nonsingular rank-one Hermitian form.

Replacing the chosen generator by $f\in R^{\times}$ replaces~$g$ by $f\ol{f}g$, so the class $[g]\in\cU(R)$ is independent of the choice.
The following result is due to Fogel~\cite{Fogel_one}.
\begin{lemma}\label{lem:ub}
    For $\epsilon\in\{\pm1\}$, the following conditions are equivalent.
    \begin{enumerate}
        \item $\Bl$ is represented by the one-dimensional matrix $[\epsilon\Delta]$;
        \item $[\epsilon g]=0$ in~$\cU(R)$, or equivalently $\epsilon g\in\Nm(R^{\times})$.
    \end{enumerate}
\end{lemma}
\begin{proof}
    The matrix $[\epsilon\Delta]$ gives the pairing $\epsilon a\ol{b}/\Delta$ on~$R$.
    An $R$-module automorphism of~$R$ is multiplication by a unit~$f$, and it gives an isometry from~$\Bl$ to this pairing
    exactly when $g=\epsilon f\ol{f}$. Equivalently, $\epsilon g$ is a unit norm; compare Lemma~\ref{lem:rank_one}.
\end{proof}
Recall that a positive crossing change turns a negative crossing into a positive crossing, and a negative crossing change does the reverse.
The following result is given in
\cite[Theorem~1.1]{BF1} and~\cite[Theorem~5.1]{BF3}.
\begin{theorem}\label{thm:fogel}
    Let $K$ have Nakanishi index~$1$, and use $R$ and~$g$ as above. For $\epsilon\in\{\pm1\}$,
    a single crossing change of sign~$-\epsilon$ turns~$K$ into a knot with Alexander polynomial~$1$
    if and only if $[\epsilon g]=0$ in~$\cU(R)$.
\end{theorem}
We write $u_\Z^\epsilon(K)=1$ when these conditions hold, and $u_\Z^\epsilon(K)>1$ otherwise.
Thus the superscript records the sign of the presenting matrix at~$t=1$; the corresponding crossing change has the opposite sign.

\begin{notation}\label{def:signed_S}
    For a conjugation-invariant subring $S\subset\C$ containing~$1$, put
    \[
        R_S=R\otimes_\Z S\cong S[t^{\pm1}]/(\Delta).
    \]
    We write $u_S^\epsilon(K)=1$ if $\Bl_S$ is represented by $[\epsilon\Delta]$, equivalently if the image of~$\epsilon g$
    belongs to $\Nm(R_S^{\times})$, and write $u_S^\epsilon(K)>1$ otherwise.
\end{notation}
The equivalence follows from the proof of Lemma~\ref{lem:ub} after extension of scalars.
We use this signed notation only to distinguish whether a one-dimensional presentation of the specified sign exists.

\subsection{Realization theorem}\label{sub:realize}
If $\Delta\in\Z[t^{\pm1}]$ is a symmetric polynomial with $\Delta(1)=\pm 1$, then $\Delta$ is the Alexander polynomial of a knot;
see \cite[Section 3]{Seifert}.
Generalizing this, one can ask which Hermitian forms can be realized as Blanchfield pairings. We have the following result.

\begin{proposition}\label{prop:trotter}
Let $R=\Z[t^{\pm1}]/(\Delta)$ where $\Delta$ is a symmetric polynomial with $\Delta(1)=\pm 1$. Write $R^\tau$ for the ring of elements of $R$
fixed under the involution. For any $g\in (R^\tau)^\times$, there exists a knot $K$ such that $\Delta(K)=\Delta$, $M(K)\cong R$, and
$\Bl(1,1)=\frac{g}{\Delta}$.
\end{proposition}
The result in the present form was stated in Kearton's paper \cite[Theorem 8.3]{Kearton_real}. It is also a consequence of \cite{Levine_KM1}, which deals exclusively with higher-dimensional links, combined with \cite{Trotter2}, which shows that these results are applicable to classical $1$-dimensional knots.

\section{Number-theoretic preliminaries}\label{sec:number_theoretic}

We fix the notation used below and recall the local facts needed for norm computations.
For background, see Lang~\cite[Chapters I--II]{Lang}, Neukirch~\cite[Chapters I--II]{Neukirch}, and Milne~\cite[Sections 2--3 and 7--8]{MilneANT}.

Let $F$ be a field. A \emph{field extension} is another field $E$ containing $F$ as a subfield; we
denote it by $E/F$. We define the \emph{degree} of the extension $[E:F]$ as the dimension of $E$
as an $F$-vector space.

An example of an extension is obtained by adjoining a root of a polynomial over $F$. If $p(X)\in F[X]$ is irreducible,
then $E:=F[X]/(p)$ is a field, extending $F$. Such $E$ is an example of an \emph{algebraic} extension.
An element $\xi\in E$ is \emph{algebraic} over $F$ if it is a zero of a nontrivial polynomial with coefficients in $F$.

Among all polynomials with that property, we choose those with smallest degree. We call a polynomial of smallest
degree the \emph{minimal polynomial}: it is well defined up to an overall normalization. An extension $E/F$ is \emph{algebraic}
if every element of $E$ is algebraic.

An algebraic extension is \emph{separable} if for any $\xi\in E$, the minimal polynomial $f$ of $\xi$ has no multiple roots in the
algebraic closure of $F$. An algebraic extension is \emph{normal} if for any $\xi\in E$, the minimal polynomial $f$ of $\xi$ factors
as a product of linear terms over $E$. An algebraic extension is \emph{Galois} if it is both separable and normal. In such a case,
we define the group $\Gal(E/F)$ of automorphisms of $E$ fixing $F$. It is called the \emph{Galois group}.

Throughout the paper, unless specified otherwise, all extensions will be finite algebraic. The quadratic extensions $E/F$ studied below
are always Galois.

A \emph{number field} is a finite algebraic extension of $\Q$. Each number field $F$ contains a \emph{ring of integers}, $\cO_F$,
consisting of elements of $F$ for which the minimal polynomial over $\Q$ can be chosen to be monic with integer coefficients.

\subsection{Valuations and places}\label{sub:places}

\begin{definition}
Let $E$ be a number field. An \emph{absolute value} is a function
$|\cdot | \colon E \to \R_{\geq 0}$ such that
\begin{enumerate}
      \item $|x| = 0 \iff x = 0$,
      \item $|x \cdot y| = |x| \cdot |y|$,
      \item $|x+y| \leq |x| + |y|$.
\end{enumerate}
We say that the absolute value~$|\cdot |$ is \emph{non-archimedean} if it satisfies a stronger version of the triangle inequality, namely
\[|x+y| \leq \max(|x|,|y|),\]
and we call $|\cdot|$ \emph{archimedean} otherwise.
\end{definition}
The \emph{trivial} absolute value is given by $|x|=1$ for every $x\ne0$ and $|0|=0$.
An absolute value defines a metric $d(x,y):=|x-y|$ on $E$. We say that two
absolute values are \emph{equivalent} if the corresponding metrics induce the same topologies.
It can be shown, see \cite[Proposition II.3.3]{Neukirch}, that two absolute values $|\cdot|_1$ and $|\cdot|_2$ are equivalent
if and only if there exists $c\in\R_{>0}$ such that $|x|_1=|x|_2^c$ for all $x\in E$.
\begin{definition}
The equivalence class of a nontrivial absolute value is called a \emph{place} on $E$.
\end{definition}
A place is archimedean or non-archimedean depending on whether the corresponding absolute value is.
Archimedean places are often called \emph{infinite}, while non-archimedean \emph{finite}.
The following result is proved in \cite[Proposition II.3.7]{Neukirch}.
\begin{example}\label{eq:q_place}
  For $E=\Q$, each finite place is represented by $|\cdot|_p$ for some prime $p$, where for $a=\frac{b}{c}\in\Q^\times$, with $b,c\in\Z$,
  we write $a=p^m\frac{b'}{c'}$ with $\gcd(p,b'c')=1$, and set $|a|_p=\frac{1}{p^m}$ and $|0|_p=0$.

  The field $\Q$ also has a unique infinite place, corresponding to the standard absolute value on $\Q$.
\end{example}
We denote by $\Omega(E)$ the set of all places on $E$.
Suppose $s\in\Omega(E)$ is a place. We denote by $\wh{E}_s$ the completion
of $E$ with respect to the absolute value representing $s$. The following result is due to Ostrowski, see \cite[Theorem II.4.2]{Neukirch}.
\begin{theorem}\label{thm:ostrowski}
  If $s$ is infinite, its completion $\wh{E}_s$ is either $\R$ or $\C$.
\end{theorem}
Moreover, every embedding of $E$ into $\R$ or $\C$ defines an infinite place. Conversely, every infinite place arises in this way, and two embeddings define the same place if and only if they are equal or complex conjugate. Therefore, infinite places are in bijection
with real embeddings and with pairs of complex conjugate non-real embeddings. Often, we say an infinite place is \emph{real} or \emph{complex} according to whether $\wh{E}_s$ is $\R$
or $\C$.

A notion tightly related to a finite place is that of a valuation.
\begin{definition}
  A \emph{discrete valuation} on $E$, shortly a \emph{valuation},
is a function $\nu\colon E^{\times}\to\R$, extended by $\nu(0)=\infty$, such that:
\begin{enumerate}
      \item $\nu(x \cdot y) = \nu(x) + \nu(y)$,
      \item $\nu(x+y) \geq \min(\nu(x),\nu(y))$,
      \item the image $\nu(E^{\times})$ is a nonzero discrete additive subgroup of~$\R$.
\end{enumerate}
\end{definition}
For each finite place $s$ of $E$, we use throughout the paper the unique associated
\emph{normalized} discrete valuation $\nu_s$, so that $\nu_s(E^\times)=\Z$.
For any absolute value $|\cdot |$ representing $s$, there is $C>0$
such that

\[|x| =
    \begin{cases}
      e^{-C \cdot \nu_{s}(x)}, & x \neq 0, \\
      0, & x = 0.
    \end{cases}
\]
Conversely, any valuation defines a non-archimedean absolute value. Therefore, we sometimes use
a finite place and discrete valuation interchangeably.

Suppose $E$ is a field, and $\nu$ a valuation on $E$. We define
\[\cO_\nu=\{x\in E^{\times}\colon \nu(x)\geq 0\}\cup\{0\}.\]
The ring $\cO_\nu$ is a discrete valuation ring, i.e.\ a principal ideal domain with a
unique nonzero maximal ideal. This maximal ideal is
\[\mfp_{\nu} = \{x \in E^{\times} \colon \nu(x) > 0\} \cup \{0\}.\]
The group of units of $\cO_\nu$ has a simple description in terms of $\nu$:
\[\cO_{\nu}^{\times} = \{x \in \cO_{\nu} \colon \nu(x) = 0\}.\]
\begin{definition}
  The residue field of $\nu$ is
  \[k_{\nu} := \cO_{\nu} / \mfp_{\nu}.\]
\end{definition}
\begin{remark}
  Since we use interchangeably finite places and discrete valuations, we sometimes write $\cO_s$, $\mfp_s$ or $k_s$
  instead of $\cO_\nu$, $\mfp_\nu$ and $k_\nu$.
\end{remark}

Suppose $\nu$ is normalized, that is, $\nu(E^\times)=\Z$. A \emph{uniformizer} is an element $\pi$ of valuation one.

\subsection{Extensions of places}
Suppose that \(F \subset E\) is a separable extension of degree~\(n\).
Pick places \(s_{E} \in \Omega(E)\) and \(s_{F} \in \Omega(F)\).
\begin{definition}
We say that~\(s_{E}\) \emph{divides}~\(s_{F}\) if the restriction of~\(s_{E}\) to~\(F\) is equivalent to~\(s_{F}\).
We denote this fact by writing~\(s_{E} \mid s_{F}\).
If $s_E$ divides $s_F$, we say that $s_E$ lies above $s_F$, and that $s_F$ lies below $s_E$.
\end{definition}

Suppose we are given a pair of places \(s_{E} \mid s_{F}\).
\begin{definition}
  The \emph{local degree} of the pair is the degree of the following field extension
\[d(s_{E} \mid s_{F}) = \left[ \adic{E}{s_{E}} \colon \adic{F}{s_{F}} \right].\]
\end{definition}
Furthermore, if~\(s_{F}\) and~\(s_{E}\) are finite, we can define two more invariants characterizing the pair of places.
Firstly, if \(\nu_{E} \colon E^{\times} \to \R\) is a discrete valuation associated to~\(s_{E}\), then its restriction to~\(F\) is also a discrete valuation.
The images \(\nu_{E}(F^{\times})\) and~\(\nu_{E}(E^{\times})\) are subgroups of~\(\R\) and \(\nu_{E}(F^{\times}) \subset \nu_{E}(E^{\times})\).
\begin{definition}
  We define the \emph{ramification index} of the pair to be the index of the subgroup~\(\nu_{E}(F^{\times})\) in~\(\nu_{E}(E^{\times})\)
\[e(s_{E} \mid s_{F}) = [\nu_{E}(E^{\times}) \colon \nu_{E} (F^{\times})].\]
If $e(s_E|s_F)>1$, we say that $s_E$ \emph{ramifies over $s_F$}. If there is a place $s_E$ over $s_F$ with $e(s_E|s_F)>1$, we say that $s_F$ \emph{ramifies in $E$}.
\end{definition}
An alternative description of ramification is via the ideals. Let $\cO_F$ be the ring of integers of $F$ and $\cO_E$ the ring of integers
of $E$. A finite place on $F$ corresponds to a prime ideal $\mfp$ of $\cO_F$. Consider the ideal $\mfp\cO_E$. It decomposes as a product
$\mfp\cO_E=\mP_1^{e_1}\dots\mP_r^{e_r}$. Then, $\mP_1,\dots,\mP_r$ correspond to places on $E$ over $s_F$, and $e_i$ is the ramification
index of $\mP_i$.

\begin{definition}
Let $s_E$ be a place in $E$ over a finite place $s_F$ of $F$.
Let
\(\cO_{s_{E}} \subset E\) and~\(\cO_{s_{F}} \subset F\) be the respective valuation rings, and let \(\mfp_{s_{E}} \subset \cO_{s_{E}}\), \(\mfp_{s_{F}} \subset \cO_{s_{F}}\) be the corresponding valuation ideals. The \emph{degree of inertia} of the pair $(s_E|s_F)$ is
the degree of the field extension
\[f(s_{E} \mid s_{F}) = \left[ k_{s_E} \colon k_{s_F}\right].\]
We say that $s_F$ is \emph{inert} in $E$ if $f(s_E|s_F)=[E:F]$.
\end{definition}
Finally, we need another piece of terminology.
\begin{definition}
  Let $s_F$ be a place in $F$. We say that $s_F$ \emph{splits in $E$} if there are at least two places on $E$ over $s_F$. If such
  a place is unique, we say that $s_F$ does not split.
\end{definition}

The following result gathers the information about extensions that we need. Its pieces can be found in \cite[Section II.8]{Neukirch},
mostly Lemma~II.8.4 and Proposition~II.8.5 for the first three items, Proposition II.9.1 for the last item.
\begin{theorem}
  Let $E$ be a finite separable extension of a number field $F$. Let $s$ be a place on $F$. Write $n=[E:F]$.
  Then, there exists at least one place over $s$. Let $\Omega_{E/F}$ be the set of places over $s$, then
  \begin{itemize}
    \item With $\wh{F}_s$ denoting the completion of $F$:
      \[E\otimes_F \wh{F}_s=\prod_{s'\in\Omega_{E/F}} \wh{E}_{s'},\]
      where the product is understood in the category of $\wh{F}_s$-algebras;
    \item $\sum_{s'}d(s' \mid s) = n$;
    \item if $s$ is finite, $s'\in\Omega_{E/F}$, $e_{s'}=e(s'|s)$ is the ramification index, and $f_{s'}=f(s'|s)$ 
      the degree of inertia, then $d(s'\mid s)=e_{s'}f_{s'}$;
    \item if $E$ is Galois, then the Galois group acts transitively on $\Omega_{E/F}$, so $e_{s'},f_{s'}$ are constant on $\Omega_{E/F}$. Consequently $e_{s'}f_{s'}m=n$, where $m$ is the cardinality of $\Omega_{E/F}$.
  \end{itemize}
\end{theorem}
As a corollary, we state a result for quadratic extensions, i.e.\ if $[E:F]=2$.
\begin{corollary}
  Suppose $[E:F]=2$ and $E$ is Galois. Let $s$ be a finite place on $F$. There are three possibilities.
\begin{itemize}
  \item $s$ splits, and $\Omega_{E/F}$ contains two elements both with $e_{s'}=f_{s'}=1$;
  \item $s$ is ramified, with $\Omega_{E/F}$ containing a single element with $e_{s'}=2$ and $f_{s'}=1$;
  \item $s$ is inert, with $\Omega_{E/F}=\{s'\}$ and $e_{s'}=1$ and $f_{s'}=2$.
\end{itemize}
\end{corollary}

In Proposition~\ref{prop:nus} we will need a distinction between dyadic places and places of odd residue characteristic.
Suppose $s_F$ is a finite place on a number field $F$. As $F$ is a finite extension of $\Q$, $s_F$ lies above
a finite place in $\Q$. There is a unique prime number $p\in\Z$ such that $s_F$ lies above $(p)$.
\begin{definition}
  We define $p$ to be the \emph{residue characteristic of $s_F$}. We say that $s_F$ is \emph{dyadic} if $p=2$, otherwise $s_F$ is non-dyadic.
\end{definition}

\subsection{Hilbert symbols}\label{sub:hilbert_symbol}
We now recall a few basic facts on Hilbert symbols. References include \cite[Section V.3]{Neukirch}. Our approach
is restricted to quadratic extensions, that is, we set $n=2$; Neukirch's book deals with general $n$.
For $s\in\Omega(E)$, the Hilbert symbol is a symmetric bimultiplicative pairing
\[
  (-,-)_s:\wh{E}_s^\times/(\wh{E}_s^\times)^2
  \times\wh{E}_s^\times/(\wh{E}_s^\times)^2\longrightarrow\{\pm1\}.
\]
For $a,b\in E^\times$, its norm interpretation is
\[
  (a,b)_s=1
  \quad\Longleftrightarrow\quad
  a\in\Nm_{L_s/\wh{E}_s}(L_s^\times),
  \qquad L_s=\wh{E}_s[X]/(X^2-b).
\]
Here $L_s$ is a quadratic field or the split algebra $\wh{E}_s\times\wh{E}_s$.
Equivalently, $(a,b)_s=1$ precisely when $ax^2+by^2=1$ has a solution in $\wh{E}_s$.
At complex places the symbol is always $1$; at real places it is $-1$ exactly when both $a$ and $b$ are negative.

For a finite non-dyadic place, put $r=|k_s|$ and write $a=u\pi^k$, $b=v\pi^\ell$, with $u,v\in\wh{\cO}_s^\times$.
Writing $\left(\frac{u}{\mfp_s}\right)$ for the quadratic character of the residue of $u$ in $k_s^\times$, the explicit formula is
\begin{equation}\label{eq:Hilbert-symbols-local-formula}
    (a,b)_s=(-1)^{k\ell(r-1)/2}
    \left(\frac{u}{\mfp_s}\right)^\ell
    \left(\frac{v}{\mfp_s}\right)^k.
\end{equation}

\subsection{Hilbert reciprocity}\label{sub:hasse_norm}
Let $E/F$ be a quadratic extension of number fields with nontrivial automorphism~$\tau$, and let $s\in\Omega(F)$ be a place of $F$.
Set
\[E_s:=E\otimes_{F}\wh{F}_s,\] 
where $\wh{F}_s$ is the completion.
The involution $\tau$ extends to the map $\tau\colon E_s\to E_s$. 
The following result is a special case of \cite[Proposition II.8.3]{Neukirch}.
\begin{lemma}\label{lem:splitsE_s}
  If $s$ splits, then $E_s=\wh{F}_s\times\wh{F}_s$. If $s$ does not split, $E_s=\wh{E}_{s'}$, where $s'$ is the unique place of $E$ above $s$.
\end{lemma}
We write $\Nm_s\colon E_s^{\times}\to \wh{F}_s^{\times}$, $\Nm_s(x)=x\tau(x)$, and put
\begin{equation}\label{def:us}
  \cU_s = \wh{F}^{\times}_s/\Nm_s(E_s^{\times}).
\end{equation}
\begin{proposition}\label{prop:split}
  The group $\cU_s$ has order $1$ or $2$. If $s$ is split, then $\Nm_s$ is onto.
\end{proposition}
\begin{proof}
  The statement follows from the local reciprocity law \cite[Theorem V.1.3]{Neukirch}, stating that $\cU_s$ is isomorphic to the abelianization of the decomposition group $\Gal(\wh{E}_{s'}/\wh{F}_s)$ of a place $s'$ above $s$. In our situation, by Lemma~\ref{lem:splitsE_s}, this group is $\Z/2$ if $s$ does not split and trivial if $s$ splits.

  For the reader's convenience, we also give an independent proof of the split case.
  Suppose $s$ splits.
  Then, $E_s\cong \wh{F}_s\times\wh{F}_s$.
  The involution $\tau$ permutes the two places above $s$, so $\tau$ swaps the two $\wh{F}_s$ summands of $E_s$. Then, the fixed point set of $\tau$ is the diagonal.
  We have $(x,y)\cdot\tau(x,y)=(xy,xy)$ in the diagonal, that is, $\Nm_s((x,y))=xy$. Any $x'\in\wh{F}_s^{\times}$ is in the image of an element in $E_s^{\times}$ under the norm map: it is enough to take  $(x',1)$.
\end{proof}

The following exact sequence combines the Hasse norm theorem for the quadratic extension~$E/F$ with Hilbert reciprocity.
We refer to \cite[Corollary VI.4.5]{Neukirch} and \cite[Theorem III.4]{Serre_course}.
\begin{theorem}[Hilbert reciprocity]\label{thm:hasse_norm}
    Let $F\subset E$ be a quadratic extension of number fields, $E=F(\alpha)$, where \(\alpha^{2} \in F\).
    The following sequence is exact:
    \[1 \to F^{\times}/\Nm_{E/F}(E^{\times})\to \bigoplus_{s\in\Omega(F)}\cU_s\xrightarrow{\Pi}\{\pm 1\} \to 1.\]
    Here, $\Pi$ is the product map
    \[\Pi(\{a_s\}_{s\in\Omega(F)})=\prod_{s\in\Omega(F)}(a_s,\alpha^2)_s.\]
\end{theorem}

\section{The rational norm problem}\label{sec:qnorm}
We now check whether an element is a norm over the rational numbers, that is, in $\Q\otimes_\Z R$.
Our aim is to express this question through explicit local conditions.
For an irreducible polynomial the criterion is Theorem~\ref{thm:finite}, and Proposition~\ref{prop:reducible-norm-tests} extends it to arbitrary symmetric polynomials.

Until Subsection~\ref{sub:reducible-qnorm}, let \(p\in\Z[t^{\pm1}]\) be a symmetric Laurent polynomial, irreducible over~\(\Q\), with \(p(1)=\pm1\).
Put \(E=\Q[t^{\pm1}]/(p)\), let~\(\tau\) be the involution induced by \(t\mapsto t^{-1}\), and let \(F=E^{\tau}\) be its fixed field. Then~\(E/F\) is quadratic.
We denote by~\(\Nm\) the associated norm map
\[\Nm \colon E^{\times} \to F^{\times}, \quad \Nm(x) = x \tau(x).\]

Put \(R=\Z[t^{\pm1}]/(p)\subset E\), and denote by~\(t_E\) the image of~\(t\) in~\(E\).
Set \(u_F=t_E+t_E^{-1}\in F\).
Then,
\begin{equation}\label{eq:discE}
  E=F[X]/(X^2-u_FX+1),\ \ \disc(X^2-u_FX+1)=(u_F+2)(u_F-2).
\end{equation}

We study rank-one Hermitian forms over~\(E\).
The following special case of Landherr's classification~\cite{landherrAequivalenzHermitescherFormen1935} follows from the Hasse norm principle, Theorem~\ref{thm:hasse_norm}.

\begin{theorem}
    Let~\(h_1\) and~\(h_2\) be rank-one Hermitian forms over~\(E\).
    The following conditions are equivalent:
    \begin{enumerate}
        \item The forms~\(h_1\) and~\(h_2\) are isometric;
        \item \(\disc h_1=\disc h_2\) in~\(F^{\times}/\Nm(E^{\times})\);
        \item For every nonsplit place~\(s\) of~\(F\), finite or infinite,
        \[
            (\disc h_1,u_F^2-4)_s=(\disc h_2,u_F^2-4)_s.
        \]
    \end{enumerate}
\end{theorem}
Here the Hilbert symbols are evaluated on representatives of the discriminant classes; they are independent of those choices because norms have trivial symbol with~\(u_F^2-4\).

For our applications, the most relevant case is that of two rank-one Hermitian forms over~\(R\), whose discriminants lie in \((R^{\tau})^{\times}/\Nm(R^{\times})\).
As a first step towards comparing such forms over~\(R\), we pass to its integral closure in~\(E\).
The majority of this section is devoted to the study of equivalence over~\(E\) for forms whose discriminants have unit representatives in that integral closure.

We conclude this section with a discussion of the case of a reducible~\(p\) in Subsection~\ref{sub:reducible-qnorm}.

\subsection{Normalization of~\(R\)}
\label{sec:normalization-r}\label{sub:tE}

If \(T\) is a finite set of finite places of~\(E\), define the ring of \(T\)-integers by
\[\cO_{E,T}=\bigcap_{s\notin T}\cO_{\nu_s}=\{x\in E\colon \nu_s(x)\ge 0\textrm{ for each finite }s\notin T\}.\]
In particular, if~\(T = \emptyset\), then \(\cO_{E,T} = \cO_{E}\), the ring of algebraic integers of~\(E\).
If \(T\) is invariant under the involution, let \(S\) be the set of finite places of~\(F\) below~\(T\). Then the fixed subring is \(\cO_{E,T}^{\tau}=\cO_{F,S}=\cO_{E,T}\cap F\).

For the ring \(R\), the relevant sets of places are
\begin{equation}\label{eq:def_T}
    T=\{s\in\Omega(E):s\text{ is finite and }\nu_s(t_E)\ne0\},
    \qquad S=\{s\in\Omega(F):s'\mid s\text{ for some }s'\in T\}.
\end{equation}
The following theorem explains why the ring~\(\cO_{E,T}\) appears in the norm problem.
\begin{theorem}\label{thm:normalization}
    The set~\(T\) is finite, and the integral closure of \(R=\Z[t^{\pm1}]/(p)\) in~\(E\) is
    \[
      \wt{R}=\cO_E[t_E,t_E^{-1}],
    \]
    which is equal to $\cO_{E,T}$.
    Moreover, \(\wt{R}\) is a finite \(R\)-module.
\end{theorem}
\begin{proof}
  A non-zero element $t_E$ can have $\nu(t_E)\neq 0$ only for finitely many valuations $\nu$. $T$ is the subset
  of these valuations, hence $T$ is finite.

  A \(\Z\)-basis of~\(\cO_E\) generates~\(\wt{R}=\cO_E[t_E,t_E^{-1}]\) as an \(R\)-module, so \(\wt{R}\) is finite and integral over~\(R\).

 For a finite place~\(s\), localize at its prime ideal~\(\mathfrak P_s\) of~\(\cO_E\). We claim that
 \[
   \wt{R}_{\mathfrak P_s}=\cO_s[t_E,t_E^{-1}]
        =\begin{cases}
            \cO_s,&s\notin T,\\
            E,&s\in T.
        \end{cases}
    \]
    Indeed, if \(\nu_s(t_E)=0\), both generators are already units in~\(\cO_s\).
    Otherwise one generator has negative valuation; adjoining it to the discrete valuation ring~\(\cO_{\nu_s}\) inverts a uniformizer,
    so the localization is~\(E\).
    Now
    \[
      \cO_E[t_E,t_E^{-1}]=\bigcap_s \wt{R}_{\mathfrak P_s}
         =\bigcap_{s\notin T}\cO_s         =\cO_{E,T}.
    \]
    This intersection of valuation rings is integrally closed. Since~\(\wt{R}\) is integral over~\(R\) and has fraction field~\(E\), it is the integral closure of~\(R\).
\end{proof}

The normalization theorem has the following consequences.
\begin{lemma}\label{lem:integral_closure}
    The additive quotient \(\cO_{E,T}/R\) is finite.
    Moreover, \(T=\emptyset\) if and only if, after multiplication by a unit of \(\Z[t^{\pm1}]\), the polynomial~\(p\) is monic.
\end{lemma}
\begin{proof}
    By Theorem~\ref{thm:normalization}, \(\cO_{E,T}\) is a finitely generated $R$-module. 
    Since \(E=\Q\otimes_{\Z}R\), clearing their denominators gives a nonzero integer~\(n\) with \(n\cO_{E,T}\subset R\).
    Thus \(\cO_{E,T}/R\) is a finite module over~\(R/nR\).
    This ring is finite: for each prime~\(\ell\mid n\), the reduction of~\(p\) is nonzero because \(p(1)=\pm1\), so \(R/\ell R\) is finite; the same follows for prime powers by the filtration by powers of~\(\ell\), and then for~\(n\) by the Chinese remainder theorem. This
    proves the first part.

    To prove the second part,
    choose the polynomial associate \(P\in\Z[t]\) of~\(p\) with nonzero constant term and positive leading coefficient.
    It is primitive and irreducible. If \(T=\emptyset\), then~\(t_E\) is an algebraic integer, so its primitive minimal polynomial~\(P\) is monic.
    Conversely, if~\(P\) is monic, the symmetry of~\(p\) implies that its constant term is~\(\pm1\).
    Hence \(t_E\) and~\(t_E^{-1}\) are algebraic integers, and \(T=\emptyset\).
\end{proof}

\begin{remark}\label{rem:normalization-general-factor}
  While Theorem~\ref{thm:normalization} is stated for rings $R=\Z[t^{\pm1}]/(p)$, where $p$ is symmetric, the proof of Theorem~\ref{thm:normalization} does not use the symmetry of~$p$, except for
  the part where $T$ is empty if and only if $p$ is monic.

  If $f\in\Z[t^{\pm1}]$ is primitive and irreducible over~$\Q$, $E_f=\Q[t^{\pm1}]/(f)$, $t_f$ is the image of~$t$, and $T_f$ is the set of finite places~$s$ of~$E_f$ with $\nu_s(t_f)\neq0$, then the integral closure of $R_f=\Z[t^{\pm1}]/(f)$ in~$E_f$ is $\cO_{E_f}[t_f^{\pm1}]=\cO_{E_f,T_f}$, and it is a finite $R_f$-module.
\end{remark}

\begin{lemma}\label{lem:norm_on_T}
    The set~\(T\) is invariant under~\(\tau\), and the involution acts freely on~\(T\). In particular, every place in~\(S\) splits in~\(E\).
\end{lemma}
\begin{proof}
    For \(s\in T\), the conjugate place~\(\tau(s)\) satisfies
    \[
        \nu_{\tau(s)}(t_E)=\nu_s(\tau(t_E))=-\nu_s(t_E)\ne0.
    \]
    Thus \(\tau(s)\in T\) and \(\tau(s)\ne s\), so the place below~\(s\) has two distinct extensions to~\(E\) and hence splits.
\end{proof}

\begin{corollary}\label{cor:T_not_harm}
    Every place~\(s\in S\) has trivial local norm quotient~\(\cU_s\).
\end{corollary}
\begin{proof}
    Every place in~\(S\) splits by Lemma~\ref{lem:norm_on_T}, so Proposition~\ref{prop:split} applies.
\end{proof}
Thus the places below~\(T\) impose no local norm conditions.

Every $\tau$-invariant unit of~\(R\) is a unit of the fixed subring~\(\cO_{F,S}\) of its normalization~\(\cO_{E,T}\).
This motivates the following definition.

\begin{definition}\label{def:T-unimodular}
    A rank-one Hermitian form~\(h\) over~\(E\) is \emph{\(T\)-unimodular} if its discriminant class in \(F^{\times}/\Nm(E^{\times})\) has a representative in~\(\cO_{F,S}^{\times}\).
\end{definition}
In particular, every rank-one Hermitian form over~\(R\) becomes \(T\)-unimodular after extension of scalars to~\(E\).
The following theorem gives a finite set of local tests for isometry of \(T\)-unimodular forms over~\(E\).

\begin{theorem}\label{thm:finite}
    Assume, as in the standing hypotheses of this section, that \(p(1)=\pm1\).
    Let \(T\) and~\(S\) be as above, and let \(h_1,h_2\) be \(T\)-unimodular rank-one Hermitian forms over~\(E\).
    Choose representatives \(q_i\in\cO_{F,S}^{\times}\) of their discriminant classes in \(F^{\times}/\Nm(E^{\times})\).
    For a root~\(\xi\) of~\(p\), write \(q_i(\xi)\) for the image of~\(q_i\) under the embedding \(E\hookrightarrow\C\) sending \(t_E\) to~\(\xi\).
    The following conditions are equivalent:
    \begin{enumerate}
        \item The forms~\(h_1\) and~\(h_2\) are isometric;
        \item \([q_1]=[q_2]\) in \(F^{\times}/\Nm(E^{\times})\);
        \item The following conditions are satisfied:
        \begin{itemize}
            \item for every root~\(\xi\in S^1\) of~\(p\),
            \[q_1(\xi)q_2(\xi)>0;\]
            \item for every finite place~\(s\) of~\(F\) outside~\(S\) lying over a prime divisor of~\(p(-1)\),
            \[(q_1,u_F^2-4)_s=(q_2,u_F^2-4)_s.\]
        \end{itemize}
    \end{enumerate}
\end{theorem}

For discriminant representatives coming from units of~\(R\), a congruence gives a sufficient condition for isometry.
\begin{theorem}\label{thm:finite-square-test}
    Let \(h_1,h_2\) be rank-one Hermitian forms over~\(E\).
    Suppose that \(g_1,g_2\in\Z[t^{\pm1}]\) are symmetric Laurent polynomials whose images \(q_i=g_i(t_E)\) belong to \((R^{\tau})^{\times}\) and represent the discriminant classes of~\(h_i\).
    If
    \begin{enumerate}
        \item \(g_1(\xi)g_2(\xi)>0\) for every root~\(\xi\in S^1\) of~\(p\), and
        \item \(g_1(-1)g_2(-1)\) is a square modulo~\(|p(-1)|\),
    \end{enumerate}
    then \(h_1\) and~\(h_2\) are isometric.
\end{theorem}

We study infinite places in Subsection~\ref{sub:infinite}, finite inert places in Subsection~\ref{sub:inert}, and finite ramified places in Subsection~\ref{sec:ramified-places}.
Theorem~\ref{thm:finite-square-test} is proved in Subsection~\ref{sec:proof-Lickorish-NT}.

\subsection{Norms at infinite places}\label{sub:infinite}

The infinite places of~\(F\) can be read off from the roots of~\(p\).
\begin{lemma}\label{lem:characterizationR}
    The infinite places of~\(F\) correspond to the orbits of roots of~\(p\) under inversion and complex conjugation. 
    They fall into the following three types:
    \begin{enumerate}
        \item Real split places correspond to pairs \(\{\xi,\xi^{-1}\}\) with \(\xi\in\R\setminus\{\pm1\}\).
        \item Real ramified places correspond to pairs \(\{\xi,\overline\xi\}\) with \(\xi\in S^1\setminus\{\pm1\}\).
        \item Complex places, which always split, correspond to four-element orbits with \(\xi\in\C\setminus(S^1\cup\R)\).
    \end{enumerate}
\end{lemma}
\begin{proof}
  It is a classical fact that an Archimedean place on $E$ corresponds to embedding of $E$ either to $\R$ or to $\C$.
    An embedding \(\iota:E\hookrightarrow\C\) is determined by \(\iota(t_E)=\xi\), where~\(\xi\) is a root of~\(p\).
    Its restriction to \(F=\Q(u_F)\) sends \(u_F\) to \(\xi+\xi^{-1}\).
This implies that two embeddings of~\(E\), corresponding to $\xi_1\neq \xi_2$ 
restrict to the same embedding of~\(F\) if and only if $\xi_1=\xi_2^{-1}$.

Take an infinite place $s'$ of $E$.
 Let~\(s\) be the place of~\(F\) below $s'$. In the three cases, respectively, the pair \((\wh F_s,E\otimes_F\wh F_s)\) is
    \[
        (\R,\R\times\R),\qquad (\R,\C),\qquad (\C,\C\times\C).
    \]
    Thus the first and third cases split, and the second ramifies.
\end{proof}

Among infinite places, only the real ramified places can obstruct an element of~\(F^{\times}\) from being a norm.
\begin{lemma}\label{lem:realxi}
    Let \(g\in F^{\times}\), let \(\xi\in S^1\) be a root of~\(p\), and let~\(s\) be the real place of~\(F\) obtained by restricting the embedding \(t_E\mapsto\xi\).
    Then \(g\) is a local norm at~\(s\) if and only if \(g(\xi)>0\).
\end{lemma}
\begin{proof}
    Here \(E\otimes_F\wh F_s\cong\C\) and \(\wh F_s\cong\R\).
    The local norm is \(\C^{\times}\to\R^{\times}\), \(z\mapsto z\overline z=|z|^2\), whose image is \(\R_{>0}\).
\end{proof}

Consequently, we obtain the following criterion.
\begin{proposition}\label{prop:real_places}
    An element \(g\in F^{\times}\) is a local norm at every infinite place of~\(F\) if and only if \(g(\xi)>0\) for every root~\(\xi\in S^1\) of~\(p\).
\end{proposition}
\begin{proof}
    At split places the norm map is surjective by Proposition~\ref{prop:split}; at real ramified places apply Lemma~\ref{lem:realxi}.
\end{proof}

\begin{remark}
  The distinction between the three cases of Lemma~\ref{lem:characterizationR} is parallel to the distinction
  of basic linking forms in \cite[Section 2]{BCP1}. The ramifying roots of Lemma~\ref{lem:characterizationR} correspond precisely
  to $\Xi^{\C}_{\mathfrak{e}}$ of \cite{BCP1}. The condition $g(\xi)>0$ controls the signature jumps at $\xi$.
\end{remark}

\subsection{Inert places}\label{sub:inert}
We will now show that inert places do not provide any restriction on being a norm. The statement is a special case of the following result; see~\cite[Section V.2, Proposition 2]{Serre_local}.
\begin{theorem}\label{thm:serre_local}
    Let \(L/K\) be a finite extension of number fields, let~\(s\) be a finite place of~\(K\), and let~\(s'\mid s\) be a place of~\(L\).
    If the extension of completions \(\wh L_{s'}/\wh K_s\) is unramified, then
    \[
        \Nm_{\wh L_{s'}/\wh K_s}:\wh{\cO}_{L,s'}^{\times}\longrightarrow\wh{\cO}_{K,s}^{\times}
    \]
    is surjective, where the two rings are the valuation rings of the completions.
\end{theorem}

Recall that \(T=\{s:\nu_s(t_E)\ne0\}\) and that~\(S\) is its image in~\(\Omega(F)\).
\begin{lemma}\label{lem:inert1}
    If a finite place~\(s\) of~\(F\) is inert in~\(E\), its unique extension~\(s'\) does not belong to~\(T\); in particular, \(s\notin S\).
\end{lemma}
\begin{proof}
    Every place in~\(S\) splits in~\(E\), by Lemma~\ref{lem:norm_on_T}.
\end{proof}

\begin{proposition}\label{prop:inert}
    Let \(g\in\cO_{F,S}^{\times}\). Then \(g\) is a local norm at every unramified finite place~\(s\) of~\(F\).
    In particular, this applies to every \(g\in(R^{\tau})^{\times}\), where \(R=\Z[t^{\pm1}]/(p)\).
\end{proposition}
\begin{proof}
    If~\(s\) splits, the norm map is surjective by Proposition~\ref{prop:split}.
    Otherwise~\(s\) is inert, so \(s\notin S\) by Lemma~\ref{lem:inert1}. Thus \(\nu_s(g)=0\), and Theorem~\ref{thm:serre_local} shows that~\(g\) is a local norm.
    For the last assertion, \(R\subset\cO_{E,T}\) because \(t_E\) and~\(t_E^{-1}\) are integral at every finite place outside~\(T\).
    Hence an invariant unit of~\(R\) belongs to \(\cO_{E,T}^{\times}\cap F=\cO_{F,S}^{\times}\).
\end{proof}

\subsection{Ramified places}
\label{sec:ramified-places}

Unlike split and inert places, finite ramified places provide potential obstructions to being a norm. Therefore our goal is
to locate the finite places which ramify in~\(E/F\).
We use the relative field discriminant~\(\md_{E/F}\), whose localizations are generated by the trace determinants
\[
    \det\left(\Tr_{E/F}(b_i b_j)\right)
\]
for integral bases over the corresponding valuation rings.
It is the norm of the relative different~\(\mD_{E/F}\); see~\cite[Theorem III.2.9]{Neukirch}.
\begin{equation}\label{eq:discriminant}
    \md_{E/F}=\Nm_{E/F}(\mD_{E/F}).
\end{equation}
The following ramification criterion is due to Dedekind; see~\cite[Corollary III.2.12]{Neukirch}.
\begin{theorem}\label{thm:discrim}
    A finite place~\(s\) of~\(F\) ramifies in~\(E\) if and only if its prime ideal~\(\mfp_s\) divides~\(\md_{E/F}\).
\end{theorem}

The field discriminant must be distinguished from the polynomial discriminant
\(\disc(X^2-u_FX+1)=u_F^2-4\).
\begin{lemma}\label{lem:invert}
    Let~\(s\) be a finite place of~\(F\), and let \(f\in\cO_{\nu_s}[X]\) be monic and irreducible over~\(F\).
    If \(\disc(f)\in\cO_{\nu_s}^{\times}\), then~\(s\) is unramified in \(L=F[X]/(f)\).
\end{lemma}
\begin{proof}
  By \cite[Theorem III.2.5]{Neukirch}, the different $\mD_{L/F}$ contains $f'(x)$, and so
  $\md_{L/F}=\Nm_{L/F}\mD_{L/F}$ divides $\Nm_{L/F}(f'(x))$. The latter is the discriminant up to sign.
  We conclude by Theorem~\ref{thm:discrim}.
\end{proof}

Lemma~\ref{lem:invert} gives the following valuation criterion for ramification.
\begin{proposition}\label{prop:nus}
    Let~\(s\) be a finite place of~\(F\), put \(d=u_F^2-4\), and let~\(\ell\) be the residue characteristic of~\(s\).
    \begin{enumerate}
        \item If \(\nu_s(d)\) is odd, then~\(s\) ramifies in~\(E\).
        \item If \(\nu_s(d)\) is even and \(\ell\ne2\), then~\(s\) is unramified in~\(E\).
        \item If \(u_F\in\cO_{\nu_s}\) and \(\nu_s(d)=0\), then~\(s\) is unramified in~\(E\).
    \end{enumerate}
\end{proposition}
\begin{proof}
    Recall that \(E=F(\sqrt d)\).

    For~(1), let~\(s'\mid s\). The normalized valuations satisfy
    \[
        e(s'\mid s)\nu_s(d)=\nu_{s'}(d)=2\nu_{s'}(\sqrt d).
    \]
    If \(\nu_s(d)\) is odd, then \(e(s'\mid s)\) must be even. Since~\(E/F\) is quadratic, this forces \(e(s'\mid s)=2\).

    For~(2), write \(\nu_s(d)=2m\), choose a uniformizer \(\pi\in F\) at~\(s\), and put \(d'=d\pi^{-2m}\).
    Then \(d'\in\cO_{\nu_s}^{\times}\) and \(E=F(\sqrt{d'})\).
    The polynomial \(X^2-d'\) is monic with coefficients in~\(\cO_{\nu_s}\), and its discriminant~\(4d'\) is a unit in that ring because \(\ell\ne2\).
    Lemma~\ref{lem:invert} therefore shows that~\(s\) is unramified.

    For~(3), the polynomial \(X^2-u_FX+1\) is monic with coefficients in~\(\cO_{\nu_s}\) and defines~\(E/F\).
    Its discriminant is \(u_F^2-4=d\), which is a unit in~\(\cO_{\nu_s}\) by hypothesis.
    Apply Lemma~\ref{lem:invert} again.
\end{proof}
\begin{remark}
  It can be checked that for the local field $\QQ{2}$, the extension $\QQ{2}(\sqrt{3})$ ramifies over $\QQ{2}$
  (meaning that the unique prime ideal in $\ZZ{2}$ ramifies). That is, the distinction between items (2) and (3)
  of Proposition~\ref{prop:nus} is necessary.
\end{remark}

We now use \(p(1)=\pm1\) to restrict further the number of places that ramify.
In the factorization \(u_F^2-4=(u_F-2)(u_F+2)\), the first factor is a unit at every finite place outside~\(S\).
\begin{lemma}\label{lem:ue1}
    Let~\(s'\notin T\) be a finite place of~\(E\), and let~\(s\) be the place of~\(F\) below it.
    Then \(\nu_{s'}(t_E-1)=0\) and \(\nu_s(u_F-2)=0\).
\end{lemma}
\begin{proof}
    Write \(p(t)-p(1)=(t-1)h(t)\), where \(h\in\Z[t^{\pm1}]\).
    Evaluating at~\(t_E\) gives
    \[
        (t_E-1)h(t_E)=-p(1)=\pm1.
    \]
    Since \(t_E\) and~\(t_E^{-1}\) are units at~\(s'\), both factors on the left are integral there. Hence \(\nu_{s'}(t_E-1)=0\).
    The identity \(u_F-2=t_E^{-1}(t_E-1)^2\) now gives
    \[
        e(s'\mid s)\nu_s(u_F-2)=\nu_{s'}(u_F-2)=0,
    \]
    so \(\nu_s(u_F-2)=0\).
\end{proof}

Next, we consider \(u_F+2\).
\begin{lemma}\label{lem:ue2}
    Let~\(s\notin S\) be a finite place of~\(F\) which ramifies in~\(E\).
    Then
    \[
        0<\nu_s(u_F+2)\leq\nu_s(p(-1)),
    \]
    and \(\nu_s(u_F+2)\) is odd. In particular, \(s\) lies above an odd prime divisor of~\(p(-1)\).
\end{lemma}
\begin{proof}
    Since~\(s\notin S\), every place~\(s'\mid s\) lies outside~\(T\), so \(t_E\) and~\(t_E^{-1}\) are units at~\(s'\).
    Thus \(u_F=t_E+t_E^{-1}\in\cO_{\nu_s}\), and Lemma~\ref{lem:ue1} gives \(\nu_s(u_F-2)=0\).
    If \(\nu_s(u_F+2)=0\), then \(u_F^2-4\) is a unit at~\(s\), contradicting ramification by Lemma~\ref{lem:invert}.
    Hence \(\nu_s(u_F+2)>0\).

    Since~\(p\) is a symmetric Laurent polynomial, write \(p(t)=f(t+t^{-1})\) with \(f\in\Z[X]\).
    We have \(f(u_F)=0\) and \(f(-2)=p(-1)\).
    Writing \(f(X)-f(-2)=(X+2)H(X)\), with \(H\in\Z[X]\), and evaluating at~\(u_F\) yields
    \[
        (u_F+2)H(u_F)=-p(-1).
    \]
    Since \(H(u_F)\in\cO_{\nu_s}\), this proves \(\nu_s(u_F+2)\leq\nu_s(p(-1))\).
    Consequently, the residue characteristic of~\(s\) divides~\(p(-1)\), which is odd because \(p(-1)\equiv p(1)=\pm1\pmod2\).
    Finally,
    \[
        \nu_s(u_F^2-4)=\nu_s(u_F+2).
    \]
    This valuation must be odd, since an even valuation at a non-dyadic place would imply that~\(s\) is unramified by Proposition~\ref{prop:nus}(2).
\end{proof}

\begin{proof}[Proof of Theorem~\ref{thm:finite}]
    Put \(g=q_1/q_2\in\cO_{F,S}^{\times}\).
    The forms are isometric precisely when \(g\in\Nm(E^{\times})\), and the Hasse norm principle reduces this to the local norm conditions.
    At infinite places, these conditions are exactly \(q_1(\xi)q_2(\xi)>0\) at the roots on~\(S^1\), by Proposition~\ref{prop:real_places}.
    Every place in~\(S\) splits, so it imposes no condition.
    At a finite place outside~\(S\), the element~\(g\) is a unit. If the place is inert, it is a local norm by Theorem~\ref{thm:serre_local}; if it splits, it is a local norm by Proposition~\ref{prop:split}.
    Every remaining finite place is ramified and lies over a prime divisor of~\(p(-1)\), by Lemma~\ref{lem:ue2}.
    At such a place~\(s\), the local norm condition is
    \[
        (g,u_F^2-4)_s=1,
    \]
    which is equivalent to the equality of the two Hilbert symbols in the statement, by multiplicativity in the first argument.
    Thus the stated conditions hold if and only if~\(g\) is a norm at every place.
    By Theorem~\ref{thm:hasse_norm}, this is equivalent to \(g\in\Nm(E^{\times})\), and hence to the two forms being isometric.
\end{proof}

\subsection{Proof of Theorem~\ref{thm:finite-square-test}}
\label{sec:proof-Lickorish-NT}

Setting  \(t=-1\) gives a necessary condition for a unit to be a norm in~\(R\):
\[
    R/(t_E+1)R\cong\Z/|p(-1)|\Z.
\]
The involution is trivial on this quotient, so the image of \(x\tau(x)\) is the square of the image of~\(x\).
The following proposition shows that this congruence condition implies all the finite local norm conditions.

\begin{proposition}\label{prop:lickorish_NT}
    Let \(b\in(R^{\tau})^{\times}\).
    If its image in \(R/(t_E+1)R\) is a square, then~\(b\) is a local norm at every finite place of~\(F\).
    If, in addition, \(b\) is a local norm at every infinite place, then \(b\in\Nm_{E/F}(E^{\times})\).
\end{proposition}
\begin{proof}
    Since \(b\in\cO_{F,S}^{\times}\), Proposition~\ref{prop:inert} shows that $b$ is a local norm at
    all inert and split places.

    Let~\(s\) be a ramified finite place of~\(F\), and let~\(s'\) be the unique place of~\(E\) above it.
    The places in~\(S\) split, so \(s\notin S\) and \(s'\notin T\).
    By Lemma~\ref{lem:ue2}, the residue characteristic~\(\ell\) is odd and \(\nu_s(u_F+2)>0\).
    Since \(u_F+2=t_E^{-1}(t_E+1)^2\) and \(\nu_{s'}(t_E)=0\), we have \(\nu_{s'}(t_E+1)>0\).

    The inclusion \(R\subset\cO_{s'}\) therefore induces a residue map which factors as
    \[
        R\longrightarrow R/(t_E+1)R
        \cong\Z/|p(-1)|\Z\longrightarrow\F_{\ell}\hookrightarrow k_{s'}=k_s.
    \]
    Here \(k_{s'}=k_s\) because a ramified quadratic extension has residue degree one.
    This map sends the square specialization of~\(b\) to a square in~\(k_s\).
    Its residue is nonzero because~\(b\) is a unit at~\(s\).
    Thus \(X^2-b\) has a nonzero root in~\(k_s\), and its derivative at that root is nonzero since \(\ell\ne2\).
    Hensel's lemma~\cite[Lemma II.4.6]{Neukirch} gives \(d\in\wh F_s^{\times}\) with \(d^2=b\).
    Viewing~\(d\) as an element of~\(\wh E_{s'}\), we obtain
    \[
        b=d^2=\Nm_{\wh E_{s'}/\wh F_s}(d).
    \]
    This proves the finite local norm conditions. As we have assumed that $b$ is a local norm at infinite
    places, Theorem~\ref{thm:hasse_norm} gives \(b\in\Nm_{E/F}(E^{\times})\).
\end{proof}
\begin{corollary}\label{cor:lickorish-specialization}
    Let \(b\in(R^\tau)^\times\), and let \(g\in\Z[t^{\pm1}]\) represent~\(b\).
    The image of~\(b\) in \(R/(t_E+1)R\) is a square if and only if \(g(-1)\) is a square modulo~\(|p(-1)|\).
\end{corollary}
\begin{proof}
    Evaluation at~\(t=-1\) induces the isomorphism
    \[
        R/(t_E+1)R\xrightarrow{\ \cong\ }\Z/|p(-1)|\Z,
    \]
    and sends the image of~\(b\) to the residue class of~\(g(-1)\).
\end{proof}
We are now in position to prove Theorem~\ref{thm:finite-square-test}.
\begin{proof}[Proof of Theorem~\ref{thm:finite-square-test}]
    Put \(b=q_1q_2\in(R^\tau)^\times\).
    Since \(q_2^2=\Nm(q_2)\), the elements \(b\) and \(q_1/q_2\) have the same class in \(\cU(E)\).
   By Theorem~\ref{thm:finite}, the forms \(h_1\) and~\(h_2\) are isometric if and only if \(b\in\Nm(E^\times)\).
    By Proposition~\ref{prop:real_places}, condition~(1) says that \(b\) is a local norm at every infinite place of~\(F\).
    Condition~(2) says that the image \(g_1(-1)g_2(-1)\) of~\(b\) in \(R/(t_E+1)R\cong\Z/|p(-1)|\Z\) is a square.
    Proposition~\ref{prop:lickorish_NT} therefore gives \(b\in\Nm(E^\times)\).
\end{proof}

\subsection{Reducible polynomials}\label{sub:reducible-qnorm}
We now allow~\(p\in\Z[t^{\pm1}]\) to be any nonconstant symmetric Laurent polynomial with \(p(1)=\pm1\).
Put \(R=\Z[t^{\pm1}]/(p)\) and \(\cA=\Q[t^{\pm1}]/(p)\), with the involution~\(\tau\) induced by \(t\mapsto t^{-1}\), and write~\(t_R\) for the image of~\(t\) in~\(R\).
Let~\(p_{\mathrm{red}}\) be the product of the distinct irreducible factors of~\(p\), and set \(\cA_{\mathrm{red}}=\Q[t^{\pm1}]/(p_{\mathrm{red}})\). We have
\[
    p_{\mathrm{red}}=\prod_{i=1}^{r}p_i\prod_{j=1}^{m}q_j\tau(q_j),
\]
where the~\(p_i\) are symmetric and \(q_j\) is coprime with~\(\tau(q_j)\).
We normalize $p_i$ so that $p_i(1)=1$. 
Write \(E_i=\Q[t^{\pm1}]/(p_i)\), \(F_i=E_i^{\tau}\), and \(E'_j=\Q[t^{\pm1}]/(q_j)\).
Since \(p(1)\ne0\) and \(p(-1)\) is odd, neither~\(1\) nor~\(-1\) is a root, so each~\(E_i/F_i\) is quadratic.
Clearly,
\begin{equation}\label{eq:field_of_fractions}
    \cA_{\mathrm{red}}\cong\prod_{i=1}^{r}E_i\times\prod_{j=1}^{m}(E'_j\times E'_j).
\end{equation}
Here~\(\tau\) acts on each~\(E_i\) as before and exchanges the two factors in each pair, after identifying the factor belonging to~\(\tau(q_j)\) with~\(E'_j\).
In particular, \(\cA_{\mathrm{red}}^{\tau}\cong\prod_iF_i\times\prod_jE'_j\).

\begin{lemma}\label{lem:qnorm}
  Projection $\cA\to \prod E_i$ induces an isomorphism
    \[
        (\cA^{\tau})^{\times}/\Nm(\cA^{\times})
        \cong\prod_{i=1}^{r}F_i^{\times}/\Nm_{E_i/F_i}(E_i^{\times}).
    \]
    Thus a unit of~\(\cA^\tau\) is a norm if and only if its image in every~\(F_i\) is a norm from~\(E_i\).
\end{lemma}
\begin{proof}
  First we pass to the reduced part of $\cA$.
  Let $I$ be the kernel of $\pi\colon \cA\to\cA_{\mathrm{red}}$.
  It is clearly a nilpotent ideal, and $\tau I=I$.
  It is well-known that if $x$ is nilpotent, then $1+x$
  is a unit. This observation leads to
  the short exact sequence
  \begin{equation}\label{eq:nilpotent}
  1\to 1+I\to\cA^\times\to \cA_{\mathrm{red}}^\times\to 1.\end{equation}
  We claim that $\wH^*(\Z_2;1+I)=0$.

  To prove the claim notice that the map $x\mapsto x^2$ is an automorphism of $1+I$. This follows from Hensel's lemma, but can be proved using elementary means, by considering the expansion of $(1+y)^{1/2}$ into power series: the expansion is finite if $y\in I$.
  The squaring map is a multiplication by $2$ on $\wH^*(\Z_2;1+I)$, which is an automorphism (because it is an automorphism of $1+I$) and the zero map (as Tate cohomology groups are annihilated by~$2$), hence
  $\wH^*(\Z_2;1+I)=0$ proving the claim.

  From Lemma~\ref{lem:snake} applied to \eqref{eq:nilpotent},
  we conclude that 
  the map \(\pi:\cA\to\cA_{\mathrm{red}}\) 
    induces an isomorphism
    \[
        (\cA^{\tau})^{\times}/\Nm(\cA^{\times})
        \longrightarrow
        (\cA_{\mathrm{red}}^{\tau})^{\times}/\Nm(\cA_{\mathrm{red}}^{\times}).
    \]
    Now by \eqref{eq:field_of_fractions},
    $\wH^0(\Z_2;\cA_{\mathrm{red}}^\times)$ 
    is the product of the groups corresponding to the factors. Now
    by Lemma~\ref{lem:shapiro}, the primed factors do
    not contribute anything.
\end{proof}

\begin{proposition}\label{prop:reducible-norm-tests}
    Let \(b\in(R^{\tau})^{\times}\).
    For each symmetric factor~\(p_i\), write \(b_i\) for the image of~\(b\) in \(R_i=\Z[t^{\pm1}]/(p_i)\), let~\(t_{E_i}\) be the image of~\(t\), and put \(u_i=t_{E_i}+t_{E_i}^{-1}\in F_i\).
    For a root~\(\xi\) of~\(p\), let~\(b(\xi)\) denote the image of~\(b\) under the evaluation map \(R\to\C\), \(t_R\mapsto\xi\).
    \begin{enumerate}
        \item The element~\(b\) is a norm in~\(\cA\) if and only if
        \[
            b(\xi)>0\quad\text{for every root }\xi\in S^1\text{ of }p,
        \]
        and, for every~\(i\) and every finite place~\(s\) of~\(F_i\) ramified in~\(E_i\),
        \[
            (b_i,u_i^2-4)_s=1.
        \]
        \item If \(b(\xi)>0\) at every root~\(\xi\in S^1\) of~\(p\), and the image of~\(b\) in
        \[
            R/(t_R+1)R\cong\Z/|p(-1)|\Z
        \]
        is a square, then~\(b\) is a norm in~\(\cA\).
    \end{enumerate}
\end{proposition}
\begin{proof}
    Each~\(b_i\) is an invariant unit of~\(R_i\).
    Lemma~\ref{lem:qnorm} reduces the norm problem to the extensions~\(E_i/F_i\).
    Proposition~\ref{prop:inert} handles all unramified finite places in each component, and the remaining finite local norm conditions are precisely the displayed Hilbert-symbol conditions.
    At infinite places, Proposition~\ref{prop:real_places} gives positivity at the roots on~\(S^1\).
    Every such root belongs to a symmetric factor: its inverse is its complex conjugate, so its irreducible factor is invariant under~\(\tau\).
    The Hasse norm principle, Theorem~\ref{thm:hasse_norm}, now proves~(1).

    For~(2), the quotient maps \(R\to R_i\) induce maps
    \[
        R/(t_R+1)R\longrightarrow R_i/(t_{E_i}+1)R_i,
    \]
    which carry the square image of~\(b\) to the square image of~\(b_i\).
    Proposition~\ref{prop:lickorish_NT} gives all finite local norm conditions in each component, so~(1) applies.
\end{proof}

\section{The integral norm problem}\label{sec:zorm}
We now pass to computing obstructions to being a norm in $R$. We define the group $\Ob(R)$ measuring the obstruction. Then, we
provide results that allow for computing that group. These tools are equation \eqref{eq:short} and Lemma~\ref{lem:exact_K0}. While they do not
give a direct answer, the two results are enough to reduce the computation of~$\Ob(R)$ to more classical objects.

We keep the setup of Subsection~\ref{sub:reducible-qnorm}: let \(p\in\Z[t^{\pm1}]\) be a nonconstant symmetric Laurent polynomial with \(p(1)=\pm1\), possibly with repeated factors, and put
\[
    R=\Z[t^{\pm1}]/(p),\qquad \cA=\Q\otimes_{\Z}R=\Q[t^{\pm1}]/(p).
\]
The involution~\(\tau\) is induced by \(t\mapsto t^{-1}\), and \(\Nm(x)=x\tau(x)\).
For \(g\in(R^\tau)^\times\), Section~\ref{sec:qnorm} determines whether \(g\in\Nm(\cA^\times)\); we now ask whether \(g\in\Nm(R^\times)\).

We first remove repeated factors of~\(p\), then study the obstruction to a rational norm being the norm of a unit of~\(R\).
\subsection{Reduction by the radical}
Let~\(p_{\mathrm{red}}\) be the product of the distinct irreducible factors of~\(p\), as in Subsection~\ref{sub:reducible-qnorm}, and put
\[
    R_{\mathrm{red}}=\Z[t^{\pm1}]/(p_{\mathrm{red}}).
\]
There is a quotient map \(R\to R_{\mathrm{red}}\).

\begin{lemma}\label{lem:reduction}
	The quotient map induces an isomorphism $\cU(R)\cong \cU(R_{\mathrm{red}})$.
\end{lemma}
\begin{proof}
    Let \(\ma=(p_{\mathrm{red}})\subset R\). This ideal is invariant under~\(\tau\), satisfies \(R/\ma\cong R_{\mathrm{red}}\), and is nilpotent: if~\(n\) is at least the largest multiplicity of a factor of~\(p\), then \(\ma^n=0\).
    Units lift modulo a nilpotent ideal, as in the proof of Lemma~\ref{lem:qnorm}, so there is a short exact sequence of groups with involution
    \[
        1\longrightarrow1+\ma\longrightarrow R^\times
        \longrightarrow R_{\mathrm{red}}^\times\longrightarrow1.
    \]
    By Lemma~\ref{lem:snake}, it suffices to prove that \(\wH^*(\Z_2;1+\ma)=0\).

    For \(k\geq1\), the map \(1+x\mapsto x\bmod\ma^{k+1}\) gives an equivariant exact sequence
    \[
        1\longrightarrow1+\ma^{k+1}\longrightarrow1+\ma^k
        \longrightarrow\ma^k/\ma^{k+1}\longrightarrow0,
    \]
    whose last term is an additive \(R\)-module with semilinear involution.
    Since \(p(1)=\pm1\), Lemmas~\ref{lem:Delta_one_half} and~\ref{lem:one_half} give
    \(\wH^*(\Z_2;\ma^k/\ma^{k+1})=0\).
    Starting from \(1+\ma^n=1\), downward induction using Lemma~\ref{lem:snake} proves the required vanishing.
\end{proof}
The map $\cU(R)\to\cU(R_{\mathrm{red}})$ sends the class of $g\in(R^\tau)^{\times}$ to the class of its image in $(R_{\mathrm{red}}^\tau)^{\times}$. In particular, we have the following corollary.
\begin{corollary}\label{cor:is_norm}
  The element $g$ is a norm in $R$ if and only if its image is a norm in $R_{\mathrm{red}}$.
\end{corollary}

\begin{remark}\label{rem:redQ}
    The proof of Lemma~\ref{lem:qnorm} establishes the rational counterpart
    \(\cU(\cA)\cong\cU(\cA_{\mathrm{red}})\).
    Together with Lemma~\ref{lem:reduction}, this gives the commutative diagram
    \[
        \begin{tikzcd}
            \cU(R) \ar[r] \ar[d,"\cong"'] & \cU(\cA) \ar[d,"\cong"] \\
            \cU(R_{\mathrm{red}}) \ar[r] & \cU(\cA_{\mathrm{red}}).
        \end{tikzcd}
    \]
    The horizontal maps are induced by inclusion, and the vertical maps by reduction.
    In particular, the vertical maps identify the kernels of the horizontal maps.
\end{remark}

\subsection{The obstruction group}\label{sub:ob_group}
The inclusion \(R\hookrightarrow\cA\) induces
\[
    \iota\colon\cU(R)=\wH^0(\Z_2;R^\times)
    \longrightarrow\cU(\cA)=\wH^0(\Z_2;\cA^\times).
\]
\begin{definition}\label{def:ob}
    The \emph{full obstruction group} is \(\cU(R)\), as in~\eqref{eq:ur_def}.
    The \emph{obstruction group} is \(\Ob(R)=\ker\iota\).
\end{definition}
Remark~\ref{rem:redQ} gives an isomorphism
\(\Ob(R)\cong\Ob(R_{\mathrm{red}})\).
Hence, from now on we may assume that~\(p\) is reduced.
Use the notation of Subsection~\ref{sub:reducible-qnorm}; then \(\cA=\cA_{\mathrm{red}}\) in~\eqref{eq:field_of_fractions}.
Lemma~\ref{lem:qnorm} identifies \(\cU(\cA)\) with the product of the norm quotients for the symmetric field factors.

For each~\(i\), let \(\Omega_i\) be the set of finite places of~\(F_i\) ramified in~\(E_i\), together with the real places of~\(F_i\) that become complex in~\(E_i\).
For \(s\in\Omega_i\), write \(\cU_{i,s}\) for the local norm quotient~\eqref{def:us} associated to~\(E_i/F_i\).
Localization induces
\[
    \lambda_i\colon F_i^\times/\Nm_{E_i/F_i}(E_i^\times)
    \longrightarrow\prod_{s\in\Omega_i}\cU_{i,s}.
\]
Using the identification of Lemma~\ref{lem:qnorm}, their product defines
\[
    \lambda\colon\cU(\cA)=(\cA^\tau)^\times/\Nm(\cA^\times)
    \longrightarrow\prod_i\prod_{s\in\Omega_i}\cU_{i,s}.
\]
When the discriminant is represented by an invariant unit of~\(R\), no local condition is needed at an unramified finite place.
Indeed, Proposition~\ref{prop:inert} applies to each symmetric component.

\begin{lemma}\label{lem:kernels}
    We have \(\ker(\lambda\circ\iota)=\ker\iota=\Ob(R)\).
\end{lemma}
\begin{proof}
    Let \([g]\in\cU(R)\), with \(g\in(R^\tau)^\times\).
    The class \(\lambda(\iota([g]))\) is trivial precisely when~\(g\) satisfies the finite ramified local norm conditions in every symmetric component and the positivity conditions at roots on~\(S^1\).
    By Proposition~\ref{prop:reducible-norm-tests}(1), these conditions are equivalent to \(g\in\Nm(\cA^\times)\), or equivalently \(\iota([g])=0\).
\end{proof}

We begin with a finiteness result.
\begin{proposition}
  The groups $\cU(R)$ and $\Ob(R)$ are finite.
\end{proposition}
\begin{proof}
  The ring $R=\Z[t^{\pm1}]/(p)$ is generated, as a $\Z$-algebra, by $t,t^{-1}$. In particular, it is a finitely generated commutative $\Z$-algebra. Since $p$ is reduced, $R$ embeds into the product of the domains $\Z[t^{\pm1}]/(f)$ over the irreducible factors~$f$ of~$p$, and Samuel's theorem \cite{Samuel}, applied to each factor, shows that $R^\times$ is finitely generated. Hence, $(R^\tau)^\times$ is finitely generated as a subgroup of $R^\times$. In particular, $\cU(R)$ is finitely generated as a quotient group of $(R^\tau)^\times$.

  Also, $\cU(R)$ is an elementary abelian $2$-group, so it is a finite group. Now, $\Ob(R)$ is a subgroup of $\cU(R)$, so it is finite.
\end{proof}

\subsection{Computation of $\Ob(R)$}\label{sub:obcomp}
To compute \(\Ob(R)\), we compare~\(R\) with its integral closure~\(\wt{R}\) in~\(\cA\).
Let~\(f\) range over all distinct irreducible factors of~\(p\), including both members of each pair \(q_j,\tau(q_j)\) in~\eqref{eq:field_of_fractions}.
Put \(E_f=\Q[t^{\pm1}]/(f)\) and \(R_f=\Z[t^{\pm1}]/(f)\), and define~\(T_f\) as in Remark~\ref{rem:normalization-general-factor}.
The normalization of a reduced ring is the product of the normalizations of its irreducible components.
Thus that remark gives
\[
    \wt{R}=\prod_f\wt{R}_f
    =\prod_f\cO_{E_f,T_f}\subset\prod_f E_f=\cA.
\]
Each~\(R_f\) is a quotient of~\(R\), and each~\(\wt{R}_f\) is finite as an \(R_f\)-module, so~\(\wt{R}\) is a finite \(R\)-module.

For a pair \(q_j,\tau(q_j)\), 
with $\tau(q_j)$ not associated with $q_j$,
the corresponding factor of $\wt{R}$ is $\wt{R}_{q_j}\times\wt{R}_{\tau(q_j)}\cong B_j\times B_j$, where
$B_j=\wt{R}_{q_j}$. The involution swaps the
components. In particular $\tau$ interchanges the two factors of $(B_j\times B_j)^\times$,
so by Lemma~\ref{lem:shapiro}
we get \(\cU(B_j\times B_j)=0\).
Thus these pairs contribute no summands to~\(\cU(\wt{R})\), although a norm preimage chosen in~\(\wt{R}^\times\) need not belong to~\(R^\times\).

Let~\(\mf\) be the conductor of~\(R\) in~\(\wt{R}\):
\[
    \mf=\{x\in\wt{R}:x\wt{R}\subseteq R\}
    =\Ann_R(\wt{R}/R).
\]
It is an ideal of both~\(R\) and~\(\wt{R}\).
\begin{lemma}\label{lem:finite_index}
    The additive quotient \(\wt{R}/R\) and the rings \(R/\mf\) and \(\wt{R}/\mf\) are finite.
\end{lemma}
\begin{proof}
    Choose finitely many \(R\)-module generators of~\(\wt{R}\).
    Since \(\cA=\Q\otimes_{\Z}R\), clearing their denominators gives a nonzero integer~\(n\) such that \(n\wt{R}\subset R\).
    Moreover, \((n\wt{R})\wt{R}\subset n\wt{R}\subset R\), so \(n\wt{R}\subset\mf\).

    As in the proof of Lemma~\ref{lem:integral_closure}, the condition \(p(1)=\pm1\) implies that \(R/nR\) is finite; this argument does not require irreducibility.
    Hence \(\wt{R}/n\wt{R}\) is finite, since it is a finite module over~\(R/nR\).
    Both \(\wt{R}/R\) and \(\wt{R}/\mf\) are quotients of this finite group, and \(R/\mf\) embeds into \(\wt{R}/\mf\).
\end{proof}

Let $j\colon R\to\wt{R}$ be
the embedding. We have the following Milnor square, see \cite[Example 1.2.6]{Weibel_K}:
\[
  \begin{tikzcd}
    R \ar[r,"j"]\ar[d] & \wt{R}\ar[d]\\
    R/\mf \ar[r] & \wt{R}/\mf.
  \end{tikzcd}
\]
From that square, we deduce the Unit-Pic long exact sequence, see \cite[Theorem I.3.10]{Weibel_K}:
\begin{equation}\label{eq:unit_pic}
  1\to R^\times\to\wt{R}^\times\times(R/\mf)^\times \to (\wt{R}/\mf)^\times\to \Pic(R)\to \Pic(\wt{R})\times\Pic(R/\mf)\to\Pic(\wt{R}/\mf).
\end{equation}
Now, by Lemma~\ref{lem:finite_index}, $R/\mf$ and $\wt{R}/\mf$ are finite rings, so $\Pic(R/\mf)$ and $\Pic(\wt{R}/\mf)$ are zero.
Write $\Xi=(\wt{R}/\mf)^{\times}/(R/\mf)^{\times}$. Then \eqref{eq:unit_pic} becomes
\[1\to R^{\times}\to\wt{R}^{\times}\xrightarrow{\rho}\Xi\to\Pic(R)\to\Pic(\wt{R})\to 1.\]
Define
\[\Xi_0=\im\rho=\ker(\Xi\to\Pic(R)).\]
Exactness shows that \(u\in\wt{R}^\times\) belongs to~\(R^\times\) precisely when \(\rho(u)=1\), or equivalently when its residue modulo~\(\mf\) lies in~\((R/\mf)^\times\).

The inclusions \(R\hookrightarrow\wt{R}\hookrightarrow\cA\) induce
\[
    \cU(R)\xrightarrow{j_*}\cU(\wt{R})
    \xrightarrow{\widetilde\iota}\cU(\cA),
    \qquad \iota=\widetilde\iota\circ j_*.
\]
Extending Definition~\ref{def:ob} to the normalization, put
\[
    \Ob(\wt{R})
    =\ker\left(\widetilde\iota:\cU(\wt{R})\to\cU(\cA)\right).
\]
The involution~$\tau$ preserves $\wt{R}$ and its conductor~$\mf$, since it preserves~$R$.
It therefore induces involutions on $\Xi$ and $\Xi_0$, and $\rho$ is equivariant.
Apply Lemma~\ref{lem:snake} to the equivariant short exact sequence
\begin{equation}\label{eq:conductor-units}
    1 \longrightarrow
    R^{\times} \longrightarrow
    \widetilde{R}^\times \xrightarrow{\rho}
    \Xi_0 \longrightarrow
    1.
\end{equation}
Using the two-periodicity of Tate cohomology, we obtain the exact sequence
\[
    \wH^1(\Z_2;\wt{R}^\times)
    \longrightarrow\wH^1(\Z_2;\Xi_0)
    \xrightarrow{\partial}\cU(R)
    \xrightarrow{j_*}\cU(\wt{R})
    \xrightarrow{\rho_*}\wH^0(\Z_2;\Xi_0).
\]
The connecting map has an explicit description.
Represent a class in $\wH^1(\Z_2;\Xi_0)$ by $\xi\in\Xi_0$ with $\xi\tau(\xi)=1$, and choose $u\in\wt{R}^\times$ such that $\rho(u)=\xi$.
Then
\[
    \partial[\xi]=[u\tau(u)]\in\cU(R).
\]
Indeed, $\rho(u\tau(u))=1$, so $u\tau(u)$ is an invariant unit of~$R$.
Replacing $u$ by another lift multiplies its norm by a norm from~$R^\times$.
Replacing $\xi$ by $\xi\eta/\tau(\eta)$, with $\eta\in\Xi_0$, allows us to replace
$u$ by $uv/\tau(v)$ for a lift $v$ of~$\eta$, leaving its norm unchanged.
Thus the formula is well-defined.
Furthermore, using the unit coordinates of Remark~\ref{rem:algorithm-effectivity} and the finite group~$\Xi_0$, these formulas compute the induced maps and the connecting map by integer linear algebra and multiplication of explicit units.

Define
\[
    \begin{aligned}
        K_0&=\coker\left(\wH^1(\Z_2;\wt{R}^\times)
            \longrightarrow\wH^1(\Z_2;\Xi_0)\right),\\
        K_1&=\ker\left(\rho_*:\cU(\wt{R})
            \longrightarrow\wH^0(\Z_2;\Xi_0)\right).
    \end{aligned}
\]
The connecting map~\(\partial\) identifies~\(K_0\) with~\(\ker j_*\), and exactness gives \(K_1=\im j_*\).
Concretely, $K_0$ consists of invariant units of~$R$ that are norms of units
of~$\wt{R}$, modulo norms of units of~$R$.
Such units are also norms in~$\cA$, which explains the inclusion $K_0\subset\Ob(R)$ below.
Thus there is a short exact sequence
\begin{equation}\label{eq:short}
    0\longrightarrow K_0\longrightarrow\cU(R)
    \xrightarrow{j_*}K_1\longrightarrow0.
\end{equation}
\begin{lemma}\label{lem:exact_K0}
    The group~\(K_0\) is naturally a subgroup of~\(\Ob(R)\), and restriction of~\(j_*\) gives a short exact sequence
    \[
        0\longrightarrow K_0\longrightarrow\Ob(R)
        \xrightarrow{j_*}K_1\cap\Ob(\wt{R})\longrightarrow0,
    \]
    where the intersection is taken inside~\(\cU(\wt{R})\).
\end{lemma}
\begin{proof}
    Since \(\iota=\widetilde\iota\circ j_*\), the subgroup \(K_0=\ker j_*\) lies in~\(\ker\iota=\Ob(R)\).
    The image of~\(\Ob(R)\) under~\(j_*\) lies in \(K_1\cap\ker\widetilde\iota=K_1\cap\Ob(\wt{R})\).
    Conversely, if \(y\in K_1\cap\Ob(\wt{R})\), exactness of~\eqref{eq:short} gives \(x\in\cU(R)\) with \(j_*(x)=y\).
    Then \(\iota(x)=\widetilde\iota(y)=0\), so \(x\in\Ob(R)\).
    The kernel of this restricted map is~\(K_0\), proving exactness.
\end{proof}

\section{Computations}\label{sec:families}

Write \(\Lambda=\Z[t^{\pm1}]\) and \(\Lambda_{\Q}=\Q[t^{\pm1}]\).
Unless stated otherwise, \(\Delta\in\Lambda\) is nonconstant, irreducible over~\(\Q\), symmetric, and normalized by~\(\Delta(1)=1\).
Put \(R=\Lambda/(\Delta)\), \(E=\Lambda_{\Q}/(\Delta)\), and let \(t_E\) be the image of~\(t\).
We use the same notation~\(\tau\) for the induced actions of the involution on ideals, ideal classes and places; in particular, \(F=E^\tau\) and \(\Nm(x)=x\tau(x)\).
If the extreme exponents of~\(\Delta\) are \(-d\) and~\(d\), its ordinary polynomial representative is
\[
    P(t)=t^d\Delta(t)\in\Z[t],\qquad \deg P=2d=[E:\Q],\qquad d=[F:\Q].
\]
We use~\(P\) when discussing polynomial degree, monicity, or polynomial discriminants.
For a ring $A$ with involution,
we write $\cU(A)=\wH^0(\Z_2;A^\times)$ as usual. Also, if $E$ is the fraction
field of $A$, with $\iota_A\colon \cU(A)\to\cU(E)$, we write
$\Ob(A):=\ker\iota_A$.

\subsection{The case $R=\cO_E$ or $R=\cO_{E,T}$}\label{sub:int_clos}
Calculations simplify if $R$ is integrally closed; in fact, in many cases $\cU(R)$ can be computed explicitly; see e.g. \cite[Proposition 2.6]{EBP_uni1}. However, $R=\cO_E$ implies in particular that $T$ is empty, that is, $\Delta$ is monic. In this subsection we translate
between $\cO_E$ and $\cO_{E,T}$, which makes more cases accessible by known methods; see Lemma~\ref{lemma:T-correction}.

We start by considering the case when~\(R = \cO_{E,T}\) is integrally closed.
In particular, in this case the obstruction group can be described in terms of arithmetic invariants of the ring~\(\cO_{E,T}\).
For \(A=R,\cO_E\), or~\(\cO_{E,T}\), write
\[
    \lambda_A=\lambda\circ\iota_A\colon\cU(A)\longrightarrow\prod_{s\in\Omega}\cU_s,
\]
where \(\Omega\) consists of the finite ramified places of~\(F\) and the real places that become complex in~\(E\), and \(\lambda\) is the localization map of Subsection~\ref{sub:ob_group}.
By Lemma~\ref{lem:kernels} and Proposition~\ref{prop:inert}, \(\ker\lambda_A=\Ob(A)\).
We identify each local norm quotient with~\(\{\pm1\}\) via its Hilbert symbol, so the trivial local invariant is~\(1\).
Recall that the \emph{class group} of~\(\cO_{E,T}\) is the quotient of the group of fractional ideals~\(I_{T}\) of~\(\cO_{E,T}\) by the group of principal fractional ideals~\(P_{T}\).
In other words, there is an exact sequence
\begin{equation}\label{eq:class-gp-def}
    1 \to P_{T} \to I_{T} \to \Cl(\cO_{E,T}) \to 1.
\end{equation}
Furthermore, the group~\(P_{T}\) fits into the following exact sequence
\begin{equation}\label{eq:short-exact-seq-princ-frac-id}
    1 \to \cO_{E,T}^{\times} \to E^{\times} \to P_{T} \to 1.
\end{equation}
We write \(\Cl(E)=\Cl(\cO_E)\) for the class group of~\(E\), and \(h(E)=|\Cl(E)|\) for its class number.
\begin{lemma}\label{lem:Cl-localization}
  Let $C_T\subset\Cl(E)$ be the subgroup generated by primes of $\cO_E$
  belonging to $T$. Then, the inclusion $\cO_E\to\cO_{E,T}$ induces an isomorphism
  \[\Cl(\cO_{E,T})=\Cl(E)/C_T.\]
\end{lemma}
\begin{proof}
  As $\cO_{E,T}$ is a localization of $\cO_E$, every fractional ideal
  of $\cO_{E,T}$ is of the form $\ma\cO_{E,T}$ for $\ma$ being a fractional
  ideal of $\cO_E$. Then, $\ma\cO_{E,T}=\cO_{E,T}$ if and only if $\ma$ is a product of powers of primes in $T$. This shows that $I_T$ is the quotient of the group of fractional ideals of~$\cO_E$ by the group generated by primes in $T$.

  The inclusion $\cO_E\to\cO_{E,T}$ takes principal ideals to principal ideals.  Now $\ma\cO_{E,T}$ is a principal ideal, $\ma\cO_{E,T}=x\cO_{E,T}$ if and only if $\ma x^{-1}$ is a product of primes in $T$. This shows the desired isomorphism.
\end{proof}
We remark that the isomorphism of Lemma~\ref{lem:Cl-localization} is clearly $\tau$-equivariant.

\begin{proposition}\label{prop:Ob-integrally-closed-case}
    Let $\Delta$ be irreducible with $R=\cO_{E,T}$ integrally closed, and let $\tau$ be the involution of $E$.
    Then
    \[\Ob(R) \cong \coker \Bigl(\bigoplus_{s=\tau(s)} \Z/2 \xrightarrow{\theta} \wH^{0} \left( \Z_{2}; \Cl(\cO_{E,T}) \right) \Bigr),\]
    the sum being over the \(\tau\)-fixed finite places outside $T$.
    In particular $\Ob(R)=0$ when $|\Cl(\cO_{E,T})|$ is odd.
    
    Furthermore, \(\theta\) maps the class represented by a \(\tau\)-invariant prime~\(\PP\) of~\(\cO_E\) outside~\(T\) to the class of~\(\PP\cO_{E,T}\) in~\(\wH^{0}(\Z_{2};\Cl(\cO_{E,T}))\), that is, to the image of~\([\PP]\in\Cl(E)\) in \(\Cl(E)/C_T\cong\Cl(\cO_{E,T})\); see Lemma~\ref{lem:Cl-localization}.
\end{proposition}

\begin{proof}
  Recall that
  \[
      \Ob(R)=\ker\left(\wH^{0}(\Z_{2};\cO_{E,T}^{\times})\to \wH^{0}(\Z_{2};E^{\times})\right).
  \]
  Apply Tate cohomology to the exact sequence~\eqref{eq:short-exact-seq-princ-frac-id}.
  Since $\wH^{1}(\Z_{2};E^{\times})=0$ by Hilbert 90, the connecting map identifies $\Ob(R)$ with $\wH^{1}(\Z_{2};P_{T})$.
  Now apply Tate cohomology to the exact sequence~\eqref{eq:class-gp-def}.
  Observe that~\(I_{T}\) is a free abelian group generated by finite places outside $T$.
  The involution permutes them, and $\wH^{\ast}$ vanishes on free orbits, so $\wH^{0}(\Z_{2};I_{T})=\bigoplus_{s=\tau(s)}\Z/2$ and $\wH^{1}(\Z_{2};I_{T})=0$.
  Consequently, the long exact sequence of Tate cohomology gives $\wH^{1}(\Z_{2};P_{T}) \cong \coker(\wH^{0}(\Z_{2};I_{T})\to\wH^{0}(\Z_{2};\Cl(\cO_{E,T})))$.
  The last assertion follows since Tate cohomology is killed by $2$.
\end{proof}
\begin{remark}\label{rem:Ob-concrete}
Assuming the hypotheses of Proposition~\ref{prop:Ob-integrally-closed-case}, that is, $R=\cO_{E,T}$, we compute the map
\[\wH^0(\Z_2;\Cl(R))\xrightarrow{\delta} \wH^1(\Z_2;P_T)\to \Ob(R)\subset \wH^0(\Z_2;R^\times)\]
from that proposition explicitly. 

To this end, recall from Subsection~\ref{sub:tate} that for any abelian group $G$ acted upon by $\Z_2$, $\wH^0(\Z_2;G)=G^{\tau}/(1+\tau)G$, where
$G^\tau$ is the fixed set, and $\tau$ is the generator of the action. Given this, an element in $\wH^0(\Z_2;\Cl(R))$
can be represented by an element in $\Cl(R)^\tau$: a fractional ideal $\ma$ such that $[\ma]=[\tau\ma]\in\Cl(R)$.
The quotient $\ma/\tau\ma$ is another fractional ideal. As it represents zero in $\Cl(R)$, it is a principal fractional ideal,
that is, there exists $z\in E^\times$ such that $zR=\ma/\tau\ma$.

Note that $\tau(z)R=\tau\ma/\ma$ so $z\tau(z)R=zR\tau(z)R=R$. In particular, $(1+\tau)(zR)=R$ is trivial in $P_T$. That is,
$zR$ is an element in $\wH^1(\Z_2;P_T)$. The class is identified with $\delta\ma$, the image of $\ma$
under the connecting homomorphism $\wH^0(\Z_2;\Cl(R))\to\wH^1(\Z_2;P_T)$.

The fact that $z\tau(z)R=R$ implies that $z\tau(z)\in R^\times$. We define the class $q\in \wH^0(\Z_2;R^\times)$ as $q=z\tau(z)$. 
Note that $q$ is a norm in $E^\times$, hence $q$ is in the kernel of the map $\wH^0(\Z_2;R^\times)\to\wH^0(\Z_2;E^\times)$.
In particular, $q\in\Ob(R)$.
We claim that $[\ma]\mapsto[q]$ is the map $\wH^0(\Z_2;\Cl(R))\to\Ob(R)$ displayed above, and that every class in $\Ob(R)$ arises in this way. 

To see this, consider $q\in R^\times\cap\Nm(E^\times)$. As $q\in\Nm(E^\times)$, $q=z\tau(z)$ for some $z\in E^\times$. 
Let $s$ be a place on $E$ outside $T$ and $\nu_s$ the
corresponding valuation. As $q$ is a unit, $\nu_s(q)=0$, so $\nu_s(z)+\nu_s(\tau(z))=0$, so $\nu_s(z)+\nu_{\tau(s)}(z)=0$.

If $\tau(s)=s$, that is, $s$ lies over an inert or ramified place of $F$, then $\nu_s(z)=0$. Otherwise, $\nu_{\tau(s)}(z)=-\nu_{s}(z)$.

Define now
\[\ma=\prod \mP_s^{\nu_s(z)},\]
where $\mP_s=\{\nu_s>0\}$, and the product is taken over all conjugate pairs $\{s,\tau(s)\}$, $s\ne\tau(s)$, of places outside $T$; we pick one place $s$ from each such pair.
By the definition, $\tau\ma=\prod(\tau\mP_s)^{\nu_s(z)}$, so $\ma/\tau\ma=zR$. We get an explicit choice of $\ma$, two different
choices differ by a $\tau$-invariant fractional ideal.
This is the procedure used in Examples~\ref{ex:arithmetic-only-ob} and~\ref{ex:arithmetic-all-silent}.

Finally, $q$ is a norm of an $R$-unit if and only if $zR$ is a coboundary in $P_T$, that is,
    \[
        zR=(yR)/\tau(yR)\qquad\text{for some }y\in E^\times.
    \]
    Indeed, this identity makes $z\tau(y)/y$ an $R$-unit of norm $q$.
    Conversely, if $q=\Nm(u)$ with $u\in R^\times$, Hilbert~90 gives $z/u=y/\tau(y)$ for some $y\in E^\times$, and hence the displayed
    ideal identity.
    The same argument shows that any two norm witnesses for $q$ determine the same class in $\wH^1(\Z_2;P_T)$.
\end{remark}

We can describe the difference between~\(\cU(\cO_{E})\) and~\(\cU(\cO_{E,T})\) in terms of the group~\(C_{T}\) with its induced involution.

\begin{lemma}\label{lemma:T-correction}
    Suppose $|C_T|$ is odd.
    Then inclusion induces isomorphisms
    \[\wH^i(\Z_2;\cO_E^\times)\xrightarrow{\cong}
      \wH^i(\Z_2;\cO_{E,T}^\times)\qquad\text{for every }i\in\Z.\]
    In particular, $\cU(\cO_E)\cong\cU(\cO_{E,T})$.
\end{lemma}
\begin{proof}
 Let \(M_{T} = \bigoplus_{\pp \in T} \Z \pp\) be the group of divisors supported on~$T$.
    The divisor map \(\div_{T} \colon \cO_{E,T}^\times \to M_{T}\), \(\div_{T}(x) = \sum_{\pp \in T} \nu_{\pp}(x) \pp\), has kernel~$\cO_E^\times$. Write \(D_{T}\) for its image.
    Since $C_T$ is generated by the classes of the primes in~$T$, we obtain two exact sequences:
    \[1 \to \cO_E^\times \to \cO_{E,T}^\times \xrightarrow{\div_T} D_T \to 0,\qquad 0 \to D_T \to M_T \to C_T \to 0.\]
    The involution acts freely on $T$ by Lemma~\ref{lem:norm_on_T}, so $\wH^*(\Z_2;M_T)=0$ by Lemma~\ref{lem:shapiro}.
    Applying Lemma~\ref{lem:snake} to the second sequence, we get $\wH^{i}(\Z_2;D_T)\cong\wH^{i-1}(\Z_2;C_T)$ for every~$i$.
    Each Tate cohomology group of $C_T$ is a subquotient of the finite group $C_T$ and is annihilated by~$2$; since $|C_T|$ is odd, all these groups vanish.
    Hence $\wH^*(\Z_2;D_T)=0$, and Lemma~\ref{lem:snake} applied to the first sequence proves the assertion.
\end{proof}

\subsection{CM fields and cyclotomic fields}\label{sub:CM_fields}

In this subsection we will be concerned with a special class of polynomials having all of their roots on the unit circle. This
includes, in particular, all Alexander polynomials of algebraic knots, and hence of torus knots.
First, we describe the general theory for arbitrary CM fields. Next, we pass to the theory of cyclotomic fields. 

\subsubsection{CM fields}
\label{sec:cm-fields}

\begin{definition}\label{def:CM-field}
  A number field $E$ is a \emph{CM field} if it is a totally imaginary quadratic extension of a totally real field $F$.
  Equivalently, $E$ admits an involution $\tau \neq \id$ which, under every embedding $\iota \colon E \hookrightarrow \C$, satisfies $\iota \circ \tau = \overline{\iota}$.

  Furthermore, \(F\) is the maximal totally real subfield of~\(E\).
\end{definition}

The following elementary lemma summarizes conditions which guarantee that a field generated by a given polynomial is a CM field.
\begin{lemma}
    With \(\Delta,E,F\) as above, write
    \[\Delta(t)=f(t+t^{-1}),\qquad P(t)=t^d f(t+t^{-1}),\qquad f\in\Z[X],\quad\deg f=d.\]
    Put \(u_F=t_E+t_E^{-1}\). Then
    \[F=\Q(u_F)\cong\Q[X]/(f),\]
    where the isomorphism sends the class of~\(X\) to~\(u_F\).
    In particular, \(f\) is irreducible over~\(\Q\), since \([F:\Q]=d\).
    The following are equivalent:
    \begin{enumerate}[label=\textup{(\roman*)}]
            \item $E$ is a CM field with maximal real subfield $F$;
            \item \label{item:CM-criterion-roots-unit-circle} every root of $\Delta$ lies on the unit circle;
            \item \label{item:CM-criterion-real-roots} all $\deg f$ roots of $f$ are real and lie in the interval $(-2,2)$;
            \item every archimedean place of $F$ is ramified.
    \end{enumerate}
\end{lemma}

We will sometimes refer to polynomials satisfying~\ref{item:CM-criterion-roots-unit-circle} as \emph{CM polynomials}.

\begin{remark}
    Observe that the condition~\ref{item:CM-criterion-real-roots} can be decided by Sturm's algorithm.
\end{remark}

Let $\mu_{E}$ be the group of roots of unity of $E$.
Let $Q=[\cO_{E}^{\times}:\mu_{E}\cO_{F}^{\times}]\in\{1,2\}$ be the \emph{Hasse unit index}, see~\cite[Theorem~4.12]{washington97}.
In general, $Q$ can be computed using \emph{regulators}, see \cite[Proposition 4.16]{washington97}.
For cyclotomic polynomials, there is an explicit formula for~\(Q\).

The following result is proved in \cite[Proposition 2.3]{EBP_uni1}.
\begin{proposition}\label{prop:order-norm-gp-CM-fields}
    If \(E\) is a CM field with maximal totally real subfield~\(F\), then $\left| \cU(\cO_{E}) \right|=2^d/Q$, where $d=[F:\Q]$.
\end{proposition}



\begin{remark}
    Let $E$ be a CM field with maximal real subfield $F$.
    Then $Q\in\{1,2\}$, by~\cite[Theorem~4.12]{washington97},
    and $Q=2$ if and only if the map $u\mapsto u/\tau(u)$ carries $\cO_{E}^{\times}$ onto $\mu_{E}$.

    In particular, if~\(F\) admits a finite non-dyadic place which ramifies in~\(E\), then \(Q = 1\), see~\cite[Remark 2.4]{EBP_uni1}.
\end{remark}

As in Subsection~\ref{sub:ob_group}, inclusion gives an exact sequence
\[0 \to \Ob(\cO_E) \to \cU(\cO_E) \xrightarrow{\iota_{\cO_E}} \cU(E)=\wH^{0}(\Z_{2};E^{\times}).\]
The composite \(\lambda_{\cO_E}=\lambda\circ\iota_{\cO_E}\) is described by signatures at real places of~\(F\) and Hilbert symbols at finite ramified places of~\(F\).
Let \(\iota_{1},\iota_{2},\ldots,\iota_{d}\), where \(d=[F:\Q]\), denote all the real places of~\(F\) (recall that~\(F\) is totally real, hence every archimedean place is real).
We can consider the \emph{signature homomorphism}
\[\sgn_{F} \colon \cO_{F}^{\times} \longrightarrow \{\pm1\}^{d}, \qquad a \longmapsto \left (\sgn \iota_{1}(a), \dots, \sgn \iota_{d}(a) \right),\]
where \(\sgn \iota_{i}(a)\) denotes the sign of the real number~\(\iota_{i}(a)\).
Recall that the kernel of the signature homomorphism is the so-called group of \emph{totally positive units} \(\cO^{\times,+}_{F} = \ker \sgn_F\).
Using the identification \(\F_2\cong\{\pm1\}\), \(e\mapsto(-1)^e\), the dimension of the image is called the \emph{unit signature rank} of~\(F\) and is denoted~\(s(F)\).
Additionally, we can consider the \emph{unit signature deficiency} \(\delta(F) = d - s(F)\).
Observe that the signature homomorphism descends to a homomorphism of finite \(\F_{2}\)-vector spaces
\[\sgn_{F} \colon \cO^{\times}_{F} / (\cO^{\times}_{F})^{2} \to \{\pm1\}^{d},\]
hence \(2^{\delta(F)} = \left| \cO^{\times,+}_{F} / (\cO^{\times}_{F})^{2} \right|\).

Observe that if \(x = \Nm(y)\), for \(x \in F^{\times}\) and \(y \in E^{\times}\), then
\[\iota_{i}(x) = \iota_{i}(\Nm(y)) = \iota_{i}(y \cdot \tau(y)) = \iota_{i}(y) \cdot \overline{\iota_{i}(y)} = \left| \iota_{i}(y) \right|^{2} > 0,\]
where we use~\(\iota_{i}\) to denote the extension of~\(\iota_{i}\) to a complex place of~\(E\).
This shows that \(\Nm(\cO^{\times}_{E}) \subset \ker \sgn_{F}\), hence the signature homomorphism descends to a map
\[\sigma_{F} \colon \cU(\cO_E) \to \{\pm1\}^{d}.\]

Consider
\[\cU^{+}(\cO_E) = \ker \sigma_{F} = \cO^{\times,+}_{F} / \Nm(\cO^{\times}_{E}).\]
As a result, we obtain a filtration of the group~\(\cU(\cO_E)\)
\[0 \subset \Ob(\cO_E) \subset \cU^{+}(\cO_E) \subset \cU(\cO_E).\]

\begin{theorem}\label{thm:Ob-group-CM-fields}
    Retain the notation and hypotheses above, and suppose that all roots of~\(\Delta\) lie on the unit circle.
    Thus \(E\) is a CM field with maximal totally real subfield~\(F\), and \(d=[F:\Q]\).
    Let~\(Q\) be the Hasse unit index and let~\(r\) denote the number of finite places of~\(F\) which ramify in~\(E\).
    \begin{enumerate}
            \item \label{item:Ob-group-CM-fields-1} \(Q\) divides \(2^{\delta(F)}\).
        In particular, the condition \(Q = 2\) forces \(\delta(F) \geq 1\).
            \item \label{item:Ob-group-CM-fields-2} We have \(\left| \cU^{+}(\cO_E) \right| = \frac{2^{\delta(F)}}{Q}\).
            \item \label{item:Ob-group-CM-fields-3} If \(r \geq 1\), the following inequality is satisfied
        \[\frac{2^{\delta(F)}}{Q \cdot 2^{r-1}} \leq \left| \Ob(\cO_E) \right| \leq \frac{2^{\delta(F)}}{Q}.\]
            \item \label{item:Ob-group-CM-fields-4} In particular, if \(r = 1\), then
        \[\Ob(\cO_E) = \cU^{+}(\cO_E), \quad \text{and} \quad \left| \Ob(\cO_E) \right| = \frac{2^{\delta(F)}}{Q}.\]
            \item \label{item:Ob-group-CM-fields-5} Similarly, if \(r = 0\) (e.g., if \(\Delta(-1) = \pm 1\)), then
        \[\Ob(\cO_E) = \cU^{+}(\cO_E).\]
        Furthermore, \(d = [F : \Q]\) is even and \(\delta(F) \geq 1\).
    \end{enumerate}
\end{theorem}
\begin{proof}
    For the proof of item~\ref{item:Ob-group-CM-fields-1}, observe that \(\Nm(\cO^{\times}_{E}) / (\cO^{\times}_{F})^{2} \subset \cO^{\times,+}_{F} / (\cO^{\times}_{F})^{2}\).
    As shown in the proof of \cite[Proposition~2.3]{EBP_uni1}, we have \(\Nm(\mu_{E} \cO^{\times}_{F}) = (\cO^{\times}_{F})^{2}\).
    Hence 
    \[Q = [\cO^{\times}_{E} : \mu_{E} \cO^{\times}_{F}] = [\Nm(\cO^{\times}_{E}) : (\cO^{\times}_{F})^{2}].\]
    Since \(2^{\delta(F)} = \left| \cO^{\times,+}_{F} / (\cO^{\times}_{F})^{2} \right|\), the desired conclusion follows.

    Item~\ref{item:Ob-group-CM-fields-2} follows from
    \[\cU^{+}(\cO_E) = \cO^{\times,+}_{F} / \Nm(\cO^{\times}_{E}) \cong \frac{\cO^{\times,+}_{F} / (\cO^{\times}_{F})^{2}}{\Nm(\cO^{\times}_{E}) / (\cO^{\times}_{F})^{2}}\]
    and the previous item.

    The right-hand inequality in item~\ref{item:Ob-group-CM-fields-3} follows from \(\Ob(\cO_E) \subset \cU^{+}(\cO_E)\).
    The left-hand inequality follows from the fact that by the results of Subsection~\ref{sec:ramified-places} the quotient \(\cU(\cO_E) / \Ob(\cO_E)\) is detected by signatures at real places of~\(F\) and Hilbert symbols at finite ramified places of~\(F\).
    Furthermore, by Theorem~\ref{thm:hasse_norm} the product of all Hilbert symbols is one, hence
    \[\left| \cU(\cO_E) / \Ob(\cO_E) \right| \leq 2^{s(F) + r-1}.\]
    Therefore,
    \[\left| \Ob(\cO_E) \right| \geq \frac{|\cU(\cO_E)|}{2^{s(F)+r-1}} = \frac{2^{d}}{Q \cdot 2^{s(F)+r-1}} = \frac{2^{\delta(F)}}{Q \cdot 2^{r-1}}.\]

    When \(r=1\), the inequality from the previous item shows that \(\left| \Ob(\cO_E) \right| = \frac{2^{\delta(F)}}{Q}\), hence \(\Ob(\cO_E) = \cU^{+}(\cO_E)\).

    When \(r=0\), arguing as in the proof of item~\ref{item:Ob-group-CM-fields-3}, we obtain \(\left| \Ob(\cO_E) \right| = \frac{2^{\delta(F)}}{Q}\).
    Furthermore, using once more the product formula from Theorem~\ref{thm:hasse_norm}, we see that for any \(a \in \cO^{\times}_{F}\)
    \[1 = \prod_{i=1}^{d} (a,(t_{E}-t_{E}^{-1})^{2})_{\iota_{i}} = \prod_{i=1}^{d} \sgn(\iota_{i}(a)),\]
    where the symbol \((-,-)_{\iota_{i}}\) denotes the Hilbert symbol with respect to the real place~\(\iota_{i}\).
    Taking \(a=-1\) shows that~\(d = [F : \Q]\) is even.
    Furthermore, the above formula shows that \(s(F) \leq d-1\), hence \(\delta(F) = d - s(F) \geq 1\).
\end{proof}

\begin{remark}\label{rem:computations-defect}
    The unit signature deficiency~\(\delta(F)\) can be computed from the ordinary and narrow class numbers of the totally real field~\(F\).
    Let \(I_F\) be the group of fractional ideals of~\(F\), \(P_F\) its subgroup of principal ideals, and \(P_F^+\) the subgroup of principal ideals generated by totally positive elements.
    The ordinary and \emph{strict} (or \emph{narrow}) class groups are
    \[\Cl(F)=I_F/P_F,\qquad \Cl^+(F)=I_F/P_F^+.\]
    Write \(h(F)=|\Cl(F)|\) and \(h^+(F)=|\Cl^+(F)|\).
    By~\cite[Theorem~IV.3.1]{Janusz}, we have
    \[\frac{h^+(F)}{h(F)}=2^{\delta(F)}.\]
    Both class numbers can be computed using PARI/GP; see Subsection~\ref{app:cyclotomic}.
\end{remark}

\subsubsection{Cyclotomic fields}
\label{sec:cyclotomic-fields}

An important class of CM polynomials is given by monic
CM polynomials. There is a full classification of such polynomials.
The following result is classical; we refer to \cite{washington97}
for the proofs.
\begin{proposition}\label{cor:CM-monic-maximal}
    With \(\Delta,P,R,E,F\) as above, suppose that \(E\) is a CM field with maximal real subfield~\(F\).
    The following are equivalent:
    \begin{enumerate}[label=\textup{(\roman*)}]
            \item $P$ is monic; \label{item:CM-monic-maximal-1}
            \item $T=\emptyset$; \label{item:CM-monic-maximal-2}
            \item $R$ is a $\Z$-order in $E$; \label{item:CM-monic-maximal-3}
            \item $P=\Phi_m$, equivalently $\Delta=t^{-\varphi(m)/2}\Phi_m$, for an $m$ with at least two distinct prime divisors; \label{item:CM-monic-maximal-4}
            \item $R=\cO_{E}$, the ring of algebraic integers of $E$. \label{item:CM-monic-maximal-5}
    \end{enumerate}
\end{proposition}

Proposition~\ref{prop:order-norm-gp-CM-fields} expresses the order of the group \(\cU(\Z[\xi_{n}])\) in terms of the degree of the extension \(\Q(\xi_{n}+\xi^{-1}_{n}) / \Q\) and the Hasse unit index~\(Q\).
For cyclotomic fields, we can obtain an explicit formula for~\(Q\) leading to the following proposition.

\begin{proposition}[{\cite[Example~2.5]{EBP_uni1}}]\label{prop:order-nomr-group-cyclotomic-case}
    Let \(E_{n} = \Q(\xi_{n})\) be the \(n\)-th cyclotomic field, where \(n\) has at least two distinct prime divisors. Then
    \[\left| \cU(\Z[\xi_{n}]) \right| =
        \begin{cases}
          2^{\varphi(n)/2}, & \text{if } n \text{ is even and } n/2 \text{ is an odd prime power}, \\
          2^{\varphi(n)/2-1}, & \text{otherwise},
        \end{cases}
    \]
    where~\(\varphi\) is Euler's totient function.
\end{proposition}
\begin{proof}
  Let \(c\) be the conductor of~\(E_n\), namely the smallest positive integer such that \(E_n\subset\Q(\xi_c)\), and put
  \[n^{\ast} =
      \begin{cases}
        n/2, & \text{if } n\equiv2\pmod4, \\
        n, & \text{otherwise}.
      \end{cases}
  \]
  By the discussion preceding Proposition~2.3 in~\cite[p.~10]{washington97}, we have \(E_n=E_{n^{\ast}}\), so \(c\leq n^{\ast}\).
  Conversely, \(\xi_{n^{\ast}}\in\Q(\xi_c)\), and every root of unity in~\(\Q(\xi_c)\) has the form \(\pm\xi_c^a\) by~\cite[Exercise~2.3]{washington97}.
  Thus \(n^{\ast}\mid\lcm(2,c)\).
  Since \(n^{\ast}\) is odd or divisible by~\(4\), this implies \(n^{\ast}\mid c\), hence \(c=n^{\ast}\).

  Applying~\cite[Corollary~4.13, p.~39]{washington97} with the index~\(n^{\ast}\), which is not congruent to~\(2\) modulo~\(4\), gives \(Q=1\) when \(n^{\ast}\) is a prime power and \(Q=2\) otherwise.
  Under the hypothesis on~\(n\), the first case occurs precisely when \(n=2p^k\) for an odd prime~\(p\).
  The asserted formula now follows from Proposition~\ref{prop:order-norm-gp-CM-fields}, since \([\Q(\xi_n+\xi_n^{-1}):\Q]=\varphi(n)/2\).
\end{proof}

By combining Proposition~\ref{prop:order-nomr-group-cyclotomic-case} and Theorem~\ref{thm:Ob-group-CM-fields} we obtain the following description of obstruction groups of cyclotomic rings.

\begin{corollary}\label{cor:Ob-groups-cyclotomic-polynomials}
    Let \(\Phi_{n}\) be the \(n\)-th cyclotomic polynomial, where \(n\) is divisible by at least two distinct primes.
    Let \(E_{n} = \Q(\xi_{n})\) be the cyclotomic field and \(F_{n} = \Q(\xi_{n}+\xi_{n}^{-1})\) be the maximal totally real subfield of \(E_{n}\).
    Let \(r\) be the number of finite primes of~\(F_n\) ramified in~\(E_n\), and let \(Q\) be the Hasse unit index of~\(E_n\).
    \begin{enumerate}
            \item If \(n = 2p^{k}\), where \(p\) is an odd prime, then \(\Phi_{n}(-1) = p\).
        Consequently:
        \[r = 1, \quad Q = 1, \quad \left| \cU(\Z[\xi_{n}]) \right| = 2^{\varphi(n)/2}\]
            \item If \(n\) is not of the form above, then \(\Phi_{n}(-1) = 1\), consequently
        \[r = 0, \qquad Q = 2, \qquad \left| \cU(\Z[\xi_{n}]) \right| = 2^{\varphi(n)/2-1}.\]
            \item Furthermore, in both cases
        \[\Ob(\Z[\xi_{n}]) = \cU^{+}(\Z[\xi_{n}]), \quad \text{and} \quad \left| \Ob(\Z[\xi_{n}]) \right| = \frac{2^{\delta(F_{n})}}{Q}.\]
    \end{enumerate}
\end{corollary}

\begin{proof}
    Recall that, for \(m>1\),
    \[\Phi_m(1)=
      \begin{cases}
        p,&m=p^k\text{ for a prime }p,\\
        1,&m\text{ has at least two distinct prime divisors}.
      \end{cases}\]
    This follows inductively from \(m=\prod_{d\mid m,\,d>1}\Phi_d(1)\).
    In particular, \(\Phi_n(1)=1\) if $n$ is divisible by at least two primes.
    By definition of $\Phi_n$, we also have
    \[\Phi_n(-1)=
      \begin{cases}
        \Phi_{2n}(1),&n\text{ odd},\\
        \Phi_{n/2}(1),&n\equiv2\pmod4,\\
        \Phi_n(1),&4\mid n.
      \end{cases}\]
    Consequently, \(\Phi_n(-1)=p\) when \(n=2p^k\) with \(p\) odd, and \(\Phi_n(-1)=1\) in all other cases allowed by the hypothesis.

    To determine the finite ramification, put \(u_n=\xi_n+\xi_n^{-1}\in\cO_{F_n}\).
    Since \(\cO_{E_n}=\Z[\xi_n]=\cO_{F_n}[\xi_n]\) and \(\xi_n\) satisfies \(X^2-u_nX+1=0\), the pair \((1,\xi_n)\) is an integral basis of~\(E_n/F_n\).
    The relative different is therefore
    \[\mD_{E_n/F_n}=(2\xi_n-u_n)=(\xi_n-\xi_n^{-1}).\]
    Its absolute ideal norm is
    \begin{align*}
      \Nm_{E_n/\Q}(\mD_{E_n/F_n})
      &=\left|\Nm_{E_n/\Q}(\xi_n-\xi_n^{-1})\right|\\
      &=\left|\Nm_{E_n/\Q}\left(\xi_n^{-1}(\xi_n-1)(\xi_n+1)\right)\right|\\
      &=\Phi_n(1)\Phi_n(-1)=\Phi_n(-1).
    \end{align*}
    If \(n=2p^k\), this norm is the prime~\(p\), so the different is a single prime ideal of~\(E_n\).
    Exactly the prime of~\(F_n\) below it ramifies, giving \(r=1\).
    Otherwise the different is the unit ideal, giving \(r=0\).

    As shown in the proof of Proposition~\ref{prop:order-nomr-group-cyclotomic-case}, the conductor of~\(E_n\) is \(n^\ast\), and \(Q=1\) precisely when \(n^\ast\) is a prime power.
    Since \(n\) has at least two distinct prime divisors, this occurs precisely when \(n=2p^k\) with \(p\) odd.
    Thus \(Q=1\) in the first case and \(Q=2\) in the second.
    Proposition~\ref{prop:order-norm-gp-CM-fields}, with \([F_n:\Q]=\varphi(n)/2\), now gives both asserted orders of~\(\cU(\Z[\xi_n])\).

    Finally, \(\Z[\xi_n]=\cO_{E_n}\), and \(r\) is either~\(0\) or~\(1\).
    Theorem~\ref{thm:Ob-group-CM-fields}, items~\ref{item:Ob-group-CM-fields-4} and~\ref{item:Ob-group-CM-fields-5}, therefore identifies \(\Ob(\Z[\xi_n])\) with \(\cU^+(\Z[\xi_n])\).
    Its order is \(2^{\delta(F_n)}/Q\) by item~\ref{item:Ob-group-CM-fields-2} of that theorem.
\end{proof}

\begin{example}
    Put \(R=\Z[\xi_n]\) and \(F=F_n=\Q(\xi_n+\xi_n^{-1})\). Table~\ref{tab:ob-groups-cyclotomic-polynomials} records norm quotients and obstruction groups for these cyclotomic rings.
    Observe that \(n = 39\) is the smallest example for which we obtain nontrivial \(\Ob(\Z[\xi_{n}])\).
    
    \begin{table}
        \centering
        \begin{center}
            \begin{tabular}{c | c | c | c | c | c | c | c | c | c | c | c}
              \hline
              \(n\) & \(\varphi(n)\) & \([F : \Q]\) & \(n^{\ast}\) & \(\Phi_n(-1)\) & \(Q\) & \(h(F)\) & \(h^+(F)\) & \(\delta(F)\) & \(\cU(R)\) & \(\Ob(R)\) & cert \\
              \hline
              $6$ & $2$ & $1$ & $3$ & $3$ & $1$ & $1$ & $1$ & $0$ & $\Z/2$ & $0$ & $1$ \\
              $10$ & $4$ & $2$ & $5$ & $5$ & $1$ & $1$ & $1$ & $0$ & $(\Z/2)^{2}$ & $0$ & $1$ \\
              $12$ & $4$ & $2$ & $12$ & $1$ & $2$ & $1$ & $2$ & $1$ & $\Z/2$ & $0$ & $1$ \\
              $14$ & $6$ & $3$ & $7$ & $7$ & $1$ & $1$ & $1$ & $0$ & $(\Z/2)^{3}$ & $0$ & $1$ \\
              $15$ & $8$ & $4$ & $15$ & $1$ & $2$ & $1$ & $2$ & $1$ & $(\Z/2)^{3}$ & $0$ & $1$ \\
              $20$ & $8$ & $4$ & $20$ & $1$ & $2$ & $1$ & $2$ & $1$ & $(\Z/2)^{3}$ & $0$ & $1$ \\
              $21$ & $12$ & $6$ & $21$ & $1$ & $2$ & $1$ & $2$ & $1$ & $(\Z/2)^{5}$ & $0$ & $1$ \\
              $30$ & $8$ & $4$ & $15$ & $1$ & $2$ & $1$ & $2$ & $1$ & $(\Z/2)^{3}$ & $0$ & $1$ \\
              $39$ & $24$ & $12$ & $39$ & $1$ & $2$ & $1$ & $4$ & $2$ & $(\Z/2)^{11}$ & $\Z/2$ & $1$ \\
              $56$ & $24$ & $12$ & $56$ & $1$ & $2$ & $1$ & $4$ & $2$ & $(\Z/2)^{11}$ & $\Z/2$ & $1$ \\
              $58$ & $28$ & $14$ & $29$ & $29$ & $1$ & $1$ & $8$ & $3$ & $(\Z/2)^{14}$ & $(\Z/2)^{3}$ & $1$ \\              
              \hline
            \end{tabular}
        \end{center}
        \caption{Norm quotients and obstruction groups for cyclotomic rings. The class numbers \(h(F)\), \(h^+(F)\) and signature deficiency \(\delta(F)\) refer to the maximal real subfield~\(F=F_n\). The last column records certification independent of the Generalized Riemann Hypothesis.}
        \label{tab:ob-groups-cyclotomic-polynomials}
    \end{table}
    
    All the ingredients required to compute the groups \(\cU(\Z[\xi_{n}])\) and~\(\Ob(\Z[\xi_{n}])\) can be computed easily; see Appendix~\ref{app:cyclotomic}.
\end{example}

\subsection{Quadratic fields}
\label{sec:imag-quadr-fields}

The simplest non-trivial Alexander polynomial of a knot has degree two. In this subsection, we describe the computations of the $\cU(R)$
 and $\Ob(R)$ groups for $R=\Z[t^{\pm 1}]/(\Delta)$, where $\Delta$
 is quadratic.

 \subsubsection{General theory}\label{sub:general}
We begin with general theory. For a parameter $a\in\Z\setminus\{0\}$,
we consider
\[\Delta_a(t)=a(t+t^{-1})-(2a-1),\ P_a(t)=t\Delta_a(t)=at^2-(2a-1)t+a.\]
With this notation \(\Delta_a(1)=1\), \(D=\disc P_a=1-4a=\Delta_a(-1)\).
In most cases, we will be assuming 
that~\(D\) is not a square, so that~\(P_a\) is irreducible.
Set, as usual, \(R_a=\Lambda/(\Delta_a)=\Lambda/(P_a)\) and \(E_a=\Lambda_{\Q}/(\Delta_a)\), and denote the fixed field by~\(F_a\). Let also $R_0=R_a^\tau$ be the fixed ring. When $a$ is fixed and understood, we will drop the subscript $a$ from the notation. We write $t_E$ for the image of $t$ in $R$, and set $u_F=t_E+t_E^{-1}$. 
As usual, we denote by~\(t_{E}\) the image of the variable~\(t\) in~\(E_{a}\) and set \(u_F = t_E+t_E^{-1}\).

The element
$\beta=a(1-t_E)$
satisfies \(\beta^2-\beta+a=0\), so \(E_a=\Q(\sqrt D)\); when \(D\) is square-free, it is the discriminant of~\(E_a\).
Thus \(E_a\) is imaginary quadratic if \(a\geq1\), and real quadratic if \(a\leq-1\).
The following proposition summarizes the properties of~\(E_a\) and~\(R_a\) that we need.
\begin{proposition}\label{prop:im-quad-family}
    In the setup above, the following hold.
    \begin{enumerate}[label=\textup{(\roman*)}]
            \item \label{item:im-quad-family-1} The set~\(T\) consists of the places of~\(E_{a}\) lying over the primes dividing~\(a\), and every such prime splits in~\(E_{a}\).
            \item \label{item:im-quad-family-2} We have \(\widetilde{R}_{a} = \cO_{E,T} = \cO_{E}[1/a]\) and \(R_0=\Z[1/a]\).
        Furthermore, \(R_{a} = \widetilde{R}_{a}\) if and only if \(D\) is square-free.
    \end{enumerate}
\end{proposition}
\begin{proof}
  Let $s$ be a finite place of $E_a$. By definition, $s\in T$ if
  and only if $\nu_s(t_E)\neq 0$. Assume that $s$ lies over a prime $p\in\Z$.
    If \(p\nmid a\), then \(t_E^{\pm1}\) are roots of \(X^2-u_FX+1\in\Z_{(p)}[X]\), so they are integral at~\(s\), and \(s\notin T\).
    Conversely, if \(p\mid a\), then $2a-1$ is invertible modulo $a$
    and so $\nu_s(2a-1)=0$. Then, $u_F=(2a-1)/a$ implies that \(\nu_s(u_F)=-\nu_s(a)<0\).
    That is, it is impossible that $\nu_s(t_E)=\nu_s(t_E^{-1})=0$,
    so $s\in T$.
    Every such prime~\(p\) splits in~\(E_a\) by Lemma~\ref{lem:norm_on_T}.
    This proves item~\ref{item:im-quad-family-1}.

    To prove the second part,
    we observe that by Theorem~\ref{thm:normalization} and item~\ref{item:im-quad-family-1}, we have \(\wt{R}_a=\cO_{E,T}=\cO_E[1/a]\).
    Since \(au_F=2a-1\), we have \(1/a=2-u_F\in R_0\), so
    \(\Z[1/a]\subset R_0\subset\wt{R}_a\cap\Q=\Z[1/a]\).
    Finally, $\Z[1/a,\beta]\subset R_a$, and the relations \(t_E=1-\beta/a\), $t_E^{-1}=u_F-t_E=1+(\beta-1)/a$ give \(R_{a}\subset \Z[1/a,\beta]\), so $R_a=\Z[1/a,\beta]$.

    If \(D\) is square-free, then \(\cO_{E} = \Z[\beta]\), because \(\disc\Z[\beta]=D\); hence \(R_a=\cO_{E}[1/a]\).
    Otherwise, the index \([\cO_{E}:\Z[\beta]]>1\) is divisible only by primes~\(p\) with \(p^2\mid D\).
    Since \(\gcd(a,D)=1\), such primes are not inverted, so \(\Z[1/a, \beta] \subsetneq \cO_{E}[1/a]\).
\end{proof}

Throughout the rest of Subsection~\ref{sub:general},
we will be assuming that 
\(D\) is square-free, so that \(R=\cO_{E,T}=\cO_E[1/a]\) and \(D\) is the discriminant of~\(E\).
Since \(F=\Q\), the invariant units of~\(\cO_E\) are~\(\pm1\); hence \(\cU(\cO_E)\) is generated by~\([-1]\), and \([-1]=0\in\cU(\cO_E)\) precisely when \(\cO_E\) has a unit of norm~\(-1\). We have two cases.

\begin{itemize}
  \item  If \(a\geq1\), all norms are positive, and \(\cU(\cO_E)\cong\Z/2\).
  \item If \(a\leq-1\), then \(\cO_E^\times=\{\pm\eta_a^k\colon k\in\Z\}\), where the fundamental unit is \(\eta_a=\frac12(x+y\sqrt D)\) for the positive integer solution \((x,y)\) of \(x^{2}-Dy^{2}=\pm4\) with the smallest \(y\), choosing the smallest \(x\) if there is a tie. Thus \(\eta_a\) is the least unit greater than~\(1\) under the embedding with \(\sqrt D>0\).
Thus, for \(a\leq-1\), we have \(\cU(\cO_E)\cong\Z/2\) if \(\Nm(\eta_a)=1\), and \(\cU(\cO_E)=0\) if \(\Nm(\eta_a)=-1\).
\end{itemize}

\begin{theorem}\label{thm:norm-gps-im-quad-fields}
  Suppose $a\neq 0$ and \(D=1-4a\) is square-free.  
  Every prime number $\ell\mid a$ lies in $R_0^\times$. Then
  the group \(\cU(R_a)\) is generated by \([-1]\) and the classes~\([\ell]\). 

  With $C_T$ as in Lemma~\ref{lem:Cl-localization},
  and $\PP_\ell$ denoting a prime of $\cO_E$ above $\ell$
  the homomorphism
  \[\kappa\colon R_0^\times\longrightarrow C_T/2C_T,\qquad
    \kappa\Bigl(\pm\prod_{\ell}\ell^{n_\ell}\Bigr)=\sum_{\ell}n_\ell\,[\PP_\ell]\bmod 2C_T,\]
  vanishes on \(\Nm(R_a^\times)\), and so
  it fits into the following exact sequence
  \begin{equation}\label{eq:quad-unit-norm-sequence}
    \langle[-1]\rangle\longrightarrow\cU(R_a)
    \xrightarrow{\kappa}C_T/2C_T\longrightarrow0.
  \end{equation}
  In particular, \(\cU(R_a)\) is generated by \([-1]\) and the classes~\([\ell]\).

  The class \([-1]\) is nonzero if and only if \(-1\notin\Nm(\cO_{E,T}^{\times})\); in particular, it is nonzero if \(a\geq1\).
\end{theorem}
\begin{proof}    
    Let \(\ell_1,\ldots,\ell_k\) be the distinct prime divisors of~\(a\).
    By Proposition~\ref{prop:im-quad-family}, \(R_0=\Z[1/a]\), so
    the invertible elements of $R_0$ are generated by $-1$
    and the prime factors of $a$. Each square of such an element is clearly a norm, so  $\cU(R_a)$
    is generated by \([-1]\) and the classes~\([\ell_i]\).

    Take a prime number $\ell|a$ and let $\PP_\ell$ be a prime of $\cO_E$ above $\ell$. We have $\ell\cO_E=\PP_\ell\tau(\PP_\ell)$. Note that since $E$ is quadratic, $[E:\Q]=2$, so there are at most two ideals above $\ell$. Since $\PP_\ell\in T$, the prime $\ell$ splits by Proposition~\ref{prop:im-quad-family}, so $\PP_\ell\neq \tau(\PP_\ell)$.

    For every fractional ideal $\ma$ of $E$, $\ma\tau(\ma)=\Nm(\ma)\cO_E$ is principal, we have $\tau(\ma)=-\ma$ in $\Cl(E)$. In particular,
    $\kappa$ identifies $\PP_\ell$ with $\tau(\PP_\ell)$. Also,
    $\kappa$ is a surjection. Essentially the same argument
    shows that $\kappa$ vanishes on $\Nm(R^\times)$.
    So $\kappa$ induces a surjection from $\cU(R)$ to $C_T/2C_T$.
    Moreover, $\kappa$ vanishes on $[-1]$.

    Conversely, let \(q=\pm\prod_\ell\ell^{n_\ell}\) satisfy \(\kappa(q)=0\). With the explicit formula for $\kappa$, we see that 
    this holds if and only if \(\sum_\ell n_\ell[\PP_\ell]=2c\) for some \(c\in C_T\).
    Write \(c=[\mb]\) with \(\mb\) supported on~\(T\); then \(2c=[\mb\,\tau(\mb)^{-1}]\).
    Hence the ideal \(\ma:=\prod_\ell\PP_\ell^{n_\ell}\cdot\tau(\mb)\,\mb^{-1}\) is principal. Take $x$ such that \(\ma=x\cO_E\).
    As \(\ma\) is supported on~\(T\), we have \(x\in\cO_{E,T}^\times=R_a^\times\), and
    \[\Nm(x)\cO_E=\ma\,\tau(\ma)=\prod_\ell\left(\PP_\ell\tau(\PP_\ell)\right)^{n_\ell}=q\cO_E.\]
    Therefore \(\Nm(x)=\pm q\), and \([q]\in\langle[-1]\rangle\). This proves the exactness of~\eqref{eq:quad-unit-norm-sequence}.
\end{proof}
\begin{remark}
    By definition, \([-1]=0\) if and only if \(-1\in\Nm(\cO_{E,T}^\times)\).
    If \(a\geq1\), the field~\(E\) is imaginary quadratic, so every norm from~\(E^\times\) is positive, and \([-1]\neq0\). As
    we see in Example~\ref{ex:real-quad-sign}, $[-1]$ can be zero in $\cU(R_a)$ also if $a<0$.
\end{remark}

\subsubsection{Particular examples}
\begin{example}\label{ex:real-quad-sign}
    Suppose \(a=-55\). Then \(D=221=13\cdot17\) is square-free,
    \(E=\Q(\sqrt{221})\), and \(R=\cO_E[1/55]\) by Proposition~\ref{prop:im-quad-family}.
    The smallest solution of \(x^2-221y^2=\pm4\) with \(y>0\) is \((x,y)=(15,1)\), since \(221+4=15^2\).
    Hence \(\eta=\frac12(15+\sqrt{221})\) is a fundamental unit, of norm~\(1\), and \([-1]\) is nontrivial in \(\cU(\cO_E)\).
    On the other hand, the element
    \[z=\frac{14+\sqrt{221}}5,\qquad
      z^{-1}=\frac{\sqrt{221}-14}5,\qquad
      \Nm(z)=\frac{196-221}{25}=-1,\]
    is a unit of \(R=\cO_E[1/55]\), because \(14+\sqrt{221}\) and \(\sqrt{221}-14\) lie in~\(\cO_E\).
    In this situation, \([-1]=0\) in \(\cU(R)\).
\end{example}

\begin{example}\label{ex:74}
    Take \(a=4\). Then $\Delta_4=4(t+t^{-1})-7$, $D=-15$ is square-free,
    so $E=\Q(\sqrt{-15})$, and $T$ consists of the two places over~$2$, which splits in~$E$.
    The class number of~$E$ is two, and the primes above~$2$ are not principal: with $\omega=\frac12(1+\sqrt{-15})$, the norm form $x^2+xy+4y^2$ of $\cO_E=\Z[\omega]$ does not represent~$2$.
    Hence $C_T=\Cl(E)\cong\Z/2$, and Theorem~\ref{thm:norm-gps-im-quad-fields} gives
    \[\cU(R)=\{[1],[-1],[2],[-2]\}\cong(\Z/2)^{2}.\]
    In particular, neither $2$ nor $-2$ is the norm of a unit of~$R$.

    By Subsection~\ref{sub:int_clos}, $\Ob(R)=\ker\lambda_R$, where $\lambda_R$ records the Hilbert symbols $(z,D)_s$ at the ramified places $s=\infty,3,5$.
    For $z\in R_0^\times=\{\pm2^m\}$, these symbols are $\sgn(z)$, $\left(\frac{z}{3}\right)$ and $\left(\frac{z}{5}\right)$, by~\eqref{eq:Hilbert-symbols-local-formula}.
    Hence
    \[
      \lambda_R([2]) = (1,-1,-1),\qquad
      \lambda_R([-1]) = (-1,-1,1).
    \]
    These vectors are independent, so $\lambda_R$ is injective and $\Ob(R)=0$.
\end{example}

\begin{example}\label{ex:im-quad-fields-nontrivial-Ob}
    Consider \(a=109^2=11881\). The discriminant
    \[D=-47523=-3\cdot7\cdot31\cdot73\]
    is square-free, so 
    previous arguments apply. There are two primes of $\cO_E$ above $109$;
    call one of them $\PP$. Then $T=\{\PP,\tau(\PP)\}$ and $\PP\tau(\PP)=(109)$.

    A PARI/GP computation (see Appendix~\ref{app:quadratic}) gives
    \[\Cl(E)\cong\Z/8\oplus\Z/2\oplus\Z/2,\qquad [\PP]=(4,0,0).\]
    Hence \(C_T=\langle[\PP]\rangle\cong\Z/2\) and \(2C_T=0\), so \(\kappa([109])\ne0\), and
    \[\cU(R)=\langle[-1],[109]\rangle\cong(\Z/2)^2.\]
    
    To determine $\Ob(R)$,
    we follow Example~\ref{ex:74}. This time, the class~\([109]\) belongs to \(\Ob(R)=\ker\lambda_R\), since \(109>0\) and
    \[
      \left(\frac{109}{3}\right)=\left(\frac{1}{3}\right),\quad
      \left(\frac{109}{7}\right)=\left(\frac{4}{7}\right),\quad
      \left(\frac{109}{31}\right)=\left(\frac{16}{31}\right),\quad
      \left(\frac{109}{73}\right)=\left(\frac{36}{73}\right)
    \]
    are all equal to~\(1\).
    The classes \([-1]\) and \([-109]\) do not belong to \(\Ob(R)\), because negative numbers are not norms from the imaginary quadratic field~\(E\).
    Hence \(\Ob(R)=\langle[109]\rangle\cong\Z/2\).
\end{example}
\begin{example}\label{ex:quad-a6}
  In this example, $|C_T|$ is odd.

    Take \(a=6\), so that \(P_6=6t^2-11t+6\) and \(D=-23\).
    Then \(E=\Q(\sqrt{-23})\) and \(\cO_E=\Z[\omega]\), where \(\omega=\frac12(1+\sqrt{-23})\), and \(\Nm(x+y\omega)=x^2+xy+6y^2\).
    By Proposition~\ref{prop:im-quad-family}, we have \(R=\cO_E[1/6]\). The set $T$ consists of the primes above the factors of $6$:
    \(T=\{\PP_2,\tau(\PP_2),\PP_3,\tau(\PP_3)\}\), where \(\PP_2=(2,1+\omega)\) and \(\PP_3=(3,\omega-1)\).

    Since \(1+\omega\in\PP_2\), \(1+2\omega=3+2(\omega-1)\in\PP_3\) and \(1-\omega\in\PP_2\cap\PP_3\), while \(2\nmid1+\omega\) and \(3\nmid1+2\omega\), comparing the norms \(8\), \(27\) and~\(6\) gives
    \[(1+\omega)\cO_E=\PP_2^3,\qquad (1+2\omega)\cO_E=\PP_3^3,\qquad (1-\omega)\cO_E=\PP_2\PP_3.\]
    As \(x^2+xy+6y^2\geq\frac{23}4y^2\), no element of~\(\cO_E\) has norm~\(2\), so \(\PP_2\) is not principal.
    Hence \(C_T=\langle[\PP_2]\rangle\cong\Z/3\).
    By the Minkowski bound~\cite[Section~I.6]{Neukirch}, every ideal class of~\(E\) contains an integral ideal of norm at most \(\frac2\pi\sqrt{23}<4\); such ideals are products of primes in~\(T\), so in fact \(C_T=\Cl(E)\).
    As \(|C_T|\) is odd, we have \(C_T/2C_T=0\), and Theorem~\ref{thm:norm-gps-im-quad-fields} gives \(\cU(R)=\langle[-1]\rangle\cong\Z/2\).
    Since \(E\) is imaginary, \(-1\) is not a norm from~\(E^\times\), so \(\Ob(R)=0\).
    Explicitly, the units
    \[2=\Nm\Bigl(\frac{1+\omega}2\Bigr),\qquad 3=\Nm\Bigl(\frac{1+2\omega}3\Bigr)\]
    are norms of units of~\(R\). 
\end{example}
    In Tables~\ref{tab:im-quad-fields} and~\ref{tab:real-quad-fields} we present computations of the norm groups for imaginary and real quadratic fields; the code is listed in Appendix~\ref{app:quadratic}.
    The case~\(a = -4861\) is a negative case in which \(\Ob(R_a)\) is nontrivial.

\begin{table}
    \centering
    \begin{tabular}{cccccccc}
      \hline
      $a$ & $P_a$ & $|D|$ & $E_a$ & $h(E_a)$ & $T$-correction & $\cU(R_a)$ & $\Ob(R_a)$\\
      \hline
      $1$ & $t^{2}-t+1$    & $3$  & $\Q(\sqrt{-3})$  & $1$ & trivial     & $\Z/2$ & \(0\) \\
      $2$ & $2t^{2}-3t+2$  & $7$  & $\Q(\sqrt{-7})$  & $1$ & trivial     & $\Z/2$ & \(0\) \\
      $3$ & $3t^{2}-5t+3$  & $11$ & $\Q(\sqrt{-11})$ & $1$ & trivial     & $\Z/2$ & \(0\) \\
      $4$ & $4t^{2}-7t+4$  & $15$ & $\Q(\sqrt{-15})$ & $2$ & nontrivial & $(\Z/2)^{2}$ & \(0\) \\
      $5$ & $5t^{2}-9t+5$  & $19$ & $\Q(\sqrt{-19})$ & $1$ & trivial     & $\Z/2$ & \(0\) \\
      $6$ & $6t^{2}-11t+6$ & $23$ & $\Q(\sqrt{-23})$ & $3$ & trivial     & $\Z/2$ & \(0\) \\
      $8$ & $8t^{2}-15t+8$ & $31$ & $\Q(\sqrt{-31})$ & $3$ & trivial     & $\Z/2$ & \(0\) \\
      \(109^{2}\) & & \(3 \cdot 7 \cdot 31 \cdot 73\) & \(\Q(\sqrt{-47523})\) & \(32\) & nontrivial & \((\Z/2)^{2}\) & \(\Z/2\) \\
      \hline
    \end{tabular}
    \caption{The table presents computations of norm groups \(\cU(\cO_{E,T})\) for selected \(a\geq1\). The polynomial column lists \(P_a=t\Delta_a\), and \(|D|\) is the knot determinant.}
    \label{tab:im-quad-fields}
\end{table}

\begin{table}
    \centering
    \resizebox{\linewidth}{!}{%
      \begin{tabular}{cccccccc}
        \hline
        $a$ & $P_a$ & $|D|$ & $E_a$ & $h(E_a)$ & $T$-correction & $\cU(R_a)$ & $\Ob(R_a)$ \\
        \hline
        $-1$ & $-t^{2}+3t-1$ & $5$ & $\mathbb{Q}(\sqrt{5})$ & $1$ & trivial & $0$ & $0$ \\
        $-5$ & $-5t^{2}+11t-5$ & $21$ & $\mathbb{Q}(\sqrt{21})$ & $1$ & trivial & $\mathbb{Z}/2$ & $0$ \\
        $-55$ & $-55t^{2}+111t-55$ & $221$ & $\mathbb{Q}(\sqrt{221})$ & $2$ & nontrivial & $\mathbb{Z}/2$ & $0$ \\
        $-4861$ & $-4861t^{2}+9723t-4861$ & $19445$ & $\mathbb{Q}(\sqrt{19445})$ & $2$ & trivial & $\mathbb{Z}/2$ & $\mathbb{Z}/2$ \\
        \hline
      \end{tabular}%
    }
    \caption{Certified norm-group computations for quadratic fields with square-free $D=1-4a$. The polynomial column lists $P_a=t\Delta_a$, and $|D|$ is the knot determinant. The $T$-correction records whether localization changes the norm quotient.}\label{tab:real-quad-fields}
\end{table}

\begin{example}\label{ex:quad-conductor-16}
    Set $a=16$. Let
    \[
        \Delta_{16}(t)=16(t+t^{-1})-31,\qquad P_{16}(t)=t\Delta_{16}(t)=16t^2-31t+16.
    \]
    Then $D=\Delta_{16}(-1)=-63=3^2\cdot(-7)$ is not square-free, so Theorem~\ref{thm:norm-gps-im-quad-fields} does not apply to~$R=R_{16}$.
    We compare $R$ with its normalization~$\wt{R}$, as in Subsection~\ref{sub:obcomp}.
    We will show that
    \[\cU(R)\cong(\Z/2)^2,\qquad \Ob(R)\cong\Z/2,\qquad
      \cU(\wt R)\cong\Z/2,\qquad \Ob(\wt R)=0,\]
    with the additional class in~$\cU(R)$ represented by~$2$.

    Consider $\beta=16(1-t_E)$. It satisfies $\beta^2-\beta+16=0$, so $E\cong\Q(\sqrt{-63})=\Q(\sqrt{-7})$. With $\omega=\frac12(1+\sqrt{-7})$, we have $\omega^2-\omega+2=0$, $\omega\tau(\omega)=2$ and
    $\cO_E=\Z[\omega]$.
    The two isomorphisms $E\cong\Q(\sqrt{-7})$ differ by~$\tau$; we fix the one for which
    \[\beta=\frac{1-3\sqrt{-7}}2=2-3\omega.\]
    Put $B=\Z[1/2]$.
    By Proposition~\ref{prop:im-quad-family}, $\wt{R}=\cO_E[1/2]=B\oplus B\omega$ and $R_0=B$, so $R_0^\times=\{\pm2^m:m\in\Z\}$.
    Since $R=B[\beta]$, we get 
    \[R=B\oplus3B\omega.\]

    In particular, $[\wt{R}:R]=3$, and the conductor of~$R$ in~$\wt{R}$ is $\mf:=3\wt{R}$. Indeed, an element $x=u+v\omega$ of $\wt{R}$ satisfies $x\wt{R}\subset R$ if and only if $v\in 3B$ and 
    $x\omega=-2v+(u+v)\omega\in R$, meaning that $u+v\in 3B$, so $u\in 3B$.

    The norm of $u+v\omega$ is $u^2+uv+2v^2$, which is positive.
    Since $\omega\in\wt R^\times$ and $\Nm(\omega)=2$, we have $\Nm(\wt R^\times)=\{2^n:n\in\Z\}$.
    If $x=u+3v\omega\in R^\times$, then $\Nm(x)=2^m$ for some~$m$, and
    \[2^m=u^2+3uv+18v^2\equiv u^2\pmod3.\]
    As $2^m\not\equiv0\pmod 3$, we get $2^m\equiv1\pmod3$, so $m$ is even.
    Since $\Nm(2)=4$, it follows that $\Nm(R^\times)=\{4^n:n\in\Z\}$.
    Therefore
    \[\cU(R)=\langle[-1],[2]\rangle\cong(\Z/2)^2,
      \qquad
      \cU(\wt R)
              =\langle[-1]\rangle\cong\Z/2.\]
    Every positive~$2^m$ is the norm $\Nm(\omega^m)$, whereas negative numbers are not norms from~$E^\times$; hence
    \[\Ob(R)=\langle[2]\rangle\cong\Z/2,
      \qquad \Ob(\wt R)=0.\]
    In the notation of~\eqref{eq:short}, the map $j_*\colon\cU(R)\to\cU(\wt R)$ is onto, so $K_1=\cU(\wt R)$, and $K_0=\ker j_*=\langle[2]\rangle$: the number~$2$ is the norm of a unit of~$\wt R$, but not of a unit of~$R$.
\end{example}

\subsection{Quartic fields with nontrivial obstruction}\label{sub:quartic-obstructions}

We apply Proposition~\ref{prop:Ob-integrally-closed-case} to two integrally closed rings with nonzero obstruction group.
We keep the notation $R,\wt{R}$ and $E,F$ and $t_E$. Set $u_F=t_E+t_E^{-1}$ and $R_0=R^\tau$.
In the examples below, $P=t^2\Delta$ is the polynomial representative of the Alexander polynomial.
The knot-theoretic applications and the Lickorish criterion are given in Examples~\ref{ex:only-Ob} and~\ref{ex:all-silent} below.

\begin{example}\label{ex:arithmetic-only-ob}
    Let
    \[
        P(t)=5t^4-4t^3-t^2-4t+5,\qquad \Delta=t^{-2}P(t).
    \]
    Then $\Delta(1)=1$ and $\Delta(-1)=17$.
    Set $R=\Z[t^{\pm1}]/(P)$, and denote the image of $t$ in $R$ by $t_E$.
    We first prove that $R$ is integrally closed,
    and then compute $\Ob(R)$ using Proposition~\ref{prop:Ob-integrally-closed-case}. We will show that
    \[
        \cU(R)\cong(\Z/2)^2,\qquad \Ob(R)\cong\Z/2.
    \]
    We will exhibit a unit $q$ of $R_0$ which is a norm from $E$ but not the norm of an $R$-unit,
    whereas $-q$ is not even a norm from $E$.

    Set $E=\Q\otimes R\cong\Q[t]/(P)$ and $u_F=t_E+t_E^{-1}$. With $f(X)=5X^2-4X-11$, we have $f(u_F)=0$, $F=E^\tau=\Q(u_F)\cong\Q[X]/(f)$,
    and we put $R_0:=R^\tau$.
    The relation $f(u_F)=0$ implies that $F=\Q(a)$, where $a=5u_F-2$ satisfies $a^2=59$. In particular,
    $F=\Q(\sqrt{59})$. Since $59\equiv3\bmod4$, the ring of integers is $\cO_F=\Z[\sqrt{59}]=\Z[a]$.
    Both conjugates $(2\pm\sqrt{59})/5$ of $u_F$ lie in $(-2,2)$, so $E/F$ is a CM extension; in particular, every norm from~$E$ is totally positive.
    Moreover, $\Nm_{F/\Q}(u_F-2)\Nm_{F/\Q}(u_F+2)=\frac{f(2)}{5}\cdot\frac{f(-2)}{5}=\frac{17}{25}$ is not a square, so $u_F^2-4$ is not a square in~$F$. Hence $[E:F]=2$, the polynomial $P$ is irreducible, and $E$ is a field.

    PARI (see Appendix~\ref{app:quartic}) shows that
    \[ 1,\ 5t_E,\ 5t_E^2-4t_E,\ 5t_E^3-4t_E^2-t_E\]
    form a $\Z$-basis of $\cO_E$. As all these elements belong to $R$, we infer that $\cO_E\subset R$. By Theorem~\ref{thm:normalization}, the normalization of $R$ is $\wt{R}=\cO_E[t_E,t_E^{-1}]$; since $\cO_E\subset R$ and $t_E^{\pm1}\in R$, we get $\wt R\subset R$. Hence
    \[R=\wt{R}=\cO_{E,T},\]
    so $R$ is integrally closed.
    Since $5t_E$ and $5t_E^{-1}=4-(5t_E^3-4t_E^2-t_E)$ belong to $\cO_E$, all places in $T$ lie above~$5$.
    Modulo~$5$, we have $P\equiv t(t^2+4t+1)$, and $t^2+4t+1$ is irreducible over~$\F_5$, because its discriminant~$2$ is not a square. As $t$ is a unit in $\F_5[t^{\pm1}]$, we get $R/5R\cong\F_{25}$. Hence exactly one prime $\PP_{5,2}$ of $\cO_E$ above~$5$ lies outside~$T$; it has ramification index~$1$ and residue degree~$2$.
    The set~$T$ is nonempty by Lemma~\ref{lem:integral_closure}, as $P$ is not monic, and $\tau$ acts freely on~$T$ by Lemma~\ref{lem:norm_on_T}. Counting degrees, we obtain
    \[
        5\cO_E=\PP_{5,1}\PP_{5,2}\PP_{5,3},\qquad T=\{\PP_{5,1},\PP_{5,3}\},\qquad \tau(\PP_{5,1})=\PP_{5,3},
    \]
    where we choose $\PP_{5,1}$ to be the unique prime above $(5)$
    with $\nu_{\PP_{5,1}}(t_E)>0$.

    Our next goal is to compute $\cU(R)$ explicitly by finding the units of $R_0$. Put $r=5u_F+6=8+a$. Then $\Nm_{F/\Q}(r)=64-59=5$, and the relation $f(u_F)=0$ gives $r(2-u_F)=1$.
    Thus $r$ is a unit of $R_0$, the prime $\pp_5=r\cO_F$ lies below~$T$, and
    \[
        R_0=\cO_{F,S}=\Z[a,r^{-1}]=\Z[u_F],\qquad S=\{\pp_5\}.
    \]
    The units of $\cO_F$ are the elements $x+y\sqrt{59}$ with $x^2-59y^2=\pm1$. The smallest solution with $y>0$ is $(x,y)=(530,69)$ (the sign $-1$ is impossible modulo~$4$), so $530+69a$ is a fundamental unit. Its inverse is
    \[
        q=530-69a=668-345u_F,
    \]
    and both conjugates of $q$ are positive. Since $\pp_5=r\cO_F$ is the only prime of $F$ inverted in $R_0$, we get
    \[R_0^\times=\{\pm q^k r^\ell\}.\]

    Explicit calculations show that $w=5t_E^3+t_E^2-5$ satisfies $\Nm(w)=r$. In particular, $w$ is a unit of $R$, and $[r]=[1]$ in $\cU(R)$, so $\cU(R)$ is generated by $[-1]$ and $[q]$.
    The elements $-1$ and $-q$ are totally negative, so they are not norms from the CM field~$E$; that is, $[-1],[-q]\notin\Ob(R)$. Note that, with
    $z=((561-73a)t_E+114-15a)/31\in E^\times$, we have $q=z\tau(z)$, so $q$ is a norm from $E^\times$, that is, $[q]\in\Ob(R)$.
    Hence $\Ob(R)\subset\{[1],[q]\}$, and it remains to show that $[q]\ne[1]$.

    We use Proposition~\ref{prop:Ob-integrally-closed-case}. Using PARI, we compute that $\Cl(\cO_E)\cong\Z/2$ and $\Cl(\cO_F)=0$. By Lemma~\ref{lem:Cl-localization}, $\Cl(\cO_{E,T})\cong\Cl(\cO_E)/C_T$. PARI shows
    that the two primes in~$T$ are principal; in fact, $\PP_{5,1}=w\cO_E$ and $\PP_{5,3}=(r/w)\cO_E$. So $C_T=0$ and $\Cl(\cO_{E,T})\cong\Z/2$. Since $\tau$ acts trivially on a group of order~$2$, the target of the map $\theta$ in Proposition~\ref{prop:Ob-integrally-closed-case} is $\wH^0(\Z_2;\Cl(\cO_{E,T}))=\Cl(\cO_{E,T})\cong\Z/2$.

    We now study the domain of $\theta$, that is, $\bigoplus_{s=\tau(s)}\Z/2$, where the sum is over the $\tau$-fixed finite places outside~$T$. Such places
    lie above places of $F$ that are inert or ramified in $E$, because $\tau$ interchanges the two places above a split place.

    Suppose $\pp'$ is a place of $F$ that is inert in $E$. There is exactly one place of $E$ above $\pp'$, namely $\PP=\pp'\cO_E$. Its class in $\Cl(E)$
    is $j[\pp']$, where $j\colon\Cl(F)\to\Cl(E)$ is the map induced by $\ma\mapsto\ma\cO_E$. Since $\Cl(F)=0$, the map $j$ is trivial. In particular, $\theta=0$ on all inert places.

    The relative discriminant of $E/F$ is $\md_{E/F}=(31+4a)\cO_F$, an ideal of norm $17$; indeed, PARI gives $d_E=236^2\cdot17=d_F^2\cdot 17$. By Theorem~\ref{thm:discrim}, the prime $\pp_{17}=(31+4a)\cO_F$ is the only prime of $F$ that ramifies in $E$.
    In $R$, the prime above $\pp_{17}$ is generated by $1+t_E$, because $\Nm(1+t_E)=2+u_F=(31+4a)r^{-1}$ and $r$ is a unit of $R_0$. In $\cO_E$, a generator is $(1+t_E)r/w$. In particular,
    this ideal is principal, so $\theta$ also vanishes on the ramified place. Hence $\theta=0$.

    To conclude, Proposition~\ref{prop:Ob-integrally-closed-case} and $\theta=0$ give $\Ob(R)\cong\Z/2$. Since $\Ob(R)\subset\{[1],[q]\}$, we get $\Ob(R)=\{[1],[q]\}$.
    As $[-1]$ and $[-q]$ are not in $\Ob(R)$, the four classes $[\pm1]$, $[\pm q]$ are distinct, and $\cU(R)=\{[\pm1],[\pm q]\}\cong(\Z/2)^2$.
    In particular, $q$ is a norm from $E$ but not the norm of an $R$-unit, whereas $-q$ is not a norm from $E$.
\end{example}

The arguments of the next example are similar to
the arguments of the previous one, nevertheless,
it will reveal a different phenomenon in knot theory.
\begin{example}[A totally real quartic field]\label{ex:arithmetic-all-silent}
    Let
    \[
        P(t)=-3t^4-15t^3+37t^2-15t-3,\qquad \Delta=t^{-2}P(t).
    \]
    We have $\Delta(-1)=61$.
    As in Example~\ref{ex:arithmetic-only-ob}, we first prove that $R$ is integrally closed,
    and then compute $\Ob(R)$ using Proposition~\ref{prop:Ob-integrally-closed-case}.
    In this case, we plan to show that the full norm quotient equals the obstruction group:
    \[
        \cU(R)=\Ob(R)\cong\Z/2.
    \]
    We will exhibit a unit $q$ of $R_0$ such that neither $q$ nor $-q$ is the norm of an $R$-unit,
    although both $q$ and $-q$ are norms from $E$.

    The calculation of $R_0$ is completely analogous to Example~\ref{ex:arithmetic-only-ob}. With
    $E=\Q\otimes R\cong\Q[t]/(P)$, $u_F=t_E+t_E^{-1}$ and $f(X)=3X^2+15X-43$, we have $F=E^\tau=\Q(u_F)\cong\Q[X]/(f)$.

    Next, $F=\Q(a)$ with $a=3u_F+8$ satisfying $a^2-a-185=0$. Then
    $F\cong \Q(\sqrt{741})$, so
    $\cO_F=\Z[\frac12(1+\sqrt{741})]=\Z[a]$.

    Both conjugates $(-15\pm\sqrt{741})/6$ of $u_F$ lie outside $[-2,2]$, so $E$ is totally real.
    Moreover, $P$ is irreducible over $\Q$, and $E$ is
    a field.

    The equation $3u_F^2+15u_F-43=0$ implies that $3(u_F^2+5u_F-14)=1$, so $1/3\in R_0$. PARI (see Subsection~\ref{app:quartic}) shows that 
    \[ 1,\ -3t_E,\ -3t_E^2-15t_E,\ -3t_E^3-15t_E^2+37t_E\]
    form a $\Z$-basis of $\cO_E$. As all these elements belong to $R$, we infer that $\cO_E\subset R$.
    By Theorem~\ref{thm:normalization}, the normalization of $R$ is $\wt{R}=\cO_E[t_E,t_E^{-1}]$; since $\cO_E\subset R$ and $t_E^{\pm1}\in R$, we get $\wt R\subset R$, so $R=\wt R$. Since $3t_E$ and $3t_E^{-1}=(-3t_E^3-15t_E^2+37t_E)-15$ belong to $\cO_E$, and $1/3\in R_0$, we see that
    \[R=\wt{R}=\cO_E[1/3]=\cO_{E,T},\]
    where $T$ consists of the primes of $E$ above $3$.
    Also $R_0=\cO_F[1/3]=\Z[a,1/3]$. This proves that $R$ is integrally closed and determines $R_0$ explicitly, which will be useful in the next step.

    We find the units of $R_0$. First note that the units of $\cO_F$ are of the form $\pm \eta^k$, where $k\in\Z$ and $\eta$ is a
    fundamental unit. Explicitly, we can find $\eta$ as
    $\eta=\frac12(245+9\sqrt{741})=118+9a$.

    To pass to the units of $R_0=\cO_F[1/3]$, we look for an element of norm $\pm3$. In the same step, we determine the sets $T$ and $S$.
    Consider the ideal $3\cO_F$. Since $a^2-a-185\equiv (a+1)^2\bmod 3$, we have $3\cO_F=\pp_3^2$, where $\pp_3=(3,a+1)$. That is, the prime $3$ ramifies in $F$. Hence $S=\{\pp_3\}$, and since $\tau$ acts freely on $T$ (Lemma~\ref{lem:norm_on_T}), the prime $\pp_3$ splits in $E$: $T=\{\PP_3^+,\PP_3^-\}$ with $\tau(\PP_3^+)=\PP_3^-$. Consider
    $q=14-a$. Then $\Nm_{F/\Q}(q)=14^2-14-185=-3$ and $q=15-(a+1)\in\pp_3$, so $\pp_3=q\cO_F$. In particular, $q$ is a unit in $R_0$ (indeed, $q^{-1}=-u_F-7$). Moreover, since $\pp_3$ is the only prime of $F$ inverted in $R_0$, the element $q$ together with the
    units of $\cO_F$ generates all units of $R_0$; hence
    \[R_0^\times=\{\pm q^k\eta^\ell\}.\]

    Using PARI, we find
    \[ b = -3t_E^3 - 18 t_E^2 + 19t_E+4,\qquad w = -15t_E^3-90t_E^2+98t_E+40\]
    with the property that
    \[\Nm(b)=q\eta,\qquad \Nm(w)=-1.\]
    In particular, $b$ and $w$ are units of $R$, because their norms are units of $R_0$.
    This shows that the classes $[q\eta]$ and $[-1]$ in $\cU(R)$ are trivial. As $[q]^2=[\Nm(q)]=[1]$, we get $[q]=[\eta]$, so $\cU(R)$ is generated by $[q]$. At this point,
    we cannot claim that $\cU(R)$ is nontrivial, because $[q]$ might be zero. Note that, with
    $z=((149 + 14a)t_E - 335 - 6a)/13\in E^\times$, we have $q=z\tau(z)$, so $q$ is a norm from $E^\times$, that is, $[q]\in\Ob(R)$.

    Our next step is to show that $\Ob(R)\cong\Z/2$, so that $\cU(R)=\Ob(R)\cong\Z/2$ is generated by $[q]$. We will use Proposition~\ref{prop:Ob-integrally-closed-case}. Using PARI, we compute that $\Cl(\cO_E)\cong\Cl(\cO_F)\cong\Z/2$, and that $\Cl(\cO_F)$ is generated by the class of the ideal $\pp=(5,a)$. By Lemma~\ref{lem:Cl-localization}, $\Cl(\cO_{E,T})\cong\Cl(\cO_E)/C_T$. PARI shows
    that the two primes above $\pp_3$ are principal; in fact, $\PP_3^+=b\cO_E$, where $\PP_3^+$ is the prime with $\nu_{\PP_3^+}(t_E)<0$, and $\PP_3^-=(q/b)\cO_E$. So $C_T=0$ and $\Cl(\cO_{E,T})\cong\Z/2$. Since $\tau$ acts trivially on a group of order~$2$, the target of the map $\theta$ in Proposition~\ref{prop:Ob-integrally-closed-case} is $\wH^0(\Z_2;\Cl(\cO_{E,T}))=\Cl(\cO_{E,T})\cong\Z/2$.

    We now study the domain of $\theta$. Again, this boils down
    to checking inert and ramified places in $E$. Unlike
    in Example~\ref{ex:arithmetic-only-ob}, $\Cl(F)$ is non-trivial.

    Take a place $\pp'$ of $F$ that is inert in $E$. There is exactly one place of $E$ above $\pp'$, namely $\PP=\pp'\cO_E$. Its class in $\Cl(E)$
    is $j[\pp']$, where $j\colon\Cl(F)\to\Cl(E)$ is the map induced by $\ma\mapsto\ma\cO_E$.
    Now the ideal $\pp\cO_E$ is principal: one can check that $(15 + a) - (10 + a)t_E = -3t_E^3 - 18t_E^2 + 22t_E + 5$ is a generator,
    and $\pp$ is inert in $E$. Since $[\pp]$ generates $\Cl(F)$, the map $j$ is trivial. In particular, $\theta=0$ on all inert places.

    The ramified places, via the calculation of $\md_{E/F}=(53+4a)\cO_F$, are places over $61$. PARI shows that $(53+4a)\cO_F$ is the
    only prime of $F$ ramifying in $E$. 
    In $R$, the prime above it is generated by $1+t_E$, because $\Nm(1+t_E)=2+u_F=-q(53+4a)/3$, and $q$ and $3$ are units of $R_0$. In $\cO_E$, a generator is $(1+t_E)b$. In particular,
    this ideal is principal, so $\theta$ also vanishes on the ramified place. Hence $\theta=0$.

    To conclude, Proposition~\ref{prop:Ob-integrally-closed-case} and $\theta=0$ give $\Ob(R)\cong\Z/2$. Since $\Ob(R)\subset\cU(R)$ and $\cU(R)$ is generated by $[q]$, we get $\cU(R)=\Ob(R)=\{[1],[q]\}\cong\Z/2$,
    and $[-q]=[q]$, because $[-1]=[1]$. In particular, neither $q$ nor $-q$ is the norm of an $R$-unit, while both are norms from $E$: $q=\Nm(z)$ and $-q=\Nm(wz)$.
\end{example}

\subsection{Products of polynomials}\label{sub:product-computations}

We now allow $\Delta$ to be reducible. Lemma~\ref{lem:reduction} and Remark~\ref{rem:redQ} allow us to replace it by its radical without changing $\cU(R)$ or $\Ob(R)$, therefore we assume that $\Delta$ is reduced.

Write $\Delta=\prod_{i=1}^r\Delta_i$, where the factors are pairwise coprime over $\Q$, symmetric, and normalized by $\Delta_i(1)=1$, and put
\[
    R=\Lambda/(\Delta),\qquad R_i=\Lambda/(\Delta_i).
\]
The symmetry assumption ensures that each $R_i$ carries the induced involution; irreducible factors exchanged by $\tau$ must be grouped together.

\begin{proposition}\label{prop:product-norm-groups}
    If $(\Delta_i,\Delta_j)=\Lambda$ whenever $i\ne j$, then the Chinese remainder isomorphism induces
    \[
        \cU(R)\cong\bigoplus_{i=1}^r\cU(R_i),\qquad
        \Ob(R)\cong\bigoplus_{i=1}^r\Ob(R_i).
    \]
\end{proposition}
\begin{proof}
    The isomorphism $R\cong\prod_iR_i$ respects the involution, so fixed units and norms are computed componentwise.
    The same decomposition holds after tensoring with $\Q$; taking kernels of the maps to the rational norm quotients proves the assertion about $\Ob$.
\end{proof}

\section{Applications in knot theory}\label{sec:knot}

\subsection{Nakanishi index}\label{sub:nakanishi}

Put $\Lambda=\Z[t^{\pm1}]$, and let $A$ be an $n\times n$ matrix over~$\Lambda$ with nonzero determinant~$\Delta$.
Set $M=\coker A=\Lambda^n/A\Lambda^n$ and $R=\Lambda/(\Delta)$.
The adjugate identity $A\adj(A)=\Delta\Id$ shows that $\Delta M=0$, so $M$ is naturally an $R$-module.
If $\Delta$ is a unit, then $A$ is invertible and $M=R=0$.
\begin{question}\label{qu:MR}
  Is $M$ isomorphic to $R$?
\end{question}
\begin{example}
  If $V$ is a Seifert matrix of a knot~$K$, then $A=V-tV^T$ presents its Alexander module; see~\cite{BF1}.
  Its determinant differs from the normalized symmetric Alexander polynomial of Subsection~\ref{sub:rank-one-blanchfield} by a unit of~$\Lambda$,
  so both define the same quotient ring~$R$.
  As explained in Subsection~\ref{sub:rank-one-blanchfield}, a nonzero cyclic Alexander module is isomorphic to~$R$.
  Thus, when the Alexander polynomial is not~$1$, we have $M\cong R$ if and only if $m(K)=1$.
  If the Alexander polynomial is~$1$, then $M=R=0$ and $m(K)=0$.
\end{example}
For the rest of this subsection, assume that $\Delta$ is not a unit.
Kearton and Wilson~\cite{KeartonWilson} study Question~\ref{qu:MR} using ideal-class methods originating in the work of
Fox and Smythe~\cite{FoxSmythe}. We formulate the obstruction in terms of the Picard group and explain the presentation convention below.

We use Fitting ideals, see~\cite[Section 20.2]{Eisenbud}. For a finitely presented $R'$-module~$N$ with a presentation
$(R')^m\xrightarrow{B}(R')^n\to N\to0$, the ideal $\Fitt_i^{R'}(N)$ is generated by the $(n-i)\times(n-i)$ minors of~$B$
for $0\le i<n$, and $\Fitt_i^{R'}(N)=R'$ for $i\ge n$.
These ideals are independent of the presentation. For our module~$M$, we abbreviate $\Fitt_i^\Lambda(M)$ to~$\Fitt_i(M)$;
in particular, $\Fitt_0(M)=(\Delta)$. We use two basic properties, the first of which follows from
\cite[Corollary-Definition 20.4]{Eisenbud}.

\begin{lemma}\label{lem:fit_gen}
  If a finitely presented $R'$-module~$N$ can be generated by $k$ elements, then $\Fitt_i^{R'}(N)=R'$ for $i\ge k$.
\end{lemma}
The next result describes base change for Fitting ideals; see~\cite[Corollary 20.5]{Eisenbud}.
\begin{lemma}\label{lem:fit_ten}
  If $R'\to R''$ is a ring map and $N$ is a finitely presented $R'$-module, then
  \[
      \Fitt_i^{R''}(N\otimes_{R'}R'')=\Fitt_i^{R'}(N)R'',
  \]
  where the right-hand side is the ideal generated by the image of $\Fitt_i^{R'}(N)$ in~$R''$.
\end{lemma}
\begin{proof}
  Tensoring a presentation with~$R''$ gives a presentation of $N\otimes_{R'}R''$.
  Its minors are the images of the original minors, which proves the formula.
\end{proof}
In particular, Fitting ideals localize. 

If $M\cong R$, then the matrix $[\Delta]$ presents $M$ over~$\Lambda$, so $\Fitt_1(M)=\Lambda$.
This gives the first obstruction to Nakanishi index~$1$.
\begin{corollary}\label{cor:first_nak}
  If $\Fitt_1(M)\ne\Lambda$, then $M\not\cong R$.
\end{corollary}
From now on, assume that $\Fitt_1(M)=\Lambda$.

\begin{proposition}
  Suppose $\Fitt_1(M)=\Lambda$. Then $M$ is an invertible $R$-module.
\end{proposition}
\begin{proof}
  Since $\Delta M=0$, we have $M\otimes_\Lambda R\cong M$.
  By Lemma~\ref{lem:fit_ten}, $\Fitt_0^R(M)=(\Delta)R=0$ and $\Fitt_1^R(M)=\Fitt_1(M)R=R$.
  By \cite[Proposition 20.8]{Eisenbud}, $M$ is a projective $R$-module of constant rank $1$, hence invertible.
\end{proof}
An invertible module defines an element of the Picard group~$\Pic(R)$.
\begin{definition}
  Suppose $\Fitt_1(M)=\Lambda$. The \emph{Fox--Smythe class} of $M$ is its isomorphism class
  $\FS(M)=[M]\in\Pic(R)$ as an invertible $R$-module.
\end{definition}
The following criterion is immediate from the definition of the Picard group.
\begin{corollary}\label{cor:second_nak}
  Suppose $\Fitt_1(M)=\Lambda$. Then $M\cong R$ if and only if $\FS(M)=0\in\Pic(R)$.
\end{corollary}
We next represent $\FS(M)$ by an ideal and compare it with the classical \emph{row class}.

To represent $M$ as a fractional ideal, we use the total quotient ring of~$R$. This notion is a generalization
of a field of fractions, so it is applicable if $R$ is not a domain, i.e.\ if $\Delta$ is not irreducible.

Let $\cS$ be the set of non-zero-divisors of~$R$, and put $K_R=\cS^{-1}R$.
Write $\Delta=u\prod_{i=1}^r f_i^{e_i}$, where $u\in\Lambda^\times$ and $f_i$ are distinct irreducible factors.
The minimal primes of~$R$ are $\mfp_i=(f_i)/(\Delta)$.
\begin{lemma}\label{lem:m_in_KR}
  The ring $K_R$ is Artinian, with decomposition
  \[
      K_R\cong\prod_{i=1}^r R_{\mfp_i}
      \cong\prod_{i=1}^r\Lambda_{(f_i)}/(f_i^{e_i}).
  \]
  Under the assumption $\Fitt_1(M)=\Lambda$, we have $M\otimes_RK_R\cong K_R$.
  The localization map $M\to M\otimes_RK_R$ is injective, so a choice of this isomorphism gives an embedding $M\hookrightarrow K_R$.
\end{lemma}
\begin{proof}
  Unique factorization in~$\Lambda$ shows that the zero-divisors of~$R$ form $\bigcup_i\mfp_i$.
  By prime avoidance, the primes of~$R$ disjoint from~$\cS$ are precisely the~$\mfp_i$.
  Hence $K_R$ is Noetherian of dimension~$0$, and therefore Artinian. Its decomposition into local factors gives the desired isomorphisms.

  The preceding proposition makes $M\otimes_RK_R$ locally free of rank~$1$.
  Over each local factor it is free of rank~$1$, so it is free of rank~$1$ over~$K_R$.
  Finally, $M$ is projective and hence flat. Tensoring the injection $R\hookrightarrow K_R$ with~$M$ proves injectivity of the localization map.
\end{proof}
The embedding of~$M$ depends on the chosen trivialization of $M\otimes_RK_R$; two choices differ by multiplication by an element of~$K_R^\times$.
\begin{corollary}\label{cor:frac}
  There exists a row vector $(\psi_1,\dots,\psi_n)$ of elements in $K_R$
  such that
  \[M\cong \psi_1R+\dots+\psi_nR.\]
\end{corollary}
\begin{proof}
  Choose an embedding $M\hookrightarrow K_R$ as in Lemma~\ref{lem:m_in_KR}.
  Let $\psi_1,\dots,\psi_n$ be the images of the standard generators under $R^n\to M\hookrightarrow K_R$.
  They generate the image ideal, which is isomorphic to~$M$.
\end{proof}
Corollary~\ref{cor:frac} represents $\FS(M)$ by an invertible fractional ideal in~$K_R$.
Changing the trivialization in Lemma~\ref{lem:m_in_KR} multiplies this ideal by an element of~$K_R^\times$, preserving its class.
We now describe the maps that produce such ideals from the presentation matrix.

Write $A^{\ad}=\adj(A)$ for the adjugate matrix, so
$AA^{\ad}=A^{\ad}A=\Delta\Id$. The reductions of its rows modulo~$\Delta$ satisfy $\psi A=0$.
Every such row defines an $R$-linear map
\[
    \phi_\psi\colon M\longrightarrow R,\qquad [v]\longmapsto\psi v,
\]
whose image is $I_\psi=\psi_1R+\cdots+\psi_nR$.
Thus $\psi$ specifies an element of $\Hom_R(M,R)$.
Conversely, every row satisfying $\psi A=0$ is an $R$-linear combination of the reduced rows of~$A^{\ad}$.
Indeed, a lift $\widehat\psi\in\Lambda^n$ satisfies $\widehat\psi A=\Delta\eta$ for some $\eta\in\Lambda^n$.
Multiplying by~$A^{\ad}$ and cancelling the nonzero element~$\Delta$ gives $\widehat\psi=\eta A^{\ad}$.

\begin{lemma}\label{lem:adjugate_row}
  Keep the assumption $\Fitt_1(M)=\Lambda$, and let $\psi\in R^n$ satisfy $\psi A=0$.
  The map $\phi_\psi$ induces an isomorphism $M\cong I_\psi$ if and only if
  $\psi\bmod\mfp_i\ne0$ in $(R/\mfp_i)^n$ for every minimal prime~$\mfp_i$ of~$R$.
  When these conditions hold, $I_\psi$ is an invertible ideal representing~$\FS(M)$.
\end{lemma}
\begin{proof}
  By Lemma~\ref{lem:m_in_KR}, $K_R=\prod_i R_{\mfp_i}$ and $M\otimes_RK_R$ is free of rank~$1$.
  At each~$\mfp_i$, choose a generator of~$M_{\mfp_i}$; then $(\phi_\psi)_{\mfp_i}$ is multiplication by some $c_i\in R_{\mfp_i}$.
  The condition $\psi\bmod\mfp_i\ne0$ says that the induced map over the residue field is nonzero, since
  $R_{\mfp_i}^n\to M_{\mfp_i}$ is surjective. It is therefore equivalent to $c_i$ being a unit.

  If all these conditions hold, $\phi_\psi\otimes_RK_R$ is an isomorphism.
  The injectivity of $M\to M\otimes_RK_R$ then implies that $\phi_\psi$ is injective, hence an isomorphism onto~$I_\psi$.
  Conversely, if $\phi_\psi$ is injective, its localization at each~$\mfp_i$ is injective.
  An injective endomorphism of the finite-length module $R_{\mfp_i}$ is surjective, so every~$c_i$ is a unit.
  Finally, $I_\psi\cong M$ is locally free of rank~$1$, and its localization is all of~$K_R$;
  thus it is an invertible ideal with class~$\FS(M)$.
\end{proof}

For comparison with the classical Fox--Smythe construction~\cite{FoxSmythe}, suppose first that $\Delta$ is irreducible, so $R$ is a domain.
In the row-presentation convention, the relation matrix for our module is $A'=A^T\bmod\Delta$.
The row class is represented by the ideal of maximal minors from $n-1$ linearly independent rows of~$A'$;
see~\cite[Section 3.5, Theorem 3.17]{HillmanLinks}.
Deleting row~$i$ of~$A'$ gives, up to signs, the entries of row~$i$ of~$A^{\ad}\bmod\Delta$.
Thus any nonzero such row represents $\FS(M)$ by Lemma~\ref{lem:adjugate_row}.
For reducible or nonreduced~$\Delta$, the definition $\FS(M)=[M]\in\Pic(R)$ 
extends this description to invertible modules over~$R$.

We now explain how the conductor square of Subsection~\ref{sub:obcomp} answers
Question~\ref{qu:MR}. Assume additionally that $\Delta(1)=\pm1$, as for an Alexander polynomial.
First, repeated factors can be removed: if $\ma$ is the nilradical of $R$, then
\[
    M\cong R\quad\Longleftrightarrow\quad M/\ma M\cong R/\ma.
\]
Indeed, a lift $x\in M$ of a generator modulo $\ma$ generates $M$.
Since $M / Rx=\ma (M / Rx)$ and $\ma$ is nilpotent, Nakayama's Lemma implies that \(M / Rx = 0\).
The resulting surjection $R\to M$ is an isomorphism because $M$ is locally free of rank~$1$.
Thus we may replace $(R,M)$ by $(R/\ma,M/\ma M)$ below.
For the norm problem, the corresponding reduction is Lemma~\ref{lem:reduction}.

Write $B=\wt{R}=\prod_f\cO_{E_f,T_f}$ for the normalization of the reduced ring,
and let $\mf$ be its conductor. Put
\[
    C=R/\mf,\qquad D=B/\mf,\qquad
    \Xi=D^\times/C^\times,\qquad
    \Xi_0=\im\left(\rho:B^\times\longrightarrow\Xi\right),
\]
where $\rho$ is reduction modulo $\mf$ followed by passage to the quotient.
Since $C$ and $D$ are finite, their Picard groups vanish.
The Units--Pic sequence \eqref{eq:unit_pic} of the conductor square, see~\cite[Theorem I.3.10]{Weibel_K}, therefore gives an exact sequence
\begin{equation}
    \label{eq:conductor_pic}
    0 \longrightarrow
    \Xi / \Xi_0 \longrightarrow
    \Pic(R) \longrightarrow
    \Pic(B) \longrightarrow
    0.
\end{equation}

Sequence~\eqref{eq:conductor_pic} gives a two-stage test for $[M]=0$.
First test whether $M\otimes_RB$ is free; this is a principal-ideal problem in each
$\cO_{E_f,T_f}$, whose class group is
\[
    \Cl(\cO_{E_f,T_f})\cong
    \Cl(\cO_{E_f})/\langle[\mathfrak p]:\mathfrak p\in T_f\rangle.
\]
If it is free, choose generators $a$ of $M\otimes_RB$ and $b$ of $M\otimes_RC$, where the latter module is free because $\Pic(C)=0$.
In $M \otimes_{R} D$ they satisfy $b=\lambda a$ for a unique $\lambda\in D^\times$.
Changing either generator changes $\lambda$ by the residue of a unit of $B$ or $C$.
Milnor patching therefore gives
\[
    M\cong R\quad\Longleftrightarrow\quad
    [\lambda]=0\text{ in }\Xi/\Xi_0.
\]
This is a membership test in a finite group. When it succeeds, the chosen generators
can be multiplied by units so that their images in $M\otimes_RD$ agree; they then
patch to a generator of $M$. Thus the procedure also constructs the isomorphism
required in Step~3 of Subsection~\ref{sub:algorithm}.

\subsection{Real and rational unknotting number one}\label{prop:qunkn}
We consider the real and rational variants of the algebraic unknotting number, with base rings $\R$ and $\Q$, respectively; see Definition~\ref{def:signed_S}.
Throughout this subsection, let $K$ have Nakanishi index $m(K)=1$, and let $\Delta$ be its symmetrized Alexander polynomial,
normalized by $\Delta(t)=\Delta(t^{-1})$ and $\Delta(1)=1$.
Put $R=\Lambda/(\Delta)$, fix an isomorphism $M(K)\cong R$, and write
\[
    g=\wt{\Bl}(1,1)\in(R^\tau)^\times.
\]
Here $\tau$ is induced by $t\mapsto t^{-1}$.
For $S\in\{\Q,\R\}$, write
\[
    R_S=R\otimes_\Z S\cong S[t^{\pm1}]/(\Delta),
    \qquad \Nm(x)=x\tau(x),
\]
and use $g$ also for its image in~$R_S$.
Use the symmetric field factors $E_i/F_i$ of Subsection~\ref{sub:reducible-qnorm}, with $p=\Delta$,
and let $g_i\in F_i^\times$ be the image of~$g$ in the corresponding fixed field.
We begin with the real case. We connect the real norm criterion to the infinite-place tests of Section~\ref{sec:qnorm}.

\begin{proposition}\label{prop:real_unknotting}
  Let $\epsilon=\pm 1$, and use the setup above.
  Then $\epsilon g\in\Nm(R_{\R}^{\times})$ if and only if, for every~$i$, $\epsilon g_i$ is a local norm
  from $E_i/F_i$ at every infinite place of~$F_i$.
\end{proposition}
\begin{proof}
  By the argument of Lemma~\ref{lem:qnorm} over~$\R$ instead of~$\Q$, the norm quotient of~$R_\R$ equals that of its reduction, which is a product of copies of~$\R$ and~$\C$.
  Factors exchanged by~$\tau$ contribute nothing by Lemma~\ref{lem:shapiro}.
  The remaining factors are copies of~$\C$ on which $\tau$ is complex conjugation; they correspond to the pairs $\{\xi,\ol\xi\}$ of roots of~$\Delta$ on the unit circle, and the norm there is $z\mapsto|z|^2$.
  Hence $\epsilon g\in\Nm(R_\R^\times)$ if and only if $\epsilon g(\xi)>0$ for every root $\xi\in S^1$ of~$\Delta$, which is the infinite-place condition by Proposition~\ref{prop:real_places} applied to each $E_i/F_i$.
\end{proof}

The next proposition specializes~\cite[Theorem~1.1]{BF2} to our signed rank-one setting, using the norm description below.
Recall that for a Seifert matrix~$V$
\[
    \sigma_K(z)=\sign\left((1-z)V+(1-\overline z)V^T\right),\qquad z\in S^1.
\]

\begin{proposition}\label{prop:real_number}
  Let $\epsilon=\pm 1$, and use the setup above. The following conditions are equivalent.
  \begin{itemize}
    \item $\sigma_K(z)\in\{-2\epsilon,-\epsilon,0\}$ for every $z\in S^1$;
    \item $\epsilon g\in\Nm(R_{\R}^{\times})$;
    \item for any root $\xi$ of $\Delta$ such that $|\xi|=1$, we have $\epsilon g(\xi)>0$;
    \item \(u_{\R}^{\epsilon}(K) = 1\).
  \end{itemize}
\end{proposition}
\begin{proof}
  The statement is proved implicitly in \cite{BF2}. Proposition~\ref{prop:real_unknotting}
  provides the equivalence of the second and the third conditions.
\end{proof}

We now turn to the rational case, retaining the setup above. For $\epsilon=\pm1$, we ask whether $\epsilon g\in\Nm(R_{\Q}^{\times})$.
The exact local criterion is the following specialization of Proposition~\ref{prop:reducible-norm-tests}(1).
\begin{proposition}\label{prop:rational_unk}
  Let $\epsilon=\pm 1$. Then $\epsilon g\in\Nm(R_{\Q}^{\times})$, that is, \(u_{\Q}^{\epsilon}(K) = 1\), if and only if both conditions hold:
  \begin{enumerate}
    \item $\epsilon g(\xi)>0$ for every root $\xi\in S^1$ of~$\Delta$;
    \item for every~$i$ and every finite place~$s$ of~$F_i$ ramified in~$E_i/F_i$,
    \[
        (\epsilon g_i,u_i^2-4)_s=1,\qquad u_i=t_{E_i}+t_{E_i}^{-1}.
    \]
  \end{enumerate}
  Here $t_{E_i}$ is the image of~$t$ in~$E_i$, and the second condition says that $\epsilon g_i$ is a local norm at~$s$.
\end{proposition}
\begin{proof}
  Apply Proposition~\ref{prop:reducible-norm-tests}(1) to the invariant unit $b=\epsilon g$ of~$R$.
\end{proof}

The Lickorish criterion gives a sufficient condition for rational unknotting number one once the real condition is satisfied.
Put $d=\Delta(-1)$ and $D=|d|=\det(K)$, so $d$ is the signed determinant and $D$ is the order of $H_1(\Sigma(K);\Z)$,
where $\Sigma(K)$ is the double branched cover.
We use the linking convention of~\cite[Section~1.1]{BF1}: for a Seifert matrix~$V$, the linking pairing on
$\coker(V+V^T)$ is
\[
    \lk([v],[w])=v^T(V+V^T)^{-1}w\pmod\Z.
\]
In particular, $\lk$ is completely determined by a Seifert matrix.
We recall Lickorish's theorem~\cite{Lickorish} in the convention of~\cite[Theorem~4.5]{BF1}.
\begin{theorem}[Lickorish] \label{thm:lickorish}
  Suppose $K$ can be unknotted by a single crossing change of sign $\delta\in\{\pm1\}$.
  Then there exists a generator $h\in H_1(\Sigma(K);\Z)$ such that $\lk(h,h)=-2\delta/d$.
\end{theorem}
\begin{remark}
  A similar criterion was considered by Jabuka \cite{Jabuka}. In \cite[Section 4]{BF1} the two criteria are compared. It turns
  out that Jabuka's criterion follows from the Lickorish criterion combined with signature obstructions.
\end{remark}
For $\delta=-\epsilon$, the condition of Theorem~\ref{thm:lickorish} reads $\lk(h,h)=2\epsilon/d$.
The following result compares this condition with the square test; it was proved in~\cite[Lemma~3.3 and Section~4.5]{BF1}.
\begin{proposition}\label{prop:square}
Let $\epsilon=\pm 1$, and use $K$, $R$, and $g$ from the setup above. The following conditions are equivalent.
\begin{itemize}
  \item there exists a generator $h\in H_1(\Sigma(K);\Z)$ such that $\lk(h,h)=2\epsilon/d$;
  \item $\epsilon g(-1)$ is a square modulo~$D$.
\end{itemize}
\end{proposition}

The following result combines Proposition~\ref{prop:square} with Proposition~\ref{prop:reducible-norm-tests}(2).
\begin{theorem}\label{thm:lick_rat}
  Let $\epsilon=\pm1$, and suppose $u_{\R}^{\epsilon}(K)=1$.
  If there exists a generator $h\in H_1(\Sigma(K);\Z)$ such that
  $\lk(h,h)=2\epsilon/d$, then $u_{\Q}^{\epsilon}(K)=1$.
\end{theorem}
\begin{proof}
  Put $b=\epsilon g\in(R^\tau)^\times$.
  Proposition~\ref{prop:real_number} gives $b(\xi)>0$ at every unit-circle root of~$\Delta$,
  while Proposition~\ref{prop:square} says that the image of~$b$ in $R/(t+1)R$ is a square.
  Proposition~\ref{prop:reducible-norm-tests}(2) therefore gives $b\in\Nm(R_{\Q}^{\times})$.
  By Definition~\ref{def:signed_S}, this is $u_{\Q}^\epsilon(K)=1$.
\end{proof}

\subsection{Algebraic unknotting number}\label{sub:algorithm}
The criteria of Sections~\ref{sec:qnorm} and~\ref{sec:zorm} 
provide a four-step algorithm for determining whether $u_a(K)=1$.

Let $\Delta$ be the symmetrized Alexander polynomial, normalized by $\Delta(t)=\Delta(t^{-1})$ and $\Delta(1)=1$.
If $\Delta=1$, then $u_a(K)=0$, and no further tests are needed.
Henceforth assume $\Delta\ne1$.

\emph{Step 1.} Check whether the Alexander module is cyclic, using the procedure described at the end of Section~\ref{sub:nakanishi}.
If it is not, then $u_a(K)\ge m(K)>1$ and nothing more is necessary.
Otherwise retain the generator constructed in Subsection~\ref{sub:nakanishi},
which gives an explicit isomorphism $M(K)\cong R_\Delta:=\Lambda/(\Delta)$.

We use the notation of Definition~\ref{def:signed_S}, that is,
we consider $\epsilon=\pm 1$. By Theorem~\ref{thm:fogel},
$u_\Z^\epsilon(K)=1$ corresponds to a crossing change of sign $-\epsilon$. We will apply Steps 2--4 for both signs $\epsilon$.
If one choice of $\epsilon$ passes all three steps, then $u_a(K)=1$;
otherwise $u_a(K)>1$.

\emph{Step 2.} By Proposition~\ref{prop:real_number}, $u_{\R}^\epsilon(K)=1$ exactly when
\[
    \sigma_K(z)\in\{0,-\epsilon,-2\epsilon\}\qquad\text{for every }z\in S^1.
\]
This is a finite test: isolate the unit-circle roots of $\Delta$, compute the
signature at each root, and choose one algebraic sample point in each component
of the complement of these roots and $1$ in $S^1$.
The signature is constant on each such component, and $\sigma_K(1)=0$.
All these signatures, including those of the singular matrices at the roots,
can be computed using exact real-algebraic arithmetic.
If any value is outside the given set, discard this sign~$\epsilon$.

\emph{Step 3.} 
Using the generator retained in Step~1, define
$g=\wt{\Bl}(1,1)\in(R_\Delta^\tau)^\times$.
For each surviving sign, test whether $\epsilon g(-1)$ is a square
modulo~$|\Delta(-1)|$.
If it is not, discard this sign.
If the square test passes, Step~2 and Theorem~\ref{thm:lick_rat}, together with
Proposition~\ref{prop:square}, show that $\epsilon g$ is a norm in
$\Q\otimes_\Z R_\Delta$.

\emph{Step 4.} Determine if $\epsilon g\in\Nm(R_\Delta^\times)$.
Let $p=\Delta_{\mathrm{red}}$ be the product of the distinct irreducible factors of~$\Delta$, chosen up to a unit, and put
\[
    R_{\mathrm{red}}=\Lambda/(p).
\]
Write $g_{\mathrm{red}}\in(R_{\mathrm{red}}^\tau)^\times$ for the image of~$g$ under the quotient map $R_\Delta\to R_{\mathrm{red}}$.
Corollary~\ref{cor:is_norm} shows that $\epsilon g$ is a norm
in $R_\Delta$ if and only if $\epsilon g_{\mathrm{red}}$ is a norm in $R_{\mathrm{red}}$.
By Steps~2 and~3, $\epsilon g$ is a norm in $\Q\otimes_\Z R_\Delta$, so $\epsilon g_{\mathrm{red}}$ defines a class

\[
    c_\epsilon:=[\epsilon g_{\mathrm{red}}]\in\Ob(R_{\mathrm{red}}).
\]

The remaining test is whether this particular class vanishes.
Use the unit-coordinate computation of Remark~\ref{rem:algorithm-effectivity}
to test whether $\epsilon g_{\mathrm{red}}$ is a norm in $R_{\mathrm{red}}$.
This decides whether $c_\epsilon=0$ directly, without first computing the whole
group $\Ob(R_{\mathrm{red}})$.
By Corollary~\ref{cor:is_norm} and Theorem~\ref{thm:fogel},
$c_\epsilon=0$ if and only if $u_\Z^\epsilon(K)=1$.
Section~\ref{sec:examples} illustrates the obstruction calculations in concrete examples.

\begin{remark}\label{rem:algorithm-effectivity}
    Computations in Steps~1 and~4 are effective.
    Throughout this remark, set $R=R_{\mathrm{red}}$ and $B=\wt{R}$.
    First observe that Step~1 requires verification of the condition $1\in\Fitt_1(M(K))$.
    Verifiation of this condition is effective over $\Lambda\cong\Z[t,s]/(ts-1)$ by integer polynomial ideal-membership algorithms     \cite{Aschenbrenner2004IdealMembership}.

    The unit to be tested in Step~4 is $\epsilon g_{\mathrm{red}}$.
    Integral bases for the number-field factors give finite $R$-module generators of $B$; clearing denominators gives $NB\subset R$
    for some integer $N>0$, as in Lemma~\ref{lem:finite_index}.
    In the finite ring $B/NB$, the conductor modulo $NB$ consists of the elements $d$ with $d(B/NB)\subset R/NB$, so it is found by finite enumeration.
    Number-field unit and ideal-class algorithms then supply the generators, relations, and principal generators needed above, see~\cite[Sections~4--6]{Lenstra1992Algorithms}.

    Sequence~\eqref{eq:conductor-units} computes $R^{\times}$ as the kernel of a map from the finitely generated group $B^\times$ to a finite group.
    The involution $t\mapsto t^{-1}$ preserves this sequence.
    Computing its action on generators of $R^\times$ gives
    \[
        \cU(R)=\ker(1-\tau:R^\times\to R^\times)/\im(1+\tau:R^\times\to R^\times),
    \]
    where $1-\tau$ and $1+\tau$ denote $u\mapsto u/\tau(u)$ and $u\mapsto u\tau(u)$, respectively.
    In particular, expressing a specified invariant unit in these generators reduces norm membership to integer linear algebra, with torsion relations included.
    This gives a direct computation of its class, complementary to the exact sequences of Subsection~\ref{sub:obcomp}.
\end{remark}

\subsection{Connected sums}\label{sec:connected}
We apply the preceding criteria to $K=K_1\#K_2$. Throughout this subsection, assume $m(K_1)=m(K_2)=1$.
Let $\Delta_i$ be the symmetrized Alexander polynomial of~$K_i$, normalized by
$\Delta_i(t)=\Delta_i(t^{-1})$ and $\Delta_i(1)=1$, and put $\Delta=\Delta_1\Delta_2$.
Write $M_i=M(K_i)$ and $R_i=\Lambda/(\Delta_i)$, so $M_i\cong R_i$, and put $R=\Lambda/(\Delta)$.
We apply the product computations of Subsection~\ref{sub:product-computations}: first determine when $M(K)$ is cyclic, then compute the coefficient of its Blanchfield pairing.

\begin{proposition}\label{prop:nakanishi_connected}
  The Nakanishi index of $K$ is $1$ if and only if $(\Delta_1,\Delta_2)=\Lambda$.
\end{proposition}
\begin{proof}
  The Alexander module of a connected sum is the direct sum of the Alexander modules of its summands. Hence
  \[
      M(K)\cong R_1\oplus R_2
      \cong\coker\begin{pmatrix}\Delta_1&0\\0&\Delta_2\end{pmatrix}.
  \]
  In the notation of Subsection~\ref{sub:nakanishi}, $\Fitt_1(M(K))=(\Delta_1,\Delta_2)$.
  If $M(K)$ is cyclic, Lemma~\ref{lem:fit_gen} gives $\Fitt_1(M(K))=\Lambda$.
  Conversely, if $(\Delta_1,\Delta_2)=\Lambda$, the Chinese remainder theorem gives
  $R\cong R_1\times R_2$, so $M(K)\cong R$ is cyclic.
  It is nonzero because $M_1$ and $M_2$ are nonzero, and therefore $m(K)=1$.
\end{proof}

For the rest of this subsection, assume $(\Delta_1,\Delta_2)=\Lambda$.
We use Proposition~\ref{prop:product-norm-groups}, whose Chinese remainder decomposition also holds after extension of scalars.
For $S\in\{\Z,\Q,\R\}$, put
\[
    R_S=S[t^{\pm1}]/(\Delta),\qquad R_{i,S}=S[t^{\pm1}]/(\Delta_i).
\]
The same componentwise norm calculation gives
\begin{equation}\label{eq:additivity_of_U}
    \cU(R_S)\cong\cU(R_{1,S})\oplus\cU(R_{2,S}).
\end{equation}
In particular, extension of scalars from~$\Z$ to~$\Q$ respects this decomposition.

To apply the decompositions, we must track the coefficient of the Blanchfield pairing.
Choose generators $x_i\in M_i$ and let $g_i\in(R_i^\tau)^\times$ be determined by
\[
    \Bl_i(x_i,x_i)=\frac{g_i}{\Delta_i}\pmod\Lambda.
\]
The connected-sum decomposition is orthogonal for the Blanchfield pairings. 
The element $x=(x_1,x_2)$ generates $M(K)\cong R$, and
\begin{equation}\label{eq:twist}
    \Bl(x,x)=\frac{g}{\Delta}\pmod\Lambda,\qquad
    g=\Delta_2 g_1+\Delta_1 g_2\pmod\Delta,
\end{equation}
Under $R\cong R_1\times R_2$, its image is
\[
    g\longmapsto(\Delta_2g_1,\Delta_1g_2).
\]
The image of~$\Delta_2$ is an invariant unit in~$R_1$, and the image of~$\Delta_1$ is an invariant unit in~$R_2$.
Thus $g\in(R^\tau)^\times$, in agreement with nonsingularity of the pairing.

The following result shows a criterion for checking the algebraic
unknotting number of the connected sum.
\begin{proposition}\label{prop:cyclic}
  With the setup above, let $S\in\{\Z,\Q,\R\}$ and $\epsilon\in\{\pm1\}$.
  Then $u_S^\epsilon(K)=1$ if and only if
  \[
      [\epsilon g_1]=[\Delta_2]\text{ in }\cU(R_{1,S}),\qquad
      [\epsilon g_2]=[\Delta_1]\text{ in }\cU(R_{2,S}).
  \]
\end{proposition}
\begin{proof}
  By Definition~\ref{def:signed_S}, $u_S^\epsilon(K)=1$ precisely when $[\epsilon g]=0$ in~$\cU(R_S)$.
  Formula~\eqref{eq:twist} and the decomposition~\eqref{eq:additivity_of_U} identify this condition with
  \[
      [\epsilon\Delta_2g_1]=0\text{ in }\cU(R_{1,S}),\qquad
      [\epsilon\Delta_1g_2]=0\text{ in }\cU(R_{2,S}).
  \]
  We write these quotient groups additively, and each has exponent two, since an invariant unit~$a$ satisfies $a^2=\Nm(a)$.
  Thus these vanishing conditions are equivalent to the displayed equalities.
\end{proof}
\begin{remark}
  We point out that Proposition~\ref{prop:cyclic} does not simply state that $(g_1,g_2)$ is a norm in $R_1\times R_2$ if $g_1$ and $g_2$ are. In particular, the criteria of Proposition~\ref{prop:cyclic} can be satisfied even if $g_1$ and $g_2$ are not a norm themselves.
\end{remark}
The real and rational tests of Subsection~\ref{prop:qunkn} apply to  $\epsilon\Delta_2g_1$ and $\epsilon\Delta_1g_2$.
Once the rational tests pass, their classes lie in $\Ob(R_1)$ and $\Ob(R_2)$; the integral criterion asks whether both classes vanish.

\begin{remark}\label{rem:drops}
  Proposition~\ref{prop:cyclic} does not require $u_a(K_1)=u_a(K_2)=1$.
  The factors $\Delta_2$ and $\Delta_1$ can change the norm classes, allowing $u_a(K_1\#K_2)=1$ even when one summand has larger algebraic unknotting number.
\end{remark}
\begin{example}\label{ex:connected-drop}
  Let $K_1=T(2,5)$ and $K_2=-T(2,3)$, where the minus sign denotes the mirror with reversed orientation.
  Put $u=t+t^{-1}$. Then
  \[
      \Delta_1=u^2-u-1,\qquad\Delta_2=u-1,\qquad u\Delta_2-\Delta_1=1.
  \]
  Both Alexander polynomials are cyclotomic.
  Theorem~\ref{thm:Ob-group-CM-fields} allows
  us to compute the obstruction groups:
  $\Ob(R_1)=\Ob(R_2)=0$, $\cU(R_1)\cong(\Z/2)^2$ and $\cU(R_2)\cong\Z/2$. Moreover,
  the generators of $\cU(R_1)$ and $\cU(R_2)$ are detected at the infinite places, that is, at the roots of the Alexander polynomials
  on the unit circle; compare Table~\ref{tab:ob-groups-cyclotomic-polynomials}.

  Proposition~\ref{prop:product-norm-groups} allows us to
  compute that $\cU(R)\cong(\Z/2)^3$ and $\Ob(R)=0$,
  and $\cU(R)$ is again detected by the infinite
  places. That is, for $\epsilon=\pm1$, we have $u^\epsilon_\Z(K_1\#K_2)=1$
  if and only if $u^\epsilon_\R(K_1\#K_2)=1$.
  The latter is determined by the signatures
  as in Proposition~\ref{prop:real_number}. As
  the signature function of $K_1\#K_2$ takes values
  only in $\{-2,-1,0\}$, we conclude that $u_a(K_1\#K_2)=1$.

  In contrast, $\sigma_{K_1}(-1)=-4$ gives $u_a(K_1)\ge2$, and two crossing changes unknot~$T(2,5)$, so $u_a(K_1)=2$.
\end{example}
\begin{remark}
  The example $T(2,5)\#-T(2,3)$ can also be handled using Proposition~\ref{prop:cyclic}. However, for these two knots all the obstructions are determined by the signature function.
\end{remark}

\section{Examples}\label{sec:examples}

Throughout this section, unless otherwise specified, we consider a symmetric polynomial $\Delta$ and
define \(R=\Lambda/(\Delta)\).
We write \(R_0=R^\tau\) for the fixed subring. 

\subsection{Examples for simple knots}\label{sub:small}
The purpose of this subsection is to provide brief calculations
for simple knots. These examples illustrate various
techniques used later on.

In the trefoil, figure-eight and $6_1$ examples below, let $q\in(R^\tau)^\times$ be the invariant unit representing the Blanchfield pairing with respect to a chosen Alexander-module generator:
\[
  \Bl(a,b)=\frac{a q\tau(b)}{\Delta}\pmod\Lambda.
\]
Thus $q$ plays the role of $g$ in Lemma~\ref{lem:ub}; for a sign $\epsilon=\pm1$, the norm criterion is $[\epsilon q]=0$ in $\cU(R)$.
The sign $\epsilon$ is the presentation sign of Definition~\ref{def:signed_S}; the corresponding crossing change has sign $-\epsilon$, by Theorem~\ref{thm:fogel}.
Note that $u_a(K)=1$ 
is true if and only if $u_\Z^{+}(K)=1$ or $u_\Z^{-}(K)=1$.
Thus an obstruction to $u_a(K)=1$ must rule out both signs.

\begin{example}[The trefoil]\label{ex:trefoil}
    $\Delta=t-1+t^{-1}$ has ordinary representative $P=t\Delta=t^2-t+1$, which is monic with roots $e^{\pm i\pi/3}$ on the unit circle.
    Here $R = \Z[\xi_{6}]$, $E = \Q(\xi_{6})$ is CM with $F = \Q$, and $\cO_{E}^{\times}=\mu_{6}$, so $\Nm(R^{\times})=\{1\}$ and
    \[\cU(R) = \{ \pm 1 \} \cong \Z/2, \quad \text{and} \quad \Ob(R) = 0.\]
    This agrees with Corollary~\ref{cor:Ob-groups-cyclotomic-polynomials}: $Q=1$ and $[F:\Q]=1$.

    As $\cU(R)\cong\Z/2$, for any $q$, exactly one sign choice $\epsilon$ yields $[\epsilon q]=0\in\cU(R)$. This
    applies for $q=\wt{\Bl}(1,1)$, and $\epsilon$ is determined by the positivity of the trefoil via the signature jump
    at $e^{i\pi/3}$.
\end{example}

\begin{example}[The figure-eight knot]\label{ex:figure-eight}
  $\Delta=3-t-t^{-1}$ with real roots $(3\pm\sqrt5)/2$.
  Here $E=\Q(\sqrt5)$, the involution is the nontrivial element of $\Gal(E/\Q)$, $F=\Q$, and $R = \Z[(3+\sqrt5)/2] = \Z[\varphi] = \cO_{E}$ where $\varphi=(1+\sqrt5)/2$.
  Thus $R_{0}=\Z$, $R_{0}^{\times}=\{\pm1\}$ and $\Nm(\varphi)=-1$, so
  \[\cU(R)=0 .\]
  This agrees with computations in Table~\ref{tab:real-quad-fields}.

  For every invariant unit $q$ representing the Blanchfield pairing, both $[q]=0$ and $[-q]=0$ in $\cU(R)$.
  Explicitly, if $q=\Nm(v)$, then $-q=\Nm(\varphi v)$ because $\Nm(\varphi)=-1$.
  Hence both signs satisfy the criterion of Lemma~\ref{lem:ub}, giving $u_{\Z}^{+}(4_1)=u_{\Z}^{-}(4_1)=1$ and $u_a(4_1)=1$.
  The real place of $F=\Q$ splits in $E$, so there is no signature condition, in accordance with $\sigma(4_{1})=0$.
\end{example}

\begin{example}[The knot $6_{1}$]\label{ex:knot-6-1}
    Here $\Delta=-2t+5-2t^{-1}\doteq(t-2)(2t-1)$ and $\Delta(-1)=9=\det(6_{1})$. Its Alexander polynomial is $\Delta_{-2}$ of Example~\ref{ex:twist-knots}.
    The computation parallels Example~\ref{ex:quad-conductor-16}, with a few changes. The polynomial $\Delta$ is
    reducible, hence $\cA:=\Q[t^{\pm1}]/(\Delta)$ is not a field $E$: it is the algebra $\cA=\Q\times\Q$. The involution $\tau$ interchanges the factors.
    Every invariant unit $(c,c)$ of~$\cA$ is the norm of $(c,1)$, so $\cU(\cA)=0$. Hence there are no local conditions, and $\Ob(R)=\cU(R)$.

    Write $t_{\cA}$ for the image of $t$ in $\cA$. With the identification of $\cA$ with $\Q\times \Q$, we have $t_{\cA}=(2,1/2)$.
    Put $B=\Z[1/2]$. As $2(3-t_{\cA}-t_{\cA}^{-1})=1$, the diagonal copy of~$B$ lies in~$R$.
    Since $2\equiv\frac12\bmod3$, and $(3,0)=2t_{\cA}-1\in R$, we get
    \[
        R=\{(a,b)\in B\times B\colon a\equiv b\bmod 3\},\qquad \wt R=B\times B,
    \]
    where $B\times B=R[(1,0)]$ is integrally closed and the idempotent $(1,0)=(2t_{\cA}-1)/3$ is integral over~$R$.
    As in Example~\ref{ex:quad-conductor-16}, we argue that the conductor of~$R$ in~$\wt R$ is~$3\wt R$.    Since $B^\times=\{\pm2^k\}$ and $2\equiv-1\bmod3$, this congruence gives
    \[
        R^\times=\left\{(\epsilon_{1}2^{i},\epsilon_{2}2^{j})\colon\epsilon_{1}\epsilon_{2}=(-1)^{i+j}\right\},
        \qquad \Nm(a,b)=(ab,ab),
    \]
    so $\Nm(R^\times)=\{(-2)^{k}\}$, whereas $\Nm(\wt R^\times)=R_0^\times=\{\pm2^{k}\}$. Hence
    \[
        \Ob(R)=\cU(R)=\{\pm2^{k}\}/\{(-2)^{k}\}\cong\Z/2,\qquad \cU(\wt R)=0,
    \]
    and $\cU(R)$ is generated by $[-1]=[2]$. As in Example~\ref{ex:quad-conductor-16}, the whole obstruction group comes from the conductor.

    For every invariant unit $q$, exactly one of $[q]$ and $[-q]$ is zero, since they differ by the nonzero class $[-1]$. This implies
    that every knot with the Alexander polynomial $\Delta_{-2}$ and cyclic Alexander module has algebraic unknotting number $1$, but only one of the two $u^+_\Z$ and $u^-_\Z$ can be equal to $1$.
\end{example}

\begin{example}[Knots with Alexander polynomial of degree $2$]\label{ex:twist-knots}
  Let $a\neq 0$ be an integer. We consider knots with Alexander polynomial
  \[
    \Delta_a(t)=at-(2a-1)+at^{-1},\qquad a\in\Z\setminus\{0\}.
  \]
  The family includes twist knots.
  Also, $\Delta_1$ is the Alexander
  polynomial of the trefoil, $\Delta_{-1}$ of the figure
  eight knot, while $\Delta_{-2}$ is
  the Alexander polynomial of the $6_1$ knot.

  Assume that $D$ is square-free. Then Proposition~\ref{prop:im-quad-family} gives $R_a=\wt{R}_a$, so $R_a$ is integrally closed.
  Theorem~\ref{thm:norm-gps-im-quad-fields} yields
  \[
    \langle[-1]\rangle
      \longrightarrow\cU(R_a)
      \xrightarrow{\kappa}C_T/2C_T\longrightarrow0.
  \]
  If $|C_T|$ is odd, then $\kappa$ is a zero map, so $\cU(R)$ consists of the image of $[-1]$: if $[-1]=0$, then $\cU(R)$ is trivial,
  otherwise $\cU(R)$ is $\Z/2$. In any case, for any $q\in R_0^\times$ there exists a sign $\epsilon$ such that $[\epsilon q]=0\in\cU(R)$.
  That is, any knot $K$ with that Alexander polynomial, if it has Nakanishi index $1$, has algebraic unknotting number $1$: examples include $a=6$, compare Example~\ref{ex:quad-a6}.
  
  The correction from $C_T$ can be non-trivial. An example is $a=4$, discussed in Example~\ref{ex:74}. We have
  \[
    \cU(R_4)=\{[1],[-1],[2],[-2]\}\cong(\Z/2)^2,
    \qquad \Ob(R_4)=0.
  \]
  Neither $2$ nor $-2$ is a norm of a unit, so a pairing represented by the invariant unit $q=2$ fails the criterion for both signs.

  If $D$ is nonsquare but not square-free, the situation is more complex.
  For example, Example~\ref{ex:quad-conductor-16} shows
  \[
    \cU(R_{16})\cong(\Z/2)^2,\qquad
    \Ob(R_{16})=\langle[2]\rangle\cong\Z/2.
  \]
  It follows that if a knot $K$ has Alexander polynomial $\Delta_{16}$ and $\wt{\Bl}(1,1)=q$ with $q=\pm 2$, then the algebraic unknotting
  number is at least $2$.

  The parameters with square $D$ lie outside the quadratic-field hypotheses.
  They are precisely
  \[
    a=-k(k+1),\qquad k\geq1,
    \qquad
    P_a(t)=t\Delta_a(t)=-(kt-k-1)((k+1)t-k).
  \]
  The two linear factors are exchanged by the involution, so the total quotient algebra is $\Q\times\Q$ with the swapping involution.
  Its norm quotient is zero, and hence $\Ob(R_a)=\cU(R_a)$; this group may nevertheless be nonzero.
  The case $k=1$ is exactly Example~\ref{ex:knot-6-1}.
\end{example}

\subsection{Lickorish criterion vs. local conditions}\label{sub:lick_local}

We compare the conditions for the rational unknotting number
with the Lickorish criterion. We have shown in Theorem~\ref{thm:lick_rat} that the Lickorish criterion, together with the signature condition, implies that $\epsilon g$ is a norm over $\Q$. We provide examples that show the following.
\begin{itemize}
\item The Lickorish criterion does not imply real algebraic unknotting number one; see Example~\ref{ex:lickorish-gap}.
\item The Lickorish criterion is effectively stronger than the finite places obstruction for the rational unknotting number; see Example~\ref{ex:residue-degree}.
\item The Lickorish criterion does not imply that the algebraic unknotting number is $1$, that is, it is weaker than the $\Ob$ calculation;
  see Subsection~\ref{sec:exampl-where-class}.
\end{itemize}
We remark that the fact that the Lickorish criterion does not see
the signatures was already known \cite{Jabuka,BF1}.

The first example is
an obstruction invisible to the Lickorish criterion, but detected by local norms at infinite places.
\begin{example}\label{ex:lickorish-gap}
    Let $P=\Phi_{15}=t^8-t^7+t^5-t^4+t^3-t+1$ and $\Delta=t^{-4}P$.
    Then $\Delta(1)=\Delta(-1)=1$, so Lickorish's congruence is vacuous.
    Here
    \[
        \begin{gathered}
            R=\Z[\xi_{15}]=\cO_E,\qquad E=\Q(\xi_{15}),\qquad F=E^\tau=\Q(u_F),\\
            t_E=\xi_{15},\qquad u_F=t_E+t_E^{-1}.
        \end{gathered}
    \]
    Corollary~\ref{cor:Ob-groups-cyclotomic-polynomials} gives no finite ramification
    in $E/F$ and
    \[
        |\cU(R)|=2^{\varphi(15)/2-1}=8.
    \]
    Thus, by Theorem~\ref{thm:finite}, only the four real places of $F$ impose
    local norm conditions.
    For detailed computations, see Subsubsection~\ref{sec:cyclotomic-fields} and Table~\ref{tab:ob-groups-cyclotomic-polynomials}.

    An explicit invariant unit failing both signed norm tests is $q=u_F$.
    Indeed, the relation
    \[
        f(u_F)=0,\qquad f(X)=X^4-X^3-4X^2+4X+1,
        \qquad P(t)=t^4f(t+t^{-1}),
    \]
    gives $q^{-1}=-u_F^3+u_F^2+4u_F-4\in R_0$.
    At the four real embeddings indexed by $k=1,2,4,7$, the values of $q$ are
    $2\cos(2\pi k/15)$, with signs $(+,+,-,-)$.
    Hence neither $q$ nor $-q$ is totally positive, and neither is a norm from $E$.
    Both signs pass the vacuous Lickorish congruence and all finite local norm tests;
    their failure is detected at the real places.
\end{example}

\begin{example}[Lickorish is stronger than the finite local norm tests]\label{ex:residue-degree}
    Consider
    \[
        P(t)=t^4-11t^3+21t^2-11t+1,\qquad \Delta=t^{-2}P(t),
        \qquad R=\Z[t^{\pm1}]/(P),
    \]
    so that $\Delta(1)=1$ and $\Delta(-1)=45$.
    Then $P(t)=t^2f(t+t^{-1})$ with $f(X)=X^2-11X+19$, whose roots are $(11\pm3\sqrt5)/2$. Thus
    \[
        F=\Q(\sqrt5),\qquad u_F=4+3\phi,\qquad \phi=(1+\sqrt5)/2,\qquad E=F(\sqrt{u_F^2-4}).
    \]
    Since $\Nm_{F/\Q}(u_F^2-4)=f(2)f(-2)=45$ is not a square, $E/F$ is quadratic and $P$ is irreducible.
    Both real places of~$F$ split in~$E$. Also, $\phi$ is the fundamental unit of $\cO_F$.

    Since $P$ is monic with constant term~$1$, we have $T=\emptyset$, so the normalization of~$R$ is~$\cO_E$; moreover, $R=\Z[t_E]=\Z[u_F]\oplus t_E\Z[u_F]$, because $t_E^2=u_Ft_E-1$.
    As $\tau(x+yt_E)=(x+yu_F)-yt_E$ for $x,y\in\Z[u_F]$, the fixed ring is $R_0=\Z[u_F]=\Z[3\phi]$.
    From $(3\phi)^2=3\cdot3\phi+9$ we get
    \[
        R_0=\Z\oplus3\phi\Z=\Z+3\cO_F.
    \]

    The ideal $(3)$ is inert in~$F$, so $\cO_F/3\cO_F\cong\F_9$ and $R_0/3\cO_F\cong\F_3$.
    Hence a unit $v\in\cO_F^\times$ lies in~$R_0^\times$ if and only if its residue lies in~$\F_3^\times$.
    Now $\cO_F^\times=\{\pm\phi^m\}$ and $\phi^4=3\phi+2\equiv-1\bmod3$, so $\phi$ has order~$8$ in~$\F_9^\times$, 
    and the residue of $\phi^m$ lies in~$\F_3^\times$ if and only if $4\mid m$.
    With $\psi=\phi^4=u_F-2$, we obtain
    \[
        R_0^\times=\{\pm \psi^m\colon m\in\Z\}.
    \]

    Since $P(1)=1$, the element $1-t_E$ divides $P(t_E)-P(1)=-1$, so $1-t_E\in R^\times$, and $\Nm(1-t_E)=2-u_F=-\psi$.
    Hence $[\psi]=[-1]$, and $\cU(R)$ is generated by~$[-1]$.
    The map $R\to R/(t_E+1)\cong\Z/45$ sends every norm $\Nm(v)=v\tau(v)$ to the square $v(-1)^2$, and $-1$ is not a square modulo~$3$. Therefore
    \[
        \Nm(R^\times)=\langle-\psi\rangle,\qquad \cU(R)=\langle[-1]\rangle\cong\Z/2.
    \]
    On the other hand, $1-t_E$ and~$\phi$ are units of~$\cO_E$, and
    \[
        \Nm_{E/F}\left((1-t_E)/\phi^2\right)=-\psi/\phi^4=-1.
    \]
    So $-1$ is the norm of a unit of the normalization, all local invariants of~$[-1]$ are trivial, and
    \[
        \lambda_R\equiv\mathbf{1},\qquad \Ob(R)=\cU(R)\cong\Z/2.
    \]

    For the coefficient $q=-1$, the Lickorish criterion fails for $\epsilon=+1$, whereas for $\epsilon=-1$ we have $\epsilon q=1=\Nm(1)$.
    Hence a knot $K$ with cyclic Alexander module~$R$ and Blanchfield coefficient $q=-1$ satisfies both rational signed criteria, $u_\Q^+(K)=u_\Q^-(K)=1$,
    but $u_\Z^+(K)>1$ and $u_\Z^-(K)=1$.
    In particular, $u_a(K)=1$: the nontrivial class $[-1]\in\Ob(R)$ obstructs only the positive sign.
\end{example}

\begin{example}
    For the normalized Alexander polynomial \(\Delta_{16}(t)=16(t+t^{-1})-31\) from Example~\ref{ex:quad-conductor-16}, put \(R=\Lambda/(\Delta_{16})\). Then
    \[\Ob(R) = \langle [2] \rangle \cong \Z/2.\]
    Take $q=2$. Example~\ref{ex:quad-conductor-16} shows that $q$ is a field norm, so all its local norm tests pass.
    For $\epsilon=+1$, the Lickorish criterion fails because $2$ is not a square modulo~$3$.
    The prime $(3)$ does not ramify in $E_{16}=\Q(\sqrt{-7})$; this failure comes from the conductor.
    For $\epsilon=-1$, the test also fails: $-2\equiv5\pmod7$ is not a square.
    Thus neither $q$ nor $-q$ is a square modulo $|\Delta_{16}(-1)|=63$.
    A knot $K$ with cyclic Alexander module $R$ and coefficient $q=2$ therefore fails both signed integral criteria, giving $u_a(K)>1$.
\end{example}

\subsection{Obstructions undetected by local and Lickorish tests}
\label{sec:exampl-where-class}

The following examples show that the local norm conditions and the Lickorish criterion do not guarantee that $u_a=1$.

In Example~\ref{ex:only-Ob}, $q$ passes the Lickorish test for both choices of $\epsilon=\pm 1$, but its class in $\Ob(R)$ is nontrivial.  The class of $-q$ fails at the real places.

In Example~\ref{ex:all-silent}, both $q$ and $-q$ pass the Lickorish test and represent the same nontrivial class in $\Ob(R)$, so $u_a(K)>1$.

\begin{example}\label{ex:only-Ob}
    Let
    \[
        P(t)=5t^4-4t^3-t^2-4t+5,\qquad \Delta=t^{-2}P(t),\qquad
        R=\Lambda/(\Delta).
    \]
    Thus $\Delta(1)=1$ and $\Delta(-1)=17$.
    Let $t_E$ be the image of $t$ in $E=\Q\otimes_\Z R$, and put $u_F=t_E+t_E^{-1}$ and $q=668-345u_F$.
    Example~\ref{ex:arithmetic-only-ob} proves that $R$ is integrally closed,
    $\Ob(R)=\{[1],[q]\}\cong\Z/2$, and $q$ is a field norm, whereas $-q$ fails the real norm conditions.
    Neither $q$ nor $-q$ is the norm of a unit in $R$.

    For a knot $K$ with cyclic Alexander module $R$ and Blanchfield coefficient $q$,
    these facts give $u_\Q^+(K)=1$ and $u_\Q^-(K)>1$, while $u_\Z^+(K)>1$ and $u_\Z^-(K)>1$, by Definition~\ref{def:signed_S} and Lemma~\ref{lem:ub}.
    Hence $u_a(K)>1$.
    Nevertheless, the Lickorish congruence of Proposition~\ref{prop:square} holds for both signs.
    Under $R\to\Z/17$ given by $t_E\mapsto-1$, we have $u_F\mapsto-2$ and $q\mapsto1358$, so
    \[
        1358\equiv15\equiv7^2\pmod{17},\qquad
        -1358\equiv2\equiv6^2\pmod{17}.
    \]
    Thus the positive presentation sign passes both the rational norm test and the Lickorish criterion,
    but the nontrivial class $[q]\in\Ob(R)$ obstructs an integral norm.
\end{example}

\begin{example}\label{ex:all-silent}
    Let
    \[
        P(t)=-3t^4-15t^3+37t^2-15t-3,\qquad \Delta=t^{-2}P(t),\qquad
        R=\Lambda/(\Delta).
    \]
    Thus $\Delta(1)=1$ and $\Delta(-1)=61$.
    Let $t_E$ be the image of $t$ in $E=\Q\otimes_\Z R$, and put $u_F=t_E+t_E^{-1}$ and $q=6-3u_F$.
    Example~\ref{ex:arithmetic-all-silent} proves that $R$ is integrally closed and
    \[
        \cU(R)=\Ob(R)=\{[1],[q]\}\cong\Z/2,\qquad [-q]=[q].
    \]
    Both $q$ and $-q$ are field norms, but neither is the norm of an $R$-unit.
    Consequently, a knot $K$ with cyclic Alexander module $R$ and Blanchfield coefficient $q$
    satisfies $u_\Q^+(K)=u_\Q^-(K)=1$ but $u_\Z^+(K)>1$ and $u_\Z^-(K)>1$, so $u_a(K)>1$.

    The Lickorish congruence holds for both signs as well.
    Under $R\to\Z/61$ given by $t_E\mapsto-1$, we have $u_F\mapsto-2$ and $q\mapsto12$, and
    \[
        12\equiv16^2\pmod{61},\qquad
        -12\equiv49\equiv7^2\pmod{61}.
    \]
    Thus both signs pass the local norm conditions and the Lickorish criterion,
    while their common nontrivial class in $\Ob(R)$ obstructs an integral norm.
\end{example}

\appendix
\section{PARI/GP code}\label{app:pari}
This appendix lists the PARI/GP code~\cite{PARI2} behind the computations in the paper. 

\subsection{Cyclotomic rings}\label{app:cyclotomic}
The first loop prints the rows of Table~\ref{tab:ob-groups-cyclotomic-polynomials}.
The second loop returns the smallest~$n$ with two distinct prime divisors and $\Ob(\Z[\xi_n])\ne0$, that is $n=39$.
\begin{lstlisting}[basicstyle=\ttfamily,frame=none]
{
foreach([6, 10, 12, 14, 15, 20, 21, 30, 39, 56, 58], n,
  F = bnfinit(polredbest(minpoly(Mod(x + 1/x, polcyclo(n, x)))), 1);
  m = if (n % 4 == 2, n/2, n);  Q = if (isprimepower(m), 1, 2);
  h = F.no;  hp = bnfnarrow(F)[1];  d = valuation(hp/h, 2);
  U = concat([-1], F.fu);
  s = F.r1 - matrank(matrix(F.r1, #U, i, j, nfeltsign(F, U[j])[i] < 0) * Mod(1, 2));
  print([n, eulerphi(n), F.r1, m, polcyclo(n, -1), Q, h, hp, d,
         F.r1 - (Q == 2), d - (Q == 2), bnfcertify(F)], "  ", s));

for (n = 6, 39,
  if (omega(n) < 2, next);
  F = bnfinit(polredbest(minpoly(Mod(x + 1/x, polcyclo(n, x)))), 1);
  Q = if (isprimepower(if (n % 4 == 2, n/2, n)), 1, 2);
  if (valuation(bnfnarrow(F)[1] / F.no, 2) > (Q == 2), print(n); break));
}
\end{lstlisting}

\subsection{Quadratic fields}\label{app:quadratic}
For the rows of Tables~\ref{tab:im-quad-fields} and~\ref{tab:real-quad-fields}, we use the following code. Here,
the last two lines give the class group and the class of the prime above~$109$ in Example~\ref{ex:im-quad-fields-nontrivial-Ob}.
\begin{lstlisting}[basicstyle=\ttfamily, frame=none]
{
A = [1, 2, 3, 4, 5, 6, 8, 109^2, -1, -5, -55, -4861];
for (k = 1, #A + 1000,
  a = if (k <= #A, A[k], #A - k);  D = 1 - 4*a;
  if (!issquarefree(D), next);
  E = bnfinit(x^2 - x + a, 1);  B = concat([-1], factor(abs(a))[, 1]~);
  G = concat([E.tu[2]], E.fu);
  if (#B > 1,
    T = [];  for (i = 2, #B, T = concat(T, idealprimedec(E, B[i])));
    S = bnfsunit(E, T)[1];  G = concat(G, vector(#S, j, nfbasistoalg(E, S[j]))));
  N = vector(#G, j, norm(Mod(lift(G[j]), E.pol)));
  V = matrix(#B, #G, i, j, if (i == 1, N[j] < 0, valuation(N[j], B[i]) % 2)) * Mod(1, 2);
  P = factor(abs(D))[, 1]~;  e = (D < 0);
  L = matrix(e + #P, #B, i, j, if (e && i == 1, B[j] < 0, kronecker(B[j], P[i - e]) < 0));
  r = matrank(V);  dU = #B - r;
  dOb = matrank(matconcat([V, matker(L * Mod(1, 2))])) - r;
  eE = e || norm(Mod(lift(E.fu[1]), E.pol)) == 1;
  eT = matrank(matconcat([V, vectorv(#B, i, i == 1) * Mod(1, 2)])) > r;
  if (k <= #A, print([a, abs(D), E.no, if (eE == eT && dU == eE, "trivial", "nontrivial"),
                      dU, dOb, bnfcertify(E)]));
  if (k > #A && dOb > 0, print(a, "  ", factor(D)); break));
E = bnfinit(x^2 - x + 109^2, 1);
print([E.cyc, bnfisprincipal(E, idealprimedec(E, 109)[1], 0), bnfcertify(E)]);
}
\end{lstlisting}

\subsection{The quartic examples}\label{app:quartic}
The two blocks below check 
the quartic examples from Section~\ref{sub:quartic-obstructions}.
Here \texttt{t} and~\texttt{u} stand for $t_E$ and~$u_F$, and \texttt{W} is the $\Z$-basis of~$\cO_E$ used in the example; its determinant is~$\pm1$.

The block checking Example~\ref{ex:arithmetic-only-ob} is as follows.
\begin{lstlisting}[basicstyle=\ttfamily, frame=none]
{
H = x^4 - 4*x^3 - 5*x^2 - 100*x + 625;
E = bnfinit(H, 1);  t = Mod(x/5, H);  u = t + 1/t;  a = 5*u - 2;
F = bnfinit(y^2 - 59, 1);  f = 5*X^2 - 4*X - 11;
print([lift(subst(f, X, u)), lift(a^2 - 59), subst(f, X, 2)*subst(f, X, -2)/25]);
print(polrootsreal(f)~);
W = [1, 5*t, 5*t^2 - 4*t, 5*t^3 - 4*t^2 - t];
print([matdet(Mat(vector(4, i, nfalgtobasis(E, lift(W[i]))))), lift(5/t - (4 - W[4]))]);
print(factormod(5*y^4 - 4*y^3 - y^2 - 4*y + 5, 5));
foreach(idealprimedec(E, 5), P, print([P.e, P.f, nfeltval(E, lift(t), P)]));
r = 5*u + 6;  q = 530 - 69*a;
print([norm(Mod(8 + y, F.pol)), lift(r*(2 - u) - 1), lift(q - (668 - 345*u))]);
for (k = 1, 69, if (issquare(59*k^2 + 1) || issquare(59*k^2 - 1), print(k); break));
print([F.fu, F.cyc, bnfcertify(F)]);
w = 5*t^3 + t^2 - 5;  z = ((561 - 73*a)*t + 114 - 15*a)/31;
V = [w, z, 1 + t];  N = vector(3, i, V[i]*Mod(subst(lift(V[i]), x, 25/x), H));
print(lift([N[1] - r, N[2] - q, N[3] - (2 + u), 2 + u - (31 + 4*a)/r]));
print([E.cyc, bnfcertify(E)]);
foreach(idealprimedec(E, 5), P, if (P.f == 1, print([nfeltval(E, lift(t), P),
  nfeltval(E, lift(w), P), nfeltval(E, lift(r/w), P), bnfisprincipal(E, P, 0)])));
D = rnfdisc(F, x^2 - (y + 2)*x + 25)[1];
print([nfdisc(H) == F.disc^2 * 17, idealhnf(F, D) == idealhnf(F, 31 + 4*y)]);
print([idealhnf(E, lift((1 + t)*r/w)) == idealhnf(E, P) | P <- idealprimedec(E, 17), P.e == 2]);
}
\end{lstlisting}

Example~\ref{ex:arithmetic-all-silent} is handled with the following block.
\begin{lstlisting}[basicstyle=\ttfamily]
{
H = x^4 + 15*x^3 - 111*x^2 + 135*x + 81;
E = bnfinit(H, 1);  t = Mod(x/3, H);  u = t + 1/t;  a = 3*u + 8;
F = bnfinit(y^2 - y - 185, 1);  f = 3*X^2 + 15*X - 43;
print([lift(subst(f, X, u)), lift(a^2 - a - 185), subst(f, X, 2)*subst(f, X, -2)/9]);
print(polrootsreal(f)~);
W = [1, -3*t, -3*t^2 - 15*t, -3*t^3 - 15*t^2 + 37*t];
print([matdet(Mat(vector(4, i, nfalgtobasis(E, lift(W[i]))))), lift(3/t - (W[4] - 15)),
       lift(3*(u^2 + 5*u - 14))]);
for (k = 1, 9, if (issquare(741*k^2 + 4) || issquare(741*k^2 - 4), print(k); break));
print([F.fu, factormod(y^2 - y - 185, 3)]);
q = 14 - a;  ep = 118 + 9*a;
print([norm(Mod(14 - y, F.pol)), lift(q*(-u - 7) - 1)]);
b = -3*t^3 - 18*t^2 + 19*t + 4;  w = -15*t^3 - 90*t^2 + 98*t + 40;
z = ((149 + 14*a)*t - 335 - 6*a)/13;  g = (15 + a) - (10 + a)*t;
V = [b, w, z, g, 1 + t];  N = vector(5, i, V[i]*Mod(subst(lift(V[i]), x, 9/x), H));
print(lift([N[1] - q*ep, N[2] + 1, N[3] - q, N[4] - (a - 15), N[5] - (2 + u),
            2 + u + q*(53 + 4*a)/3, g - (-3*t^3 - 18*t^2 + 22*t + 5)]));
print([E.cyc, bnfcertify(E), F.cyc, bnfcertify(F),
       bnfisprincipal(F, idealhnf(F, 5, y), 0)]);
foreach(idealprimedec(E, 3), P, print([nfeltval(E, lift(t), P), nfeltval(E, lift(b), P),
  nfeltval(E, lift(q/b), P), bnfisprincipal(E, P, 0)]));
P5 = [P | P <- idealprimedec(E, 5), nfeltval(E, lift(a), P) > 0][1];
print([P5.e, P5.f, idealhnf(E, lift(g)) == idealhnf(E, P5)]);
D = rnfdisc(F, x^2 - (y - 8)*x + 9)[1];
print([nfdisc(H) == F.disc^2 * 61, idealhnf(F, D) == idealhnf(F, 53 + 4*y)]);
print([idealhnf(E, lift((1 + t)*b)) == idealhnf(E, P) | P <- idealprimedec(E, 61), P.e == 2]);
}
\end{lstlisting}

\def\MR#1{}
\bibliographystyle{amsalpha}
\bibliography{references.bib}
\end{document}